\documentclass[11pt]{amsart}
\usepackage{a4wide}
\usepackage{amsmath}
\usepackage[utf8]{inputenc}
\usepackage{amssymb}
\usepackage{amsopn}
\usepackage{epsfig}
\usepackage{amsfonts}
\usepackage{latexsym}
\usepackage{graphicx}
\usepackage{enumerate}
\usepackage[dvipsnames]{xcolor}
\usepackage{ mathrsfs }
\usepackage{tikz}
\usetikzlibrary{arrows,automata}

\newtheorem{theorem}{Theorem}[section]
\newtheorem{lemma}[theorem]{Lemma}
\newtheorem{proposition}[theorem]{Proposition}

\newtheorem{main}{Theorem}

\theoremstyle{definition}
\newtheorem{definition}[theorem]{Definition}
\newtheorem{example}[theorem]{Example}

\theoremstyle{remark}
\newtheorem{remark}[theorem]{Remark}

\numberwithin{equation}{section}

\newcommand{\N}{\ensuremath{\mathbb{N}}}
\renewcommand{\L}{\ensuremath{\mathbb{L}}}
\renewcommand{\S}{\ensuremath{\mathbb{S}}}

\newcommand{\F}{\Omega_{F}^*}

\newcommand{\LL}{\Omega_{L}^*}

\renewcommand{\c}{ {\mathbf{c}}}
\renewcommand{\a}{ {\mathbf{a}}}
\renewcommand{\d}{ {\mathbf{d}}}
\renewcommand{\b}{ {\mathbf{b}}}

\renewcommand{\t}{ {\mathbf{t}}}
\newcommand{\s}{ {\mathbf{s}}}
\newcommand{\cs}{ {\mathbf{S}}}
\newcommand{\ca}{ {\mathbf{A}}}
\renewcommand{\r}{ {\mathbf{r}}}

\newcommand{\set}[1]{\left\{#1\right\}}

\newcommand{\Ga}{\Gamma}
\newcommand{\ep}{\varepsilon}
\newcommand{\f}{\infty}
\newcommand{\de}{\delta}

\newcommand{\lle}{\preccurlyeq}
\newcommand{\lge}{\succcurlyeq}
\renewcommand{\a}{ \mathbf{a}}
\newcommand{\si}{\sigma}
\newcommand{\La}{\Lambda}

\newcommand{\K}{\mathbf K}

\begin{document}
\title[Critical values {for the} $\beta$-transformation]{Critical values for the $\beta$-transformation with a hole at $0$}

\author[P. Allaart]{Pieter Allaart}
\address[P. Allaart]{Mathematics Department, University of North Texas, 1155 Union Cir \#311430, Denton, TX 76203-5017, U.S.A.}
\email{allaart@unt.edu}

\author[D. Kong]{Derong Kong}
\address[D. Kong]{College of Mathematics and Statistics, Chongqing University, Chongqing 401331, People's Republic of China.}
\email{derongkong@126.com}

\dedicatory{}


\subjclass[2020]{Primary: 37B10, 28A78; Secondary: 68R15, 26A30,37E05}

\begin{abstract}
Given $\beta\in(1,2]$, let $T_\beta$ be the $\beta$-transformation on the unit circle $[0,1)$ such that $T_\beta(x)=\beta x\pmod 1$. For each $t\in[0,1)$ let $K_\beta(t)$ be the survivor set consisting of all $x\in[0,1)$ whose orbit $\{T^n_\beta(x): n\ge 0\}$ never hits the open interval $(0,t)$.  Kalle et al.~proved in [{\em Ergodic Theory Dynam. Systems}, 40 (9): 2482--2514, 2020] that the Hausdorff dimension function $t\mapsto\dim_H K_\beta(t)$ is a non-increasing Devil's staircase. So there exists a critical value $\tau(\beta)$ such that $\dim_H K_\beta(t)>0$ if and only if $t<\tau(\beta)$. In this paper we determine the critical value $\tau(\beta)$ for all $\beta\in(1,2]$, answering a question of Kalle et al.~(2020). For example, we find that for the {\em Komornik-Loreti constant} $\beta\approx 1.78723$ we have $\tau(\beta)=(2-\beta)/(\beta-1)$.
Furthermore, we show that (i) the function $\tau: \beta\mapsto\tau(\beta)$ is left continuous on $(1,2]$ with right-hand limits everywhere, but has countably infinitely many discontinuities; (ii) $\tau$ has no downward jumps, with $\tau(1+)=0$ and $\tau(2)=1/2$; and (iii) there exists an open set $O\subset(1,2]$, whose complement $(1,2]\setminus O$ has zero Hausdorff dimension, such that $\tau$ is real-analytic, convex and strictly decreasing on each connected component of $O$.  Our strategy to find the critical value $\tau(\beta)$ depends on certain substitutions of Farey words and a renormalization scheme from dynamical systems.
\end{abstract}

\keywords{$\beta$-transformation, survivor set, Farey word, Lyndon word, Farey interval, substitution operator, Hausdorff dimension}

\maketitle
\tableofcontents

\section{Introduction} \label{sec: Introduction}

The mathematical study of dynamical systems with holes, called \emph{open dynamical systems}, was first proposed by Pianigiani and Yorke \cite{Pianigiani-Yorke-1979} in 1979. In recent years open dynamical systems have received considerable attention from both theoretical and applied perspectives (cf.~\cite{Demers-Wright-Young-2010, Demmers-2005, Demers-Young-2006}). In the general setting, one considers a discrete dynamical system $(X, T)$, where $X$ is a compact metric space and $T: X\to X$ is a continuous map having positive topological entropy. Let $H\subset X$ be an open connected set, called the \emph{hole}. It is interesting to study the set of points $x\in X$ whose orbit $\set{T^n(x): n\ge 0}$ never hits the hole $H$. In other words, we are interested in the \emph{survivor set}
\[
K(H)=\set{x\in X: T^n(x)\notin H~\forall n\ge 0}=X\setminus\bigcup_{n=0}^\f T^{-n}(H).
\]
It is known that the size of ${K}(H)$ depends not only on the size but also on the position of the hole $H$ (cf.~\cite{Bunimovich-Yurchenko-2011}).
 In \cite{Urbanski_1986, Urbanski-87} Urba\'nski considered $C^2$-expanding, orientation-preserving circle maps with a hole of the form $(0,t)$. In particular,  he proved that for the doubling map $T_2$ on the circle $\mathbb R/\mathbb Z\sim[0,1)$, i.e., $T_2: [0,1)\to[0,1);\quad x\mapsto 2x\pmod 1$, the Hausdorff dimension of the survivor set
$
 K_2(t):=\set{x\in[0,1): T_2^n(x)\notin(0,t)~\forall n\ge 0}
 $
 depends continuously on the parameter $t\in[0,1)$. Furthermore, he showed that  the dimension function $\eta_2: t\mapsto \dim_H K_2(t)$ is a devil's staircase, and studied its bifurcation set.  Carminati and Tiozzo \cite{Carminati-Tiozzo-2017} showed that the function $\eta_2$ has an interesting analytic property: the local H\"older exponent of $\eta_2$ at any bifurcation point $t$ is equal to $\eta_2(t)$. For the doubling map $T_2$ with an arbitrary hole $(a,b)\subset[0,1)$, Glendinning and Sidorov \cite{Glendinning-Sidorov-2015} studied (i) when the survivor set $K_2(a,b)=\set{x\in[0,1): T_2^n (x)\notin(a,b)~\forall n\ge 0}$ is nonempty; (ii) when $K_2(a,b)$ is infinite; and (iii) when $K_2(a,b)$ has positive Hausdorff dimension. They proved that when the size of the hole $(a,b)$ is strictly smaller than  $0.175092$, the survivor set $K_2(a,b)$ has positive Hausdorff dimension. The work of Glendinning and Sidorov was partially extended by Clark \cite{Lyndsey-2016} to {the} $\beta$-dynamical system $([0,1), T_\beta)$ with a hole $(a,b)$, where $\beta\in(1,2]$ and $T_\beta(x):=\beta x\pmod 1$.

Motivated by the above works, Kalle et al.~\cite{Kalle-Kong-Langeveld-Li-18} considered the survivor set in the $\beta$-dynamical system $([0, 1), T_\beta)$ with a hole at zero.  More precisely, for $t\in[0,1)$ they determined the Hausdorff dimension of the survivor set
\[
K_\beta(t)=\set{x\in[0,1): T_\beta^n (x)\notin (0,t)~\forall n\ge 0},
\]
and showed that the dimension function $\eta_\beta: t\mapsto \dim_H K_\beta(t)$ is a non-increasing Devil's staircase. So there exists a critical value $\tau(\beta)\in[0,1)$ such that
$\dim_H K_\beta(t)>0$ if and only if $t<\tau(\beta)$.
Kalle et al.~\cite{Kalle-Kong-Langeveld-Li-18} gave general lower and upper bounds for $\tau(\beta)$. In particular, they showed that $\tau(\beta)\le 1-\frac{1}{\beta}$ for all $\beta\in(1,2]$, and the equality $\tau(\beta)=1-\frac{1}{\beta}$ holds for infinitely many $\beta\in(1,2]$. They left open the interesting question to determine $\tau(\beta)$ for all $\beta\in(1,2]$. In this paper we give a complete description of the critical value
 \begin{equation} \label{eq:defi-critical-value}
 \tau(\beta)=\sup\set{t: \dim_H K_\beta(t)>0}=\inf\set{t: \dim_H K_\beta(t)=0}
 \end{equation}
 for each $\beta\in(1,2]$. Qualitatively, our main result is:

\begin{main} \label{main:critical-devils-staircase}\mbox{}
\begin{enumerate}[{\rm(i)}]
\item The function $\tau: \beta\mapsto\tau(\beta)$ is left continuous on $(1,2]$ with right-hand limits everywhere {\em(c\`adl\`ag)}, and as a result has only countably many discontinuities;
\item $\tau$ has no downward jumps;
\item There is an open set $O\subset (1,2]$, whose complement $(1,2]\setminus O$ has zero Hausdorff dimension, such that $\tau$ is real-analytic, convex and strictly decreasing on each connected component of $O$.
\end{enumerate}
\end{main}

Quantitatively, the main results are Theorem \ref{main:critical-basic-intervals} and Propositions \ref{prop:cric-exception-1} and \ref{prop:cri-exception-2}. {Together with Proposition \ref{lem:critical-non-Fareyinterval}, they specify the value of $\tau(\beta)$ for all $\beta\in(1,2]$.}
In Proposition \ref{prop:discontinuities} below we give an explicit description of the discontinuities of the map $\tau$, which shows that the dimension $\dim_H K_\beta(t)$ is not jointly continuous in $\beta$ and $t$.
The closures of the connected components of the set $O$ in Theorem \ref{main:critical-devils-staircase} (iii) form a pairwise disjoint collection $\{I_\alpha\}$ of closed intervals which we call {\em basic intervals} (see Definition \ref{def:basic-intervals}). In the remainder of this introduction we describe these basic intervals by using certain substitutions on Farey words. We then give a formula for $\tau(\beta)$ on each basic interval (see Theorem \ref{main:critical-basic-intervals}) and decompose the complement $(1,2]\backslash \bigcup_\alpha I_\alpha$ into countably many disjoint subsets (see Theorem \ref{main:geometrical-structure-basic-intervals}), which are of two essentially different types. We then calculate $\tau(\beta)$ on each subset.

\begin{figure}[h!]
 \begin{center}
\begin{tikzpicture}[xscale=13,yscale=13]
\draw [->] (1,0) node[anchor=north] {$1$}  -- (2.05,0) node[anchor=west] {$\beta$};
\draw [->] (1,0) node[anchor=east] {$0$} -- (1,0.55) node[anchor=south] {$\tau(\beta)$};

\draw[dashed] (2,-0.005) node[anchor=north] {$2$} -- (2,0.5)--(1,0.5)node[anchor=east] {$\frac{1}{2}$};

\draw[variable=\q,domain=1:2,dashed,red] plot({\q},{(1-1/\q)})(1.5,0.35)node[above,scale=0.8pt](1.5,){$1-\frac{1}{\beta}$};

\draw[densely dotted,variable=\q,domain=1:1.1187] plot({\q},{1-1/\q});
\draw[thick,variable=\q,domain=1.1187: 1.12241] plot({\q},({1/\q/(pow(\q,20)-1)});	

\draw[thick,variable=\q,domain=1.12331: 1.12734] plot({\q},({1/\q/(pow(\q,19)-1)});

\draw[thick,variable=\q,domain=1.12836: 1.13275] plot({\q},({1/\q/(pow(\q,18)-1)});

\draw[thick,variable=\q,domain=1.1339: 1.13871] plot({\q},({1/\q/(pow(\q,17)-1)});

\draw[thick,variable=\q,domain=1.14003: 1.14532] plot({\q},({1/\q/(pow(\q,16)-1)});

\draw[thick,variable=\q,domain=1.14685: 1.15271] plot({\q},({1/\q/(pow(\q,15)-1)});

\draw[thick,variable=\q,domain=1.15449: 1.16102] plot({\q},({1/\q/(pow(\q,14)-1)});

\draw[thick,variable=\q,domain=1.16312: 1.17045] plot({\q},({1/\q/(pow(\q,13)-1)});

\draw[thick,variable=\q,domain=1.17295: 1.18126] plot({\q},({1/\q/(pow(\q,12)-1)});

\draw[thick,variable=\q,domain=1.18428: 1.19379] plot({\q},({1/\q/(pow(\q,11)-1)});

\draw[thick,variable=\q,domain=1.19749: 1.20851] plot({\q},({1/\q/(pow(\q,10)-1)});

\draw[thick,variable=\q,domain=1.21315: 1.22611] plot({\q},({1/\q/(pow(\q,9)-1)});

\draw[thick,variable=\q,domain=1.23205: 1.24758] plot({\q},({1/\q/(pow(\q,8)-1)});

\draw[thick,variable=\q,domain=1.25139: 1.25451] plot({\q},{(1/\q+pow(\q,{-8}))/(1-pow(\q,{-15}))-1/\q});

\draw[thick,variable=\q,domain=1.25542: 1.27444] plot({\q},({1/\q/(pow(\q,7)-1)});

\draw[thick,variable=\q,domain=1.27964: 1.28382] plot({\q},{(1/\q+pow(\q,{-7}))/(1-pow(\q,{-13}))-1/\q});

\draw[thick,variable=\q,domain=1.2852: 1.3092] plot({\q},({1/\q/(pow(\q,6)-1)});

\draw[thick,variable=\q,domain=1.31663: 1.3225] plot({\q},{(1/\q+pow(\q,{-6}))/(1-pow(\q,{-11}))-1/\q});

\draw[thick,variable=\q,domain=1.32472: 1.35626] plot({\q},({1/\q/(pow(\q,5)-1)});

\draw[thick,variable=\q,domain=1.35787: 1.36329] plot({\q},{pow(\q,-6)+(1/\q+pow(\q,-5))/(pow(\q,10)-1)});

\draw[thick,variable=\q,domain=1.36759: 1.37635]plot({\q},{(1/\q+pow(\q,{-5}))/(1-pow(\q,{-9}))-1/\q});

\draw[thick,variable=\q,domain=1.38028:1.42421] plot({\q},({1/\q/(pow(\q,4)-1)});

\draw[thick,variable=\q,domain=1.42706: 1.43532] plot({\q},{pow(\q,-5)+(1/\q+pow(\q,-4))/(pow(\q,8)-1)});

\draw[thick,variable=\q,domain=1.44327: 1.45759] plot({\q},{(1/\q+pow(\q,{-4}))/(1-pow(\q,{-7}))-1/\q});

\draw[thick,variable=\q,domain=1.46557:1.53259] plot({\q},({1/\q/(pow(\q,3)-1)});

\draw[dotted, blue]({1.46557},({1/1.46557/(pow(1.46557,3)-1)})--(1.46557,-0.01);
\draw[dotted, blue]({1.53259},({1/1.53259/(pow(1.53259,3)-1)})--(1.53259,-0.01);
\draw ({(1.46557+1.53259)/2},-0.04)node[anchor=south] {$I^{001}$};

\draw[thick,variable=\q,domain=1.5385: 1.55256 ] plot({\q}, {pow(\q,-4) +(1/\q+pow(\q,-3))/(pow(\q,6)-1)});

\draw[thick,variable=\q,domain=1.55392: 1.55759 ] plot({\q}, {pow(\q,-4) +pow(\q,-6)+(1/\q+pow(\q,-3)+pow(\q,-6))/(pow(\q,9)-1)});

\draw[thick,variable=\q,domain=1.57015: 1.5974] plot({\q},{(1/\q+pow(\q,{-3}))/(1-pow(\q,{-5}))-1/\q});

\draw[thick,variable=\q,domain=1.61803:1.73867] plot({\q},({1/\q/(pow(\q,2)-1)});

\draw[dotted, blue]({1.61803},({1/1.61803/(pow(1.61803,2)-1)})--(1.61803,-0.01);
\draw[dotted,blue]({1.73867},({1/1.73867/(pow(1.73867,2)-1)})--(1.73867,-0.01);
\draw ({(1.61803+1.73867)/2},-0.04)node[anchor=south] {$I^{01}$};

\draw[thick,variable=\q,domain=1.7437: 1.75337 ] plot({\q}, {pow(\q,-3)+pow(\q,-5)+(1/\q+pow(\q,-2)+pow(\q,-5))/(pow(\q,6)-1)});

\draw[thick,variable=\q,domain=1.75488: 1.78431 ] plot({\q}, {pow(\q,-3)+(1/\q+pow(\q,-2))/(pow(\q,4)-1)});

\draw[densely dotted,blue]({1.75488},{pow(1.75488,-3)+(1/1.75488+pow(1.75488,-2))/(pow(1.75488,4)-1)})--(1.75488,-0.005);
\draw[densely dotted,blue]({1.78431},{pow(1.78431,-3)+(1/1.78431+pow(1.78431,-2))/(pow(1.78431,4)-1)})--(1.78431,-0.005);
\draw ({(1.75488+1.78431)/2},-0.03)node[anchor=south, scale=0.6pt] {$I^{0011}$};

\draw[thick,variable=\q,domain=1.78854: 1.79758  ] plot({\q}, {pow(\q,-3)+pow(\q,-4)+(1/\q+pow(\q,-2)+pow(\q,-4))/(pow(\q,6)-1)});

\draw[thick,variable=\q,domain=1.79794: 1.80069  ] plot({\q}, {pow(\q,-3)+pow(\q,-4)+pow(\q,-6)+(1/\q+pow(\q,-2)+pow(\q,-4)+pow(\q,-6))/(pow(\q,8)-1)});

\draw[thick,variable=\q,domain=1.8124: 1.83401] plot({\q}, {((1-pow(\q,-2-2))/(\q-1)-pow(\q,-1-2))/(1-pow(\q,-2-3))-1/\q});

\draw[thick,variable=\q,domain=1.83929: 1.9097] plot({\q}, {(pow(\q,-2)-1)/(\q-1)/(pow(\q,-3)-1)-1/\q});

\draw[dotted, blue]({1.83929}, {(pow(1.83929,-2)-1)/(1.83929-1)/(pow(1.83929,-3)-1)-1/1.83929})--(1.83929,-0.01);
\draw[dotted, blue]({1.9097}, {(pow(1.9097,-2)-1)/(1.9097-1)/(pow(1.9097,-3)-1)-1/1.9097})--(1.9097,-0.01);
\draw ({(1.83929+1.9097)/2},-0.04)node[anchor=south] {$I^{011}$};

\draw[thick,variable=\q,domain=1.91118: 1.91971] plot({\q}, {pow(\q,-2)+pow(\q,-4)+pow(\q,-5)+(pow(\q,-1)+pow(\q,-2)+pow(\q,-3)+pow(\q,-5))/(pow(\q,6)-1)});

\draw[thick,variable=\q,domain=1.92213: 1.92712] plot({\q}, {((1-pow(\q,-4-2))/(\q-1)-pow(\q,-2-2))/(1-pow(\q,-4-3))-1/\q});

\draw[thick,variable=\q,domain=1.92756: 1.96223] plot({\q}, {(pow(\q,-3)-1)/(\q-1)/(pow(\q,-4)-1)-1/\q});

\draw[thick,variable=\q,domain=1.96595: 1.98274] plot({\q}, {(pow(\q,-4)-1)/(\q-1)/(pow(\q,-5)-1)-1/\q});

\draw[thick,variable=\q,domain=1.98358: 1.99177] plot({\q}, {(pow(\q,-5)-1)/(\q-1)/(pow(\q,-6)-1)-1/\q});
\draw[thick,variable=\q,domain=1.99196: 1.99598] plot({\q}, {(pow(\q,-6)-1)/(\q-1)/(pow(\q,-7)-1)-1/\q});

\end{tikzpicture}
\end{center}
\caption{The graph of  the critical value function {$\tau(\beta)$ for} $\beta\in(1,2]$. We see that $\tau(\beta)\le 1-1/\beta$ for all $\beta\in(1,2]$, and the function $\tau$ is strictly decreasing in each basic interval $I^\cs$. For example, the basic interval generated by the Farey word $01$ is given by $I^{01}=[\beta_\ell, \beta_*]\approx[1.61803,1.73867]$ with {$((10)^\f)_{\beta_\ell}=(1100(10)^\f)_{\beta_*}=1$}. Furthermore, for any $\beta\in I^{01}$ we have $\tau(\beta)=(00(10)^\f)_\beta=\frac{1}{\beta(\beta^2-1)}$; see Example \ref{ex:cri-farey} for more details. }
\label{fig:graph-tau}
\end{figure}
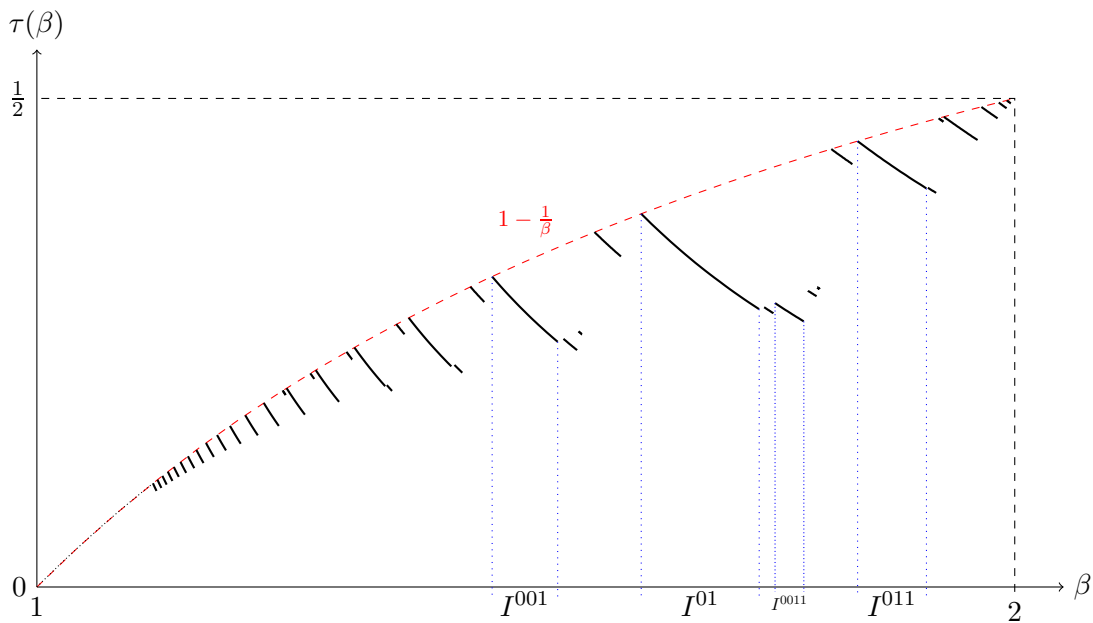

 To describe the critical value $\tau(\beta)$ we first introduce the Farey words, also called standard words (see \cite[Ch.~2.2]{Lothaire-2002}). Following a recent paper of Carminati, Isola and Tiozzo \cite{Carminati-Isola-Tiozzo-2018} we define recursively a sequence of ordered sets $(F_n)_{n=0}^\f$. Let $F_0=(0,1)$, and for $n\ge 0$ the ordered set $F_{n+1}=(v_1,\ldots, v_{2^{n+1}+1})$ is obtained from $F_n=(w_1,\ldots, w_{2^n+1})$ by {inserting} for each $1\le j\le 2^n$ the new word $w_jw_{j+1}$ between the two neighboring words $w_j$ and $w_{j+1}$. So,
 \begin{equation} \label{eq:farey-1}
 \begin{split}
 &F_1=(0,01,1),\qquad F_2=(0,001,01,011,1),\\
  F_3&=(0,0001,001,00101,01,01011,011,0111,1),
 \end{split}
 \end{equation}
 and so on (see Section \ref{sec:preliminaries and Farey words} for more details on Farey words). Set $\F:=\bigcup_{n=1}^\f F_n\setminus F_0$. Then each word in $\F$ is called a non-degenerate  \emph{Farey word}.  Note that any word in $\F$ has length at least two, and begins with digit $0$ and ends with digit $1$. We will use the Farey words as basic bricks to construct infinitely many pairwise disjoint closed intervals, so that we can explicitly  determine $\tau(\beta)$ for $\beta$ in each of these intervals. Furthermore, we will show that these closed intervals cover $(1,2]$ up to a set of zero Hausdorff dimension.

 The construction of these basic intervals depends on certain substitutions of Farey words. For this reason we need to introduce a larger class of words, called Lyndon words; see \cite[Lemma 3.2]{Kalle-Kong-Langeveld-Li-18}.

\begin{definition} \label{def:lyndon words}
  A word $\s=s_1\ldots s_m\in\set{0,1}^*$ is \emph{Lyndon} if
  \[
  s_{i+1}\ldots s_m\succ s_1\ldots s_{m-i}\quad\forall ~0<i<m.
  \]
\end{definition}

Here and throughout the paper we use lexicographical order $\succ$ between sequences and words; see Section \ref{sec:preliminaries and Farey words}. 
{The words $0$ and $1$ are (vacuously) Lyndon.}
Let $\LL$ {denote} the set of all Lyndon words of length at least two. Then by Definition \ref{def:lyndon words}, each $\s\in\LL$ has a prefix $0$ and a suffix $1$.
 It is well-known that each Farey word is Lyndon (cf.~\cite[Proposition 2.8]{Carminati-Isola-Tiozzo-2018}). Thus $\F\subset\LL$.

  Now we define a substitution operator $\bullet$ in $\LL$. This requires the following notation.  By a {\em word} we mean a finite string of zeros and ones. For any two words $\mathbf u=u_1\ldots u_m, \mathbf v=v_1\ldots v_n$ we denote by $\mathbf u\mathbf v=u_1\ldots u_m v_1\ldots v_n$ their concatenation. Furthermore, we write $\mathbf u^\f$ for the periodic sequence with periodic block $\mathbf u$.
 For a word $\mathbf w=w_1\ldots w_n\in\set{0,1}^n$ we denote $\mathbf w^-:=w_1\ldots w_{n-1}{0}$ if $w_n=1$, and $\mathbf w^+:=w_1\ldots w_{n-1}{1}$ if $w_n=0$. Furthermore, we denote by $\L(\mathbf w)$ the lexicographically largest cyclic permutation of $\mathbf w$.
  Now for two words $\s=s_1\ldots s_m\in\LL$ and $\r=r_1\ldots r_\ell\in\set{0,1}^\ell$ we define
 \begin{equation}  \label{eq:substition-operator}
   \s\bullet\r:=c_1\ldots c_{\ell m},
 \end{equation}
 where 
 \[c_1\ldots c_m=\left\{\begin{array}
   {lll}
   \s^-&\textrm{if}& r_1=0\\
   \L(\s)^+&\textrm{if}& r_1=1,
 \end{array}\right.\]
 and for $1\le j<\ell$,
 \[
 c_{jm+1}\ldots c_{(j+1)m}=\left\{\begin{array}
   {lll}
   \L(\s)&\textrm{if}& r_jr_{j+1}=00\\
   \L(\s)^+&\textrm{if}& r_jr_{j+1}=01\\
 \s^-&\textrm{if}& r_{j}r_{j+1}=10\\
    \s &\textrm{if}& r_jr_{j+1}=11.
 \end{array}\right.
 \]
For an equivalent definition of the substitution operator $\bullet$, see Section \ref{sec:substitution}.

\begin{example} \label{ex:substitution}
 Let $\r=01, \s=001$ and $\t=011$ be three words in $\F$.
 Then $\L(\r)=10$ and $\L(\s)=100$. So, by (\ref{eq:substition-operator}) it follows that
 \begin{align*}
 \r\bullet\s&=\r\bullet 001=\r^-\L(\r)\L(\r)^+=001011\in\LL,\\
  \s\bullet\t&=\s\bullet  011=\s^-\L(\s)^+\s=000\,101\,001\in\LL.
 \end{align*}
Then $\L(\r\bullet\s)=110010$, and thus
 \[
 (\r\bullet\s)\bullet\t=(\r\bullet\s)\bullet 011=(\r\bullet\s)^-\L(\r\bullet\s)^+(\r\bullet\s)=001010\,110011\,001011,
 \]
 and
 \begin{align*}
 \r\bullet(\s\bullet\t)&=\r\bullet 000101001 \\
 &= \r^-\L(\r)\L(\r)\L(\r)^+\r^-\L(\r)^+\r^-\L(\r)\L(\r)^+=00\,10\,10\,11\,00\,11\,00\,10\,11.
 \end{align*}
 Hence, $(\r\bullet\s)\bullet\t=\r\bullet(\s\bullet\t)$, suggesting that the operator $\bullet$ is associative. On the other hand, observe that
$\r\bullet\s=00\, 10\, 11\ne 000\,101=\s\bullet\r$.
So $\bullet$ is not commutative.
\end{example}

From Example \ref{ex:substitution} we see that $\F$ is not closed under the substitution operator $\bullet$, since $\r\bullet\s=001011\not\in\F$. Hence we need the larger collection $\LL$. It turns out that $\LL$ is a non-Abelian semi-group under the substitution operator $\bullet$.

  \begin{proposition}  \label{main:substitution}
    $(\LL,\bullet)$ forms a non-Abelian semi-group.
  \end{proposition}

\begin{remark}
The substitution operator $\bullet$ defined in (\ref{eq:substition-operator}) {is} similar {to} that {introduced} by Allaart \cite{Allaart-2019}, who used it to study the entropy plateaus in unique $q$-expansions.
\end{remark}
	
Let
\begin{equation} \label{eq:set-substitution-Farey-words}
\La:=\set{\cs=\s_1\bullet\s_2\bullet\cdots\bullet\s_k:~~ \s_i\in\F\textrm{ for any }1\le i\le k; ~k\in\N}
\end{equation}
be the set of all substitutions of Farey words from $\F$. Then by Proposition \ref{main:substitution} it follows that {$\F\subset\La\subset\LL$}. Moreover, both inclusions are strict. For instance, $001011=01\bullet 001\in\La\backslash \F$ by Example \ref{ex:substitution} and Proposition \ref{prop:Farey-new-form} below, and $0010111\in\LL\backslash\La$.

Given $\beta\in(1,2]$, for a sequence $(c_i)\in\set{0,1}^\N$ we write
 \[
((c_i))_{\beta}:=\sum_{i=1}^\f\frac{c_i}{\beta^i}.
\]
Now we define the basic intervals.

 \begin{definition}  \label{def:basic-intervals}
   A closed interval $I=[\beta_\ell,\beta_*]\subset(1,2]$ is called a \emph{basic interval} if there exists a word $\cs\in\La$ such that
   \[
   (\L(\cs)^\f)_{\beta_\ell}=1\quad\textrm{and}\quad (\L(\cs)^+\cs^-\L(\cs)^\f)_{\beta_*}=1.
   \]
   The interval $I=I^\cs$ is also called a basic interval generated by $\cs$.
 \end{definition}

The subscripts for the endpoints $\beta_\ell$ and $\beta_*$ of a basic interval will be clarified when we define the Lyndon intervals (see Definition \ref{def:lyndon-intervals} below). Our second main result gives a formula for $\tau(\beta)$ {when $\beta$ lies} in a basic interval $I^\cs$.

 \begin{main} \label{main:critical-basic-intervals}\mbox{}
	\begin{enumerate} [{\rm(i)}]
\item The basic intervals $I^\cs, \cs\in\La$ are pairwise disjoint.
\item If $I^\cs$ is a basic interval generated by $\cs\in\La$, then
   \begin{equation} \label{eq:tau-in-basic-interval}
   \tau(\beta)=(\cs^-\L(\cs)^\f)_\beta\quad\textrm{for {every}}\quad \beta\in I^\cs.
   \end{equation}
\item The function $\tau$ is strictly decreasing on $I^\cs$, and is real-analytic and strictly convex in the interior of $I^\cs$.
\end{enumerate}
 \end{main}

\begin{remark}
  {Note that (iii) follows immediately from (ii). For the special case when $\cs\in\F$, the formula \eqref{eq:tau-in-basic-interval} was stated without proof by Kalle et al.~\cite{Kalle-Kong-Langeveld-Li-18}.}
\end{remark}

\begin{example} \label{ex:cri-farey}\mbox{}
  \begin{enumerate}[{\rm(i)}]
  \item  Let $\s=01\in\F$. Then by Definition \ref{def:basic-intervals} the basic interval $I^{01}=[\beta_\ell,\beta_*]$ satisfies
  \[
  (\L(01)^\f)_{\beta_\ell}=((10)^\f)_{\beta_\ell}=1\quad\textrm{and}\quad(\L(01)^+(01)^-\L(01)^\f)_{\beta_*}=(1100(10)^\f)_{\beta_*}=1.
  \]
  By numerical calculation we get $I^{01} \approx[1.61803, 1.73867]$ (see Figure \ref{fig:graph-tau}). In fact, $\beta_\ell=(1+\sqrt{5})/2$. Theorem \ref{main:critical-basic-intervals} {yields} that
  \[
  \tau(\beta)=(00(10)^\f)_\beta=\frac{1}{\beta(\beta^2-1)}\qquad\forall\beta\in I^{01}.
  \]

\item Let $\s_1=\s_2=01\in\F$. Then $\cs=\s_1\bullet\s_2=01\bullet  01=0011$. By Definition \ref{def:basic-intervals} the basic interval $I^{\s_1\bullet\s_2}=I^{0011}=[\beta_\ell,\beta_*]$ is  given implicitly by
  \begin{align*}
  (\L(0011)^\f)_{\beta_\ell}&=((1100)^\f)_{\beta_\ell}=1,\\
  \quad (\L(0011)^+(0011)^-\L(0011)^\f)_{\beta_*}&=(11010010(1100)^\f)_{\beta_*}=1.
  \end{align*}
  Numerical calculation gives $I^{0011}\approx[1.75488, 1.78431]$ (see Figure \ref{fig:graph-tau}), {and} Theorem \ref{main:critical-basic-intervals} {implies}
  \[
  \tau(\beta)=(\cs^-\L(\cs)^\f)_\beta=(0010(1100)^\f)_\beta=\frac{1}{\beta^3}+\frac{1+\beta}{\beta^2(\beta^4-1)}\qquad \forall\beta\in I^{0011}.
  \]
  \end{enumerate}
\end{example}

 Next, we introduce the Lyndon intervals.

 \begin{definition} \label{def:lyndon-intervals}
 For each Lyndon word $\cs\in\LL$ the interval $J^\cs=[\beta_\ell^\cs,\beta_r^\cs]\subset(1,2]$ is called a \emph{Lyndon interval} generated by $\cs$ if
 \[
 (\L(\cs)^\f)_{\beta_\ell^\cs}=1\quad\textrm{and}\quad (\L(\cs)^+\cs^\f)_{\beta_r^\cs}=1.
 \]
If in particular $\cs\in\F$, we call $J^\cs$ a {\em Farey interval}.
 \end{definition}
We remark that the Farey intervals defined in \cite[Definition 4.5]{Kalle-Kong-Langeveld-Li-18} are half-open intervals, which is slightly different from our definition.
It turns out that the discontinuity points of $\tau$ are precisely the right endpoints of the Lyndon intervals $J^\cs$ with $\cs\in\La$:

\begin{proposition} \label{prop:discontinuities}
The function $\tau$ is continuous on $(1,2]\backslash \{\beta_r^\cs: \cs\in\La\}$. On the other hand, for each $\cs\in\La$, we have
\begin{equation} \label{eq:upward-jump}
\lim_{\beta\searrow\beta_r^\cs}\tau(\beta)=(\cs^\f)_{\beta_r^\cs}>(\cs 0^\f)_{\beta_r^\cs}=\tau(\beta_r^\cs).
\end{equation}
\end{proposition}

\begin{remark}
In a previous version of this article (as well as in the published version), we concluded erroneously from Proposition \ref{prop:discontinuities} that the dimension $\dim_H K_\beta(t)$ is not jointly continuous in $\beta$ and $t$. We also stated that when $t=\tau(\beta_r^\cs)$ for $\cs\in\La$, the function $\beta\mapsto \dim_H K_\beta(t)$ has a jump at $\beta_r^\cs$. This is incorrect. It is in fact not difficult to show that $\dim_H K_\beta(t)$ is continuous in $\beta$ for all fixed $t$. We suspect that the function $(\beta,t)\mapsto \dim_H K_\beta(t)$ is jointly continuous, but have not been able to prove this.
\end{remark}

It was shown in \cite[{Section 4}]{Kalle-Kong-Langeveld-Li-18} that the Farey intervals $J^\s, \s\in\F$ are pairwise disjoint, and the {\em {exceptional} set}
\[
E:=(1,2]\setminus\bigcup_{\s\in\F}J^\s
\]
has zero Hausdorff dimension. We strengthen this result slightly and show in Proposition \ref{prop:geometry-basic-interval-2} (i) that $E$ is uncountable and has zero packing dimension.

From Definitions \ref{def:basic-intervals} and \ref{def:lyndon-intervals} it follows that $I^\cs\subset J^{\cs}$ for any $\cs\in\La$, and the two intervals $I^\cs$ and $J^\cs$ have the same left endpoint (see Proposition \ref{prop:basic-interval-disjoint}.) In Proposition \ref{prop:geometry-basic-interval-2} (ii) we show that for any $\cs\in{\La}$ {the Lyndon intervals $J^{\cs\bullet\r}, \r\in\F$ are pairwise disjoint subsets of $J^\cs\setminus I^\cs$, and} the \emph{relative {exceptional} set}
 \[
 E^\cs:=(J^\cs\setminus I^\cs)\setminus\bigcup_{\r\in\F}J^{\cs\bullet\r}
 \]
 is also uncountable and has zero box-counting dimension. In Proposition \ref{prop:basic-interval-disjoint} we {show} that the Lyndon intervals $J^\cs, \cs\in\La$ have a tree structure. This gives rise to the {set}
\begin{equation} \label{eq:E_s}
 E_\f:=\bigcap_{k=1}^\f\bigcup_{\cs\in\La(k)}J^{\cs},
 \end{equation}
 where $\La(k):=\set{\cs=\s_1\bullet\cdots\bullet\s_k:~ \s_i\in\F~\forall 1\le i\le k}$. {We call $E_\f$ the {\em infinitely Farey set}, because its elements arise from substitutions of an infinite sequence of Farey words.}
It follows at once that $E_\f$ is uncountable; we show in Proposition \ref{prop:dim-exceptional-set} that it has zero Hausdorff dimension.

Combining the above results we obtain our last main theorem:

\begin{main} \label{main:geometrical-structure-basic-intervals}
  The interval $(1,2]$ can be partitioned as
 \[
  (1,2]=E\cup E_\f\cup\bigcup_{\cs\in\La}E^{\cs}\cup\bigcup_{\cs\in\La}I^{\cs},
  \]
and the  basic intervals $\set{I^{\cs}: \cs\in\La}$ cover $(1,2]$ up to a set of zero Hausdorff dimension.
\end{main}

\begin{remark}
  It is worth mentioning that the Lyndon intervals $J^\cs$ and the relative {exceptional} sets $E^\cs$ constructed in our paper have similar geometrical structure as the relative entropy plateaus and relative bifurcation sets studied in \cite{Allaart-Kong-2018}, where they were used to describe the local dimension of the set of univoque bases.
\end{remark}

The following result was established in the proof of \cite[Theorem D]{Kalle-Kong-Langeveld-Li-18}.

\begin{proposition}  \label{lem:critical-non-Fareyinterval}
For any $\beta\in(1,2]$ we have $\tau(\beta)\le 1-\frac{1}{\beta}$. Furthermore,
\[
\tau(\beta)=1-\frac{1}{\beta}\quad\textrm{for any}\quad \beta\in E.
\]
\end{proposition}

Thus, in view of Theorem \ref{main:geometrical-structure-basic-intervals} it remains to determine $\tau(\beta)$ for $\beta\in E^\cs$ with $\cs\in\La$ and for $\beta\in E_\f$. In Proposition \ref{prop:cric-exception-1} we compute $\tau(\beta)$ for $\beta\in E^\cs$ by relating the relative {exceptional} set $E^\cs$ to the {exceptional} set $E$ via a renormalization map $\Psi_\cs$. Proposition \ref{prop:cri-exception-2} gives an expression for $\tau(\beta)$ when $\beta\in E_\f$. As an illustration of the latter, in Proposition \ref{prop:cri-exception-3} we construct in each Farey interval $J^\s$ a transcendental base $\beta_\f^\s\in E_\f$ and give {an} explicit formula for $\tau(\beta_\f^\s)$. Here we point out an interesting connection with unique $\beta$-expansions: Let $\beta\approx 1.78723$ be the \emph{Komornik-Loreti constant} (cf.~\cite{KomornikLoreti1998}); that is, $\beta$ is the smallest base in which 1 has a unique expansion. Then it follows from Proposition \ref{prop:cri-exception-3} that $\beta=\beta_\f^{01}\in E_\f$, and $\tau(\beta)=\frac{2-\beta}{\beta-1}\approx 0.270274$.

The rest of the paper is organized as follows. In Section \ref{sec:preliminaries and Farey words} we recall some properties of Farey words {and} Farey intervals, as well as greedy and quasi-greedy $\beta$-expansions.
 In Section \ref{sec:substitution} we give an equivalent definition of the substitution operator $\bullet$, and prove Proposition \ref{main:substitution}. The proof of Theorem \ref{main:critical-basic-intervals} is given in Section \ref{sec:basic-intervals}. {At the heart of the argument is Proposition \ref{prop:countable}, which clarifies the role of the special Lyndon words $\cs\in\La$ and is used in several settings to derive the upper bound for $\tau(\beta)$.
The {relative exceptional} sets $E^\cs, \cs\in\La$ and {the infinitely Farey set} $E_\f$ are studied in detail in Section \ref{sec:geometric-structure}, where we show that all of these sets have zero Hausdorff dimension}, proving Theorem \ref{main:geometrical-structure-basic-intervals}. In Section \ref{sec:cri-bifurcation-sets} we determine the critical value $\tau(\beta)$ for $\beta$ in the relative {exceptional} sets $E^\cs$ and the {infinitely Farey set} $E_\f$. Finally, in Section \ref{sec:continuity} we show that the function $\beta\mapsto \tau(\beta)$ is {c\`adl\`ag}, and prove Proposition \ref{prop:discontinuities} and Theorem \ref{main:critical-devils-staircase}.

\section{Farey words and Farey intervals} \label{sec:preliminaries and Farey words}

In this section we recall some properties of Farey words, which {are} vital in determining the critical value $\tau(\beta)$. We also recall from \cite{Kalle-Kong-Langeveld-Li-18} the Farey intervals, {and review basic properties of greedy and quasi-greedy $\beta$-expansions.}

First we introduce some terminology from symbolic dynamics (cf.~\cite{Lind_Marcus_1995}). Let $\set{0,1}^\N$ be the set of all infinite sequences of zeros and ones. Denote by $\si$ the left shift map. Then $(\set{0,1}^\N, \si)$ is a full shift. {By} a \emph{word} we mean a finite string of zeros and ones. Let $\set{0, 1}^*$ be the set of all words over the alphabet $\set{0,1}$ together with the empty word $\epsilon$. For a word $\c\in\set{0,1}^*$ we denote its length by $|\c|$, and for a digit $a\in\set{0,1}$ we denote by $|\c|_a$ the number of occurrences of $a$ in the word $\c$. For two words $\mathbf c=c_1\ldots c_m$ {and} $\mathbf d=d_1\ldots d_n$ {in} $\set{0,1}^*$ we write $\mathbf{cd}=c_1\ldots c_md_1\ldots d_n$ for their concatenation. For $n\in\N$ we denote by $\mathbf c^n$ the $n$-fold concatenation of $\mathbf c$ with itself, and by $\mathbf c^\f$ the periodic sequence with {period} block $\mathbf c$.

Throughout the paper we will use the lexicographical order `$\prec, \lle, \succ$' or `$\lge$' between sequences and words. For example, for two sequences $(c_i), (d_i)\in\set{0,1}^\N$, we say $(c_i)\prec (d_i)$ if $c_1<d_1$, or there exists $n\in\N$ such that $c_1\ldots c_n=d_1\ldots d_n$
 and $c_{n+1}<d_{n+1}$. For two words $\mathbf c, \mathbf d$, we say $\mathbf c\prec \mathbf d$ if $\mathbf c 0^\f\prec \mathbf d 0^\f$. We also recall from Section \ref{sec: Introduction} that if $\mathbf c=c_1\ldots c_m$ with $c_m=0$, then $\mathbf c^+ =c_1\ldots c_{m-1}1$; and if $\mathbf c=c_1\ldots c_m$ with $c_m=1$, then $\mathbf c^- =c_1\ldots c_{m-1}0$. {Finally}, for a word {$\mathbf{c}=c_1 c_2\ldots c_n$ we denote its \emph{reflection} by $\overline{\mathbf{c}}:=(1-c_1)(1-c_2)\ldots(1-c_n)$}.

\subsection{Farey words}
 Farey words {have attracted much} attention in the literature due to their intimate connection with rational rotations on the circle (see \cite[Chapter 2]{Lothaire-2002}) and their one-to-one correspondence with the rational numbers in $[0,1]$ (see {(\ref{eq:kk-6})} below). In the following we adopt the definition from a recent paper of Carminati, Isola and Tiozzo \cite{Carminati-Isola-Tiozzo-2018}.

 First we recursively define a sequence of ordered sets $F_n, n=0,1,2,\ldots$. Let $F_0=(0,1)$; and for $n\ge 0$ the ordered set $F_{n+1}=(v_1, \ldots, v_{2^{n+1}+1})$ is obtained from $F_{n}=(w_1,\ldots, w_{2^n+1})$ by
 \[
 \left\{
 \begin{array}
   {lll}
   v_{2i-1}=w_i&\textrm{for}& 1\le i\le 2^{n}+1,\\
   v_{2i}=w_iw_{i+1}&\textrm{for}& 1\le i\le 2^n.
 \end{array}\right.
 \]
 In other words, $F_{n+1}$ is obtained from $F_n$ by {inserting} for each $1\le j\le 2^{n}$ the new word $w_jw_{j+1}$ between the two neighboring words $w_j$ and $w_{j+1}$. See (\ref{eq:farey-1}) for examples. Note that for each $n\ge 0$ the ordered set $F_n$ consists of $2^n+1$ words which are listed from the left to the right in lexicographically increasing order.
We call $w\in\set{0,1}^*$ a \emph{Farey word} if $w\in F_n$ for some $n\ge 0$, and we denote by $\Omega_F:=\bigcup_{n=1}^\f F_n$ the set of all Farey words. As shown in \cite[Proposition 2.3]{Carminati-Isola-Tiozzo-2018}, the set $\Omega_F$ can be bijectively mapped to $\mathbb{Q}\cap[0,1]$ via the map
 \begin{equation} \label{eq:kk-6}
 \xi: \Omega_F\to\mathbb Q\cap[0,1];\quad \s\mapsto\frac{|\s|_1}{|\s|}.
 \end{equation}
  So, $\xi(\s)$ is the frequency of the digit $1$ in $\s$.

 For each $n\ge 1$ set
 \[
 F_n^*:=F_n\setminus\set{0,1},
 \]
 and
 \[
 F_n^0:=\set{w\in F_n^*: |w|_0>|w|_1},\quad F_n^1:=\set{w\in F_n^*: |w|_1>|w|_0}.
 \]
 For example, $F_1^*=(01), F_2^*=(001, 01, 011)$, and $F_2^0=(001), F_2^1=(011)$.
 The following decomposition can be deduced from \cite[Proposition 2.3]{Carminati-Isola-Tiozzo-2018}.

 \begin{lemma} \label{lem:farey-Fn}
   For any $n\ge 2$ we have $F_n^*=F_n^0\cup F_1^*\cup F_n^1$.
 \end{lemma}

 The ordered sets $F_n^*,n\ge 1$ can also be obtained via substitutions. We define the two substitution operators  by
 \begin{equation} \label{eq:subs-01}
 U_0:\left\{\begin{array}{lll}
 0&\mapsto&0\\
 1&\mapsto&01,
 \end{array}\right.\quad\textrm{and}\quad U_1:\left\{
 \begin{array}{lll}
 0&\mapsto&01\\
 1&\mapsto&1.
 \end{array}\right.
 \end{equation}
Then $U_0$ and $U_1$ naturally induce a map  on $\set{0,1}^*$ or $\set{0,1}^\N$. For example,
 \[
 U_0: \set{0,1}^*\to\set{0,1}^*;\quad c_1\ldots c_n\mapsto U_0(c_1)\ldots U_0(c_n).
\]
 The following result was proven in \cite[Proposition 2.9]{Carminati-Isola-Tiozzo-2018}.

\begin{lemma}  \label{lem:Farey-substitution}
   For each $a\in\set{0,1}$ the map $U_a: F_n^*\to F_{n+1}^a$ is bijective.
 \end{lemma}

 By Lemmas \ref{lem:farey-Fn} and \ref{lem:Farey-substitution} it follows that the ordered sets $F_n^*$ can be obtained by the substitution operators $U_0$ and $U_1$ on the set $F_1^*=(01)$. We will clarify this in the next proposition.
Let $\F$ be the set of all non-degenerate Farey words, i.e.,
\[
\F=\bigcup_{n=1}^\f F_n^*.
\]

For a word $\mathbf c=c_1\ldots c_m\in\set{0,1}^*$ let $\S(\mathbf c)$ and $\L(\mathbf c)$ be the lexicographically smallest and largest cyclic permutations of $\mathbf c$, respectively. In other words, $\S(\mathbf c)$ is the lexicographically smallest word among
\[
c_1c_2\ldots c_m, \quad c_2\ldots c_m c_1,\quad c_3\ldots c_m c_1c_2,\quad \cdots,\quad c_mc_1\ldots c_{m-1};
\]
and $\L(\mathbf c)$ is the lexicographically largest word in the above list. The following {properties of Farey words are} well known (see, e.g., \cite[Proposition 2.5]{Carminati-Isola-Tiozzo-2018}).

\begin{lemma}  \label{lem:Farey-property}
Let $\s=s_1\ldots s_m\in\F$. Then
\begin{enumerate}[{\rm (i)}]
\item $\S(\s)=\s$ and $\L(\s)=s_ms_{m-1}\ldots s_1$.
\item ${\s^-}$ is a palindrome; that is $s_1\ldots s_{m-1}(s_m-1)=(s_m-1)s_{m-1}s_{m-2}\ldots s_1$.
\item The word $\s$ has a {\em conjugate} $\tilde{\s}\in\F$, given by
\begin{equation} \label{eq:Farey-conjugate}
\tilde{\s}:=\overline{\L(\s)}=0\,\overline{s_2\ldots s_{m-1}}\,1.
\end{equation}
\end{enumerate}
\end{lemma}

The last equality in \eqref{eq:Farey-conjugate} follows from statements (i) and (ii). In terms of the correspondence \eqref{eq:kk-6}, if $\xi(\s)=r\in\mathbb{Q}\cap[0,1]$, then $\xi(\tilde{\s})=1-r$.
Note also that the conjugate of $\tilde{\s}$ is simply $\s$ itself.

The following explicit description of $\F$ will be useful in Section \ref{sec:basic-intervals} to prove the upper bound for $\tau(\beta)$.

\begin{proposition} \label{prop:Farey-new-form}
   $\F$ consists of all words in one of the following forms:
   \begin{itemize}
    \item[{\rm(i)}] $01^p$ or $0^p 1$ for some $p\in\N$;
		\item[{\rm(ii)}] $01^p 01^{p+t_1}\ldots 01^{p+t_N}01^{p+1}$ for some $p\in\N$ and Farey word $0 t_1\ldots t_N1\in\F$;
    \item[{\rm(iii)}] $0^{p+1}10^{p+t_1}1\ldots 0^{p+t_N}10^p 1$ for some $p\in\N$ and Farey word $0 t_1\ldots t_N1\in\F$.
   \end{itemize}
\end{proposition}

\begin{proof}
Note that $01=U_1(0)=U_0(1)\in F_1^*\subset\F$. Furthermore, for $p\in\N$ and $0 t_1\ldots t_N1\in\F$ we have
\begin{align*}
01^p=U_1(01^{p-1})&=U_1^{p-1}(U_0(1)),\\
0^p1=U_0(0^{p-1}1)&=U_0^{p-1}(U_1(0)),\\
01^p 01^{p+t_1}\ldots 01^{p+t_N}01^{p+1}&=U_1^p\big(U_0(0t_1\ldots t_{N}1)\big),\\
0^{p+1}10^{p+t_1}1\ldots 0^{p+t_N}10^p 1&=U_0^p\big(U_1(0\,\overline{t_1\ldots t_N}\,1)\big).
\end{align*}
By Lemma \ref{lem:Farey-property} (iii), if $0 t_1\ldots t_N1\in\F$ then $0\,\overline{t_1\ldots t_N}\,1\in\F$ as well.
Hence by Lemma \ref{lem:Farey-substitution}, all the above words lie in $\F$.

 To prove the converse, it suffices to show that each word in $\F$ is of the form $U_0^p(U_1(\t))$ or $U_1^p(U_0(\t))$ for some $p\geq 0$ and Farey word $\t\in\Omega_F$. This is clearly true for $01=U_0^0(U_1(0))$, where $U_0^0$ denotes the identity map. Let $n\geq 1$ and suppose the statement is true for all Farey words in $F_n^*$. Take $\s\in F_{n+1}^*$ with $\s\neq 01$. By Lemmas \ref{lem:farey-Fn} and \ref{lem:Farey-substitution}, $\s=U_0(\t)$ or $\s=U_1(\t)$ for some Farey word $\t\in F_n^*$. We assume the former, as the argument for the second case is similar. By the induction hypothesis, either $\t=U_0^p(U_1(\mathbf{u}))$ for some $\mathbf{u}\in\Omega_F$ and $p\geq 0$, in which case $\s=U_0^{p+1}(U_1(\mathbf{u}))$; or $\t=U_1^p(U_0(\mathbf{u}))$ for some $\mathbf{u}\in\Omega_F$ and $p\geq 1$, in which case $\s=U_0(U_1(\mathbf{v}))$ where $\mathbf{v}=U_1^{p-1}(U_0(\mathbf{u}))\in \Omega_F$. In both cases, $\s$ is of the required form.
\end{proof}

Observe that the two types of words in Proposition \ref{prop:Farey-new-form} (i) are each others conjugates, and the conjugate of a Farey word of type (ii) is a Farey word of type (iii), and vice versa.
For more properties of Farey words we refer to the book of Lothaire \cite{Lothaire-2002} and the references therein.

\subsection{{Quasi-greedy expansions}, Farey intervals and Lyndon intervals}

 Given $\beta\in(1,2]$, let $\de(\beta)=\de_1(\beta)\de_2(\beta)\ldots\in\set{0,1}^\N$ be the \emph{quasi-greedy} $\beta$-expansion of $1$ (cf.~\cite{Daroczy_Katai_1993}), i.e., $\de(\beta)$ is the lexicographically largest sequence not ending with $0^\f$ such that $(\de_i(\beta))_\beta=1$.
The following  property of $\de(\beta)$ is well known (cf.~\cite{Baiocchi_Komornik_2007}).

\begin{lemma} \label{lem:quasi-greedy expansion}\mbox{}
\begin{enumerate}[{\rm(i)}]
\item The map $\beta\mapsto\de(\beta)$ is an increasing bijection from $\beta\in (1,2]$ to the set of sequences $(a_i)\in\set{0,1}^\N$ not ending with $0^\f$ and satisfying
\[
\si^n((a_i))\lle (a_i)\quad \forall n\ge 0.
\]
\item The map $\beta\mapsto\de(\beta)$ is left continuous everywhere on $(1,2]$ with respect to the order topology, and it  is right continuous at $\beta_0\in(1,2)$ if and only if $\de(\beta_0)$ is not periodic. Furthermore, if $\de(\beta_0)=(a_1\ldots a_m)^\f$ with {minimal} period $m$, then $\de(\beta)\searrow a_1\ldots a_m^+ 0^\f$ as $\beta\searrow \beta_0$.
\end{enumerate}
\end{lemma}

Recall from Definition \ref{def:lyndon words} that for a word $\s=s_1\ldots s_m \in\LL$ we have  $s_{i+1}\ldots s_m\succ s_1\ldots s_{m-i}$ for all $1\le i<m$.
  The following basic fact can be found in \cite[Theorem 1.5.3]{Alloche_Shallit_2003}:

\begin{lemma} \label{lem:cyclic-permutations}
Let $\c=c_1\ldots c_m\in\{0,1\}^*$, and suppose two cyclic permutations of $\c$ are equal (that is, $c_{i+1}\ldots c_m c_1\ldots c_i=c_{j+1}\ldots c_m c_1\ldots c_j$ where $i\neq j$). Then $\c$ is periodic; in other words, $\c=\b^k$ for some word $\b$ and $k\geq 2$.
\end{lemma}

In fact, the length of $\b$ in Lemma \ref{lem:cyclic-permutations} can be taken to equal $\gcd(|i-j|,m)$.

\begin{lemma} \label{lem:lyndon-equivalence}
Let $\s \in\LL$ and $\a=\L(\s)=a_1\ldots a_m$. Then
  \begin{equation} \label{eq:inequality-ai}
  a_{i+1}\ldots a_m\prec a_1\ldots a_{m-i}\quad\forall\ 1\le i<m.
  \end{equation}
Furthermore,
\begin{equation}\label{eq:inequality-2}
\si^n(\a^+{\s^-}\a^\f)\lle \a^+{\s^-}\a^\f\quad\forall n\ge 0.
\end{equation}
\end{lemma}

\begin{proof}
First we prove (\ref{eq:inequality-ai}). Since $\s$ is Lyndon, it is not periodic. Hence $\a=\L(\s)$ is not periodic, because any cyclic permutation of a periodic word is periodic. Since $\a=\L(\s)$, we have
\[
a_{i+1}\ldots a_m\lle a_1\ldots a_{m-i}\quad\forall\ 1\le i<m.
\]
Suppose equality holds for some $i$. Then
\[
a_{i+1}\ldots a_m a_1\ldots a_i=a_1\ldots a_{m-i} a_1\ldots a_i\lge a_1\ldots a_{m-i} a_{m-i+1}\ldots a_m=\a,
\]
so $a_{i+1}\ldots a_m a_1\ldots a_i=\L(\s)=\a$ by definition of $\L(\s)$. By Lemma \ref{lem:cyclic-permutations}, this cannot happen, since $\a$ is not periodic.

Next we prove (\ref{eq:inequality-2}).
 Since $\s=s_1\ldots s_m$ is a Lyndon word, any word of length $k\in\set{1, \ldots, m-1}$ occurring in $\a=\L(\s)$ is lexicographically larger than or equal to $s_1\ldots s_k$. By (\ref{eq:inequality-ai}) it follows that
  \begin{equation} \label{eq:29-2}
    a_{k+1}\ldots a_m^+ s_1\ldots s_k\lle a_1\ldots a_{m-k}a_{m-k+1}\ldots a_m\prec a_1\ldots a_m^+
  \end{equation}
  for all $0<k<m$.  Hence, by  (\ref{eq:29-2}) and (\ref{eq:inequality-ai}) we conclude that
  $\si^n(\a^+{\s^-}\a^\f)\prec \a^+{\s^-}\a^\f$ for all $n\ge 1$. This completes the proof.
\end{proof}

\begin{lemma} \label{lem:periodic-qg-expansion}
Let $\beta\in(1,2)$. Then $\de(\beta)$ is periodic if and only if $\de(\beta)=\L(\s)^\f$ for some Lyndon word $\s$ of length at least two.
\end{lemma}

\begin{proof}
Suppose $\de(\beta)=(a_1\ldots a_m)^\f$ with minimal period block $\a=a_1\ldots a_m$. Then $m\ge 2$ {since $\beta<2$}. Take $\s:=\S(\a)$. Then $\a=\L(\s)$, and
\begin{equation} \label{eq:almost-Lyndon}
s_{i+1}\ldots s_m\lge s_1\ldots s_{m-i} \qquad \forall\ 1\leq i<m.
\end{equation}
If equality holds for some $i$, then we deduce just as in the proof of Lemma \ref{lem:lyndon-equivalence} that $\s$ is periodic. But then $\a$ is also periodic, contradicting that $m$ is the minimal period of $\de(\beta)$. Hence, strict inequality holds in \eqref{eq:almost-Lyndon}, and $\s$ is Lyndon. The converse is trivial.
\end{proof}

Recall the Farey intervals and Lyndon intervals {from} Definition \ref{def:lyndon-intervals}.
The following properties of Lyndon intervals and Farey intervals were established in the proof of \cite[Theorem C]{Kalle-Kong-Langeveld-Li-18}.

\begin{lemma} \label{lem:farey-interval}\mbox{}

  \begin{itemize}
  \item[{\rm(i)}] The Farey intervals $J^\s, \s\in\F$ are pairwise disjoint, {and their union is dense in $(1,2]$}. 

  \item[{\rm(ii)}] {Any two Lyndon intervals are either disjoint}, or one is contained in the other.

  \item[{\rm(iii)}] For any Lyndon interval $J^\cs, \cs\in\LL$ there exists a unique Farey interval {$J^\r$} such that $J^\cs\subset {J^\r}$.
  \end{itemize}
\end{lemma}

{Note that (iii) follows immediately from (i) and (ii).}

\subsection{Greedy expansions and the symbolic survivor set}

Given $\beta\in(1,2]$ and $t\in[0,1)$, we call the sequence $(d_i)\in\set{0,1}^\N$ a \emph{$\beta$-expansion} of $t$ if
$((d_i))_\beta=t$.
Note that a point $t\in[0,1)$ may have  multiple $\beta$-expansions. {We denote by} $b(t,\beta)=(b_i(t,\beta))\in\set{0,1}^\N$ the \emph{greedy} $\beta$-expansion of $t$, which is the lexicographically largest expansion of $t$ in base $\beta$. {Since $T_\beta(t)=\beta t\pmod 1$, it follows that} $b(T_\beta^n(t), \beta)=\si^n(b(t,\beta))=b_{n+1}(t,\beta)b_{n+2}(t,\beta)\ldots$.
The following result was established by Parry \cite{Parry_1960} and de Vries and Komornik \cite[Lemma 2.5 and Proposition 2.6]{DeVries-Komornik-2011}.

\begin{lemma} \label{lem:greedy-expansion}
Let $\beta\in(1,2]$. The map $t\mapsto b(t,\beta)$ is an increasing bijection from $[0,1)$ to
\[
\big\{(d_i)\in\set{0,1}^\N: \si^n((d_i))\prec \de(\beta)~\forall n\ge 0\big\}.
\]
Furthermore,
\begin{itemize}
  \item [{\rm(i)}] The map $t\mapsto b(t,\beta)$ is right-continuous everywhere in $[0,1)$ with respect to the order topology in $\set{0,1}^\N$;
  \item [{\rm(ii)}] If $b(t_0,\beta)$ does not end with $0^\f$, then the map $t\mapsto b(t,\beta)$ is continuous at $t_0$;
	\item [{\rm(iii)}] If $b(t_0,\beta)=b_1\ldots b_m 0^\f$ with $b_m=1$, then $b(t, \beta)\nearrow b_1\ldots b_m^-\de(\beta)$ as $t\nearrow t_0$.
\end{itemize}
\end{lemma}

  Recall that the survivor set $K_\beta(t)$ consists of all $x\in[0,1)$ whose orbit $\{T^n_\beta(x): n\ge 0\}$ avoids the hole $(0,t)$. To describe the dimension of $K_\beta(t)$ we introduce the topological entropy of a symbolic set. For a subset $X\subset\set{0,1}^\N$ its \emph{topological entropy} $h_{top}(X)$ is defined by
\[
h_{top}(X):=\liminf_{n\to\f}\frac{\log\#B_n(X)}{n},
\]
where $\#B_n(X)$ denotes the number of all length $n$ words occurring in sequences of $X$.
The following result for the Hausdorff dimension of $K_\beta(t)$ can be essentially deduced from Raith \cite{Raith-94}  (see also \cite{Kalle-Kong-Langeveld-Li-18}).

\begin{lemma} \label{lem:dim-survivorset}
Given $\beta\in(1,2]$ and $t\in[0,1)$, the Hausdorff dimension of $K_\beta(t)$ is given by
\[
\dim_H K_\beta(t) =\frac{h_{top}(\K_\beta(t))}{\log\beta},
\]
where
\[
\K_\beta(t)=\set{(d_i)\in\set{0,1}^\N: b(t,\beta)\lle\si^n((d_i))\prec \de(\beta)~\forall n\ge 0}.
\]
\end{lemma}

In order to determine the critical value $\tau(\beta)$ for $\beta$ inside any Farey interval $J^\s$, we first need to develop some properties of the substitution operator $\bullet$ from \eqref{eq:substition-operator}. We do this in the next section.

\section{Substitution of Lyndon words} \label{sec:substitution}

In this section we give an equivalent definition of the substitution operator in $\LL$ introduced in (\ref{eq:substition-operator}), and prove that $\LL$ forms a semi-group under this substitution operator. This will play a {crucial role in the rest of the paper.}

\subsection{An equivalent definition of the substitution}
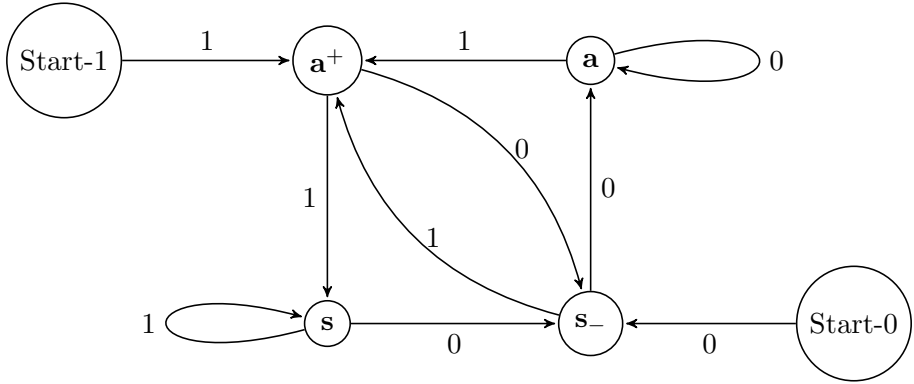
\begin{figure}[h!]
  \centering
  \begin{tikzpicture}[->,>=stealth',shorten >=1pt,auto,node distance=3.5cm, semithick,scale=5]

  \tikzstyle{every state}=[minimum size=0pt,fill=none,draw=black,text=black]

  \node[state] (A)                    { $\s$};
  \node[state]         (B) [ right of=A] {${\s^-}$ };
  \node[state]         (C) [ above of=A] {$\a^+$};
  \node[state] (E)[left of=C]{Start-$1$};
  \node[state](D)[right of=C]{$\a$};
  \node[state](F)[right of=B]{Start-$0$};

  \path[->,every loop/.style={min distance=0mm, looseness=40}]
  (E) edge[->] node{$1$} (C)
  (C) edge[->,left] node{$1$} (A)
  (C) edge[bend left,->,right] node{$0$} (B)
(D) edge[loop right,->] node{$0$} (D)
(D) edge[->,above] node{$1$} (C)
(A) edge[loop left,->] node{$1$} (A)
(A)edge[->,below] node{$0$} (B)
(B) edge[bend left,->,right] node{$1$} (C)
(B) edge[->,right] node{$0$} (D)
(F) edge[->] node{$0$} (B)
;
\end{tikzpicture}

\caption{The directed graph $G=(V, E)$ with the edge labels from $\set{0, 1}$ and vertex labels from $\set{{\s^-}, \s, \a, \a^+}$, where $\s\in\LL$ and $\a=\L(\s)$. 
}
\label{fig1}
\end{figure}

Given a Lyndon word $\s\in\LL$ with $\a=\L(\s)$, we construct a directed graph $G=(V, E)$ as in Figure \ref{fig1}. The directed graph $G$ has two starting vertices  `Start-$0$' and `Start-$1$'.  The directed edges in the graph $G$ take  labels from $\set{0,1}$, and the vertices in $G$ take labels from $\set{{\s^-},\s,\a,\a^+}$. Denote by $\mathcal L^E$ the edge labeling and by $\mathcal L^V$ the vertex labeling. Then for each directed edge $e\in E(G)$ we have $\mathcal L^E(e)\in\set{0,1}$, and for each vertex $v\in V(G)$ we have $\mathcal L^V(v)\in\set{{\s^-},\s,\a,\a^+}$. The labeling maps $\mathcal L^E$ and $\mathcal L^V$ naturally induce the maps on the infinite edge paths and infinite vertex paths in $G$, respectively. For example, for an infinite edge path $e_1e_2\ldots$ we have
\[
\mathcal L^E(e_1e_2\ldots)=\mathcal L^E(e_1)\mathcal L^E(e_2)\ldots~\in\set{0,1}^\N.
\]
Here we call $e_1e_2\ldots$  an \emph{infinite edge path} in $G$, if the initial vertex of $e_1$ is one of the starting vertices and  for any $i\ge 1$ the \emph{terminal vertex} $t(e_i)$ equals the \emph{initial vertex} $i(e_{i+1})$.
Similarly, for an infinite vertex path $v=v_1v_2\ldots$ we have
\[
\mathcal L^V(v_1v_2\ldots)=\mathcal L^V(v_1)\mathcal L^V(v_2)\ldots~\in\set{{\s^-}, \s,\a,\a^+}^\N,
\]
where we call $v_1v_2\ldots$ an \emph{infinite vertex path} in $G$, if $v_1$ is one of the starting vertices and for any $i\ge 1$ there exists a directed edge $e\in E(G)$ such that $i(e)=v_i$ and $t(e)=v_{i+1}$.

Let $X_E=X_E(G)$ be the \emph{edge shift} consisting of all {labelings} of infinite edge paths in $G$, i.e.,
\[
X_E:=\set{\mathcal L^E(e_1e_2\ldots): e_1e_2\ldots\textrm{ is an infinite edge path in }G}.
\]
One can verify easily that ${X_E=\set{0,1}^\N}$.
 Also, let $X_V=X_V(G)$ be the \emph{vertex shift} which consists of all {labelings} of infinite vertex paths in $G$, i.e.,
\[
X_V:=\set{\mathcal L^V(v_1v_2\ldots): v_1v_2\ldots\textrm{ is an infinite vertex path in }G}.
\]
Then any sequence in $X_V$ is an infinite concatenation of words from $\set{{\s^-}, \s,\a,\a^+}$.  Observe that the edge shift $X_E$ is \emph{right-resolving}, which means that out-going edges from the same vertex have different labels (cf.~\cite{Lind_Marcus_1995}). Moreover, different vertices have different labels. So for each $(d_i)\in X_E$ there is a unique infinite edge path $e_1e_2\ldots$ in $G$ such that
$d_1d_2\ldots=\mathcal L^E(e_1e_2\ldots).$ 

\begin{definition} \label{def:substitution}
The substitution map $\Phi_\s$ from $X_E$ to $X_V$ is defined by
\begin{equation*} \label{eq:block-map}
\Phi_\s: X_E\to X_V;\quad \mathcal L^E(e_1e_2\ldots)\mapsto \mathcal L^V(t(e_1)t(e_2)\ldots),
\end{equation*}
where $t(e_i)$ denotes the terminal vertex of the directed edge  $e_i$.
\end{definition}

We can extend the substitution map $\Phi_\s$ to a map from $B_*(X_E)$ to $B_*(X_V)$ by
\begin{equation} \label{eq:substition-1}
\Phi_\s: B_*(X_E)\to B_*(X_V);\quad \mathcal L^E(e_1\ldots e_n)\mapsto \mathcal L^V(t(e_1)\ldots t(e_n)),
\end{equation}
where
 $B_*(X_E)$ consists of all labelings of finite edge paths in $G$, and $B_*(X_V)$ consists of all labelings of finite vertex paths in $G$. So, by (\ref{eq:substition-operator}) and (\ref{eq:substition-1}) it follows that for any two words $\s\in\LL$ and $\r\in\set{0,1}^*$ we have
 \begin{equation} \label{eq:substitution}
 \s\bullet\r=\Phi_\s(\r).
 \end{equation}

\begin{example} \label{ex:1}
  Let $\s=01$ and $\r=001011$. Then $\s\in\F$ and $\r\in\LL\setminus\F$. Furthermore, ${\s^-}=00, \a=\L(\s)=10, \a^+=11$. So by the definition of $\Phi_\s$ it follows that
  \begin{align*}
   \Phi_\s(\r)&=\Phi_\s( 001011 )={\s^-}\a\a^+{\s^-}\a^+\s=001011001101, \\
  \Phi_\s(\r^-)&= \Phi_\s( 001010)={\s^-}\a\a^+{\s^-}\a^+{\s^-}=001011001100.
  \end{align*}
Observe that $\Phi_\s(\r^-)=\Phi_\s(\r)^-$.
  By Definition \ref{def:lyndon words} one can check that $\Phi_\s(\r)\in\LL$. Furthermore,
 \[ \Phi_\s(\L(\r))=\Phi_\s( 110010) =\a^+\s{\s^-}\a\a^+{\s^-}=110100101100=\L(\Phi_\s(\r)),\]
 and
 \[
  \Phi_\s(\r^\f)= \Phi_\s((001011)^\f)=({\s^-}\a\a^+{\s^-}\a^+\s)^\f=(001011001101)^\f=\Phi_\s(\r)^\f. \]
\end{example}

\subsection{Properties of the substitution}
Motivated by Examples \ref{ex:substitution} and \ref{ex:1} we study the properties of the substitution $\Phi_\s$. We will show that $\LL$ forms a semi-group under the substitution operator defined in Definition \ref{def:substitution}.
 First we prove the monotonicity of $\Phi_\s$.

\begin{lemma} \label{lem:monotonicity-Phi}
Let $\s\in\LL$. Then the map $\Phi_\s$ is strictly increasing in $X_E=\set{0,1}^\N$.
\end{lemma}

\begin{proof}
Let $(d_i)$ and $(d_i')$ be two sequences in $X_E$, and let $(e_i), (e_i')$ be their corresponding edge paths; thus, $(d_i)=\mathcal L^E((e_i))$ and $(d_i')=\mathcal L^E((e_i'))$. Suppose $(d_i)\prec(d_i')$. Then there exists $k\in\N$ such that $d_1\ldots d_{k-1}=d_1'\ldots d_{k-1}'$ and $d_k<d_k'$. If $k=1$, then $d_1=0$ and $d_1'=1$. So, $\mathcal L^V(t(e_1))={\s^-}$ and $\mathcal L^V(t(e_1'))=\a^+$. By Definition \ref{def:substitution} it follows that  $\Phi_\s((d_i))\prec\Phi_\s((d_i'))$.

 If $k>1$, then $e_1\ldots e_{k-1}=e_1'\ldots e_{k-1}'$, which implies that the initial vertices of $e_k$ and  $e_k'$ coincide.  Since $d_k<d_{k}'$, by the definition of $\mathcal L^V$ it follows that (see Figure \ref{fig1})
\[
\mathcal L^V(t(e_k))\prec\mathcal L^V(t(e_k')).
\]
By Definition \ref{def:substitution} we also have $\Phi_\s((d_i))\prec\Phi_\s((d_i'))$. This completes the proof.
\end{proof}

\begin{lemma} \label{lem:mon-1}
Let $\s\in\LL$. Then for any  word $\d=d_1\ldots d_k\in B_*(X_E)$ with $k\ge 2$, we have
\[
\left\{
\begin{array}{lll}
\Phi_\s(\d^-)=\Phi_\s(\d)^-&\textrm{if}& d_k=1,\\
\Phi_\s(\d^+)=\Phi_\s(\d)^+&\textrm{if}& d_k=0.
\end{array}\right.
\]
\end{lemma}

\begin{proof}
Since $\d=d_1\ldots d_k\in B_*(X_E)$, there exists a unique finite edge path $e_1\ldots e_k$ such that $\mathcal L^E(e_1\ldots e_k)=\d$. If $d_k=1$, then $\d^-=d_1\ldots d_{k-1}0$ can be represented by a unique finite edge path $e_1'\ldots e_k'$ with  $e_1'\ldots e_{k-1}'=e_1\ldots e_{k-1}$.  By the definition of $\mathcal L^V$ it follows that $\mathcal L^V(t(e_k'))=\mathcal L^V(t(e_k))^-$. Therefore, by Definition \ref{def:substitution} it follows that
\begin{align*}
\Phi_\s(\d^-)=\Phi_\s(\mathcal L^E(e_1'\ldots e_k'))
&= \Phi_\s(\mathcal L^E(e_1\ldots e_{k-1}e_k'))\\
&=\mathcal L^V(t(e_1)\ldots t(e_{k-1})t(e_k')) \\
&=\mathcal L^V(t(e_1)\ldots t(e_k))^-=\Phi_\s(\d)^-.
\end{align*}
This proves the first equality of the lemma. The second equality follows analogously.
\end{proof}



Recall the operator $\bullet$ from (\ref{eq:substitution}). In the following we prove Proposition \ref{main:substitution} by showing that $\LL$ is closed under $\bullet$ and that $\bullet$ is associative. The proof will be split into a sequence of lemmas. First we prove that $\LL$ is closed under $\bullet$.

\begin{lemma}  \label{lem:substitution-Lyndon}
  For any $\s, \r\in\LL$ we have $\s\bullet\r\in\LL$.
\end{lemma}

\begin{proof}
  Let $\s=s_1\ldots s_m\in\LL$ and $\a=\L(\s)$. Then there exists $j\in\set{1,\ldots,m-1}$ such that
  \begin{equation}\label{eq:lyndon-1}
  \a=s_{j+1}\ldots s_m s_1\ldots s_j.
  \end{equation}
  Let $\r=r_1\ldots r_\ell\in\LL$. Then we can write $\s\bullet\r=\Phi_\s(\r)=b_1\ldots b_{m\ell}$. Furthermore, there exists a finite edge path $e_1\ldots e_\ell$ representing $\r$ such that
  \[\Phi_\s(\r)=\mathcal L^V(t(e_1) \ldots t(e_\ell))=:\b_1\ldots\b_\ell,\]
  where each block $\b_i\in\set{{\s^-}, \s, \a,\a^+}$.
 Note that $\b_1={\s^-}$ since the block $\r$ begins with $r_1=0$.
  By Definition \ref{def:lyndon words} it  suffices to prove
  \begin{equation}
    \label{eq:lyndon-2}
    b_{i+1}\ldots b_{m\ell}\succ b_1\ldots b_{m\ell-i}\quad\textrm{for any }0<i<m\ell.
  \end{equation}
  We split the proof of (\ref{eq:lyndon-2}) into two cases.

  Case I. $i=km$ for some $k\in\set{1, 2,\ldots, \ell-1}$. Then $b_{i+1}\ldots b_{m\ell}=\b_{k+1}\ldots \b_{\ell}$. Since $\r$ is a Lyndon word, we have
  $r_{k+1}\ldots r_\ell\succ r_1\ldots r_{\ell-k}$. So, (\ref{eq:lyndon-2}) follows directly by Lemma \ref{lem:monotonicity-Phi}.

Case II. $i=km+p$ for some $k\in\set{0,1,\ldots, \ell-1}$ and $p\in\set{1,\ldots, m-1}$. Then $b_{i+1}\ldots b_{m\ell}=b_{i+1}\ldots b_{i+m-p}\b_{k+2}\ldots \b_{\ell}$. In the following we prove (\ref{eq:lyndon-2}) by considering  the four possible choices  of $\b_{k+1}\in\set{{\s^-}, \s, \a,\a^+}$. If $\b_{k+1}=\s$, then by using that $\s\in\LL$ we conclude that
\[b_{i+1}\ldots b_{i+m-p}=s_{p+1}\ldots s_m\succ s_1\ldots s_{m-p}=b_1\ldots b_{m-p},\] proving (\ref{eq:lyndon-2}).
Similarly, if $\b_{k+1}=\a^+$, then by  (\ref{eq:lyndon-1}) one can also prove that $b_{i+1}\ldots b_{i+m-p}\succ b_1\ldots b_{m-p}$. Now we assume $\b_{k+1}={\s^-}$. Then by using $\s\in\LL$ it follows that
\begin{equation} \label{eq:lyndon-3}
b_{i+1}\ldots b_{i+m-p}=s_{p+1}\ldots s_m^-\lge s_1\ldots s_{m-p}=b_1\ldots b_{m-p}.
\end{equation}
Observe that the word ${\s^-}$ can only be followed by $\a$ or $\a^+$ in   $G$ (see Figure \ref{fig1}). So $\b_{k+2}\in\set{\a, \a^+}$. Since $\a=\L(\s)$, we obtain that
\begin{equation}\label{eq:lyndon-4}
b_{i+m-p+1}\ldots b_{i+m}=a_1\ldots a_p\lge s_{m-p+1}\ldots s_m\succ s_{m-p+1}\ldots s_m^-=b_{m-p+1}\ldots b_m.
\end{equation}
Thus, by (\ref{eq:lyndon-3}) and (\ref{eq:lyndon-4}) we conclude that $b_{i+1}\ldots b_{i+m}\succ b_{1}\ldots b_m$, proving (\ref{eq:lyndon-2}). Finally, suppose $\b_{k+1}=\a$. Note that the word $\a$ can only be followed by $\a$ or $\a^+$ in $G$. Then by (\ref{eq:lyndon-1}) and using $\s\in\LL$ we have
\[
b_{i+1}\ldots b_{i+m}\lge s_1\ldots s_m\succ b_1\ldots b_m.
\]
This completes the proof.
\end{proof}

{
Say a finite or infinite sequence of words $\b_1,\dots, \b_n$ or $\b_1,\b_2,\dots$ is {\em connectible} if for each $i$, the last digit of $\b_i$ differs from the first digit of $\b_{i+1}$. Thus, for instance, the sequence $1101, 00111$ is connectible whereas the sequence $11010, 0111$ is not.

\begin{lemma} \label{lem:connectible}\mbox{}
\begin{enumerate}[{\rm (i)}]
\item Let $\b_1,\b_2,\dots$ be a (finite or infinite) connectible sequence of words. Then for any $\s\in\LL$,
\[
\Phi_\s(\b_1\b_2\dots)=\Phi_\s(\b_1)\Phi_\s(\b_2)\dots.
\]
\item Let $\s, \r\in\LL$. Then $\Phi_\s(\r^\f)=\Phi_\s(\r)^\f$ and $\Phi_\s(\L(\r)^\f)=\Phi_\s(\L(\r))^\f$.
\end{enumerate}
\end{lemma}

\begin{proof}
To prove (i) it suffices to show that if $\b_1,\b_2$ is a connectible sequence, then $\Phi_\s(\b_1\b_2)=\Phi_\s(\b_1)\Phi_\s(\b_2)$; the statement then extends to arbitrary connectible sequences by induction.

Without loss of generality, by the symmetry of the edge-labels in Figure \ref{fig1}, we may assume that $\b_1$ ends in the digit 0 and $\b_2$ begins with the digit 1. But note that in the directed graph in Figure \ref{fig1}, if we travel along an edge labeled 0 followed by an edge labeled 1, we always end up at the vertex labeled $\a^+$, which is also the first vertex visited after traveling along an edge labeled 1 from the ``Start-1" vertex. Thus, $\Phi_\s(\b_1\b_2)=\Phi_\s(\b_1)\Phi_\s(\b_2)$.

Statement (ii) follows from (i) since $\r$ begins with digit 0 and ends with digit 1, so $\r$ is connectible to itself; and similarly, $\L(\r)$ begins with digit 1 and ends with digit 0, so $\L(\r)$ is connectible to itself.
\end{proof}
}

To prove that $\bullet$ is associative we need the following result, which says that the two operators $\bullet$ and $\L$ commute.

\begin{lemma}  \label{lem:lyndon-commutative}
  For any $\s, \r\in\LL$ we have $\L(\s\bullet\r)=\s\bullet\L(\r)$.
\end{lemma}

\begin{proof}
 The proof is similar to that of Lemma \ref{lem:substitution-Lyndon}. Let $\r=r_1\ldots r_\ell\in\LL$. First we {show} that $\s\bullet\L(\r)$ is a cyclic permutation of $\s\bullet\r$.
    Note that $\L(\r)=r_{j+1}\ldots r_\ell r_1\ldots r_j$ for some $1<j<\ell$. Then $r_j=0$ and $r_{j+1}=1$, {so Lemma \ref{lem:connectible} (i) implies} that
  \begin{equation}  \label{eq:a11-1}
     \s\bullet\r=\Phi_\s(r_1\ldots r_\ell)=\Phi_\s(r_1\ldots r_j)\Phi_\s(r_{j+1}\ldots r_\ell).
  \end{equation}
  On the other hand, since $\r\in\LL$ we have $r_\ell=1$ and $r_1=0$, so by {Lemma \ref{lem:connectible} (i)} we obtain that
  \[
  \s\bullet\L(\r)=\Phi_\s(r_{j+1}\ldots r_\ell r_1\ldots r_j)=\Phi_\s(r_{j+1}\ldots r_\ell)\Phi_\s(r_1\ldots r_j).
  \]
  This, together with (\ref{eq:a11-1}), proves that $\s\bullet\L(\r)$ is indeed a cyclic permutation of $\s\bullet\r$. {It remains} to prove that $\s\bullet\L(\r)$ is the lexicographically largest cyclic permutation of itself.
	
  Write $\s=s_1\ldots s_m\in\LL$ with $\a=\L(\s)=a_1\ldots a_m$, and write $\L(\r)=c_1\ldots c_\ell$. Then by Lemma \ref{lem:lyndon-equivalence} it follows that
  \begin{equation}
  \begin{split}\label{eq:sep17-1}
      a_{i+1}\ldots a_m&\prec a_1\ldots a_{m-i}\quad\forall ~0<i<m;\\
      c_{i+1}\ldots c_\ell&\prec c_1\ldots c_{\ell-i}\quad \forall ~0<i<\ell.
  \end{split}
  \end{equation}
  Write $\s\bullet\L(\r)=\b_1\ldots \b_\ell=b_1\ldots b_{m\ell}$,
  where each $\b_i\in\set{{\s^-}, \s, \a,\a^+}$. Then it suffices to prove that
\begin{equation}  \label{eq:sep17-2}
  b_{i+1}\ldots b_{m\ell}\prec b_1\ldots b_{m\ell-i}\quad \forall~0<i<m\ell.
\end{equation}
  Since {$\L(\r)$} has a prefix $c_1=1$, we see that $b_1\ldots b_m=\b_1=\a^+$.
So, by using  (\ref{eq:sep17-1}) and the same argument as in the proof of Lemma \ref{lem:substitution-Lyndon} we can prove (\ref{eq:sep17-2}).
\end{proof}

The next lemma will be used in the proof of Lemma \ref{lem:geometry-basic-interval-2} and Proposition \ref{prop:cric-exception-1}.

\begin{lemma}  \label{lem:cri-exception-1}
  Let $\s\in\LL$ and {take} two sequences $(c_i), (d_i)\in\set{0,1}^\N$.
  \begin{enumerate}
    \item [{\rm(i)}] If $d_1=1$, then
    \[
    \si^n((c_i))\prec (d_i)~\forall n\ge 0\quad\Longrightarrow\quad \si^n(\Phi_\s((c_i)))\prec \Phi_\s((d_i))~\forall n\ge 0.
    \]
    \item [{\rm(ii)}] If $d_1=0$, then
    \[
    \si^n((c_i))\succ (d_i)~\forall n\ge 0\quad\Longrightarrow\quad \si^n(\Phi_\s((c_i)))\succ \Phi_\s((d_i))~\forall n\ge 0.
    \]
  \end{enumerate}
 \end{lemma}


 \begin{proof}
 (i) Suppose $d_1=1$ and $\si^n((c_i))\prec (d_i)$ for all $n\ge 0$. Then $\Phi_{\s}((d_i))$ begins with $\L(\s)^+$. If $n\equiv 0\pmod {|\s|}$, then by Lemma \ref{lem:monotonicity-Phi} it follows that $\si^n(\Phi_\s((c_i)))\prec \Phi_\s((d_i))$. If $n\ne 0\pmod {|\s|}$, then by using $\Phi_\s(d_1)=\L(\s)^+$ and the same argument as  in the proof of Lemma \ref{lem:lyndon-commutative} one can verify that $\si^n(\Phi_\s((c_i)))\prec \Phi_\s((d_i))$. The proof of (ii) is similar.
 \end{proof}

Finally, we show that $\bullet$ is associative.

\begin{lemma}  \label{lem:lyndon-associative}
  For any three words $\r, \s, \t\in\LL$ we have $(\r\bullet\s)\bullet\t=\r\bullet(\s\bullet\t)$.
\end{lemma}

\begin{proof}
Let $\r=r_1\ldots r_m, \s=s_1\ldots s_n$ and $\t=t_1\ldots t_\ell$. 
Then we can write $(\r\bullet\s)\bullet\t$ as
\begin{equation} \label{eq:26-1}
(\r\bullet\s)\bullet\t=\mathbf B_1\mathbf B_2\ldots \mathbf B_\ell,
\end{equation}
where each $\mathbf B_i\in\set{(\r\bullet\s)^-, \r\bullet\s, \L(\r\bullet\s), \L(\r\bullet\s)^+}$. Since $t_1=0$, we have $\mathbf B_1=(\r\bullet\s)^-$. Furthermore,  by the definition of $\Phi_{\r\bullet\s}$ it follows that for $1<i\le \ell$,
 \begin{equation}   \label{eq:kk-1}
   \mathbf B_i=\left\{
   \begin{array}
     {lll}
     \L(\r\bullet\s)&\textrm{if}&  t_{i-1}t_i=00,\\
     \L(\r\bullet\s)^+&\textrm{if}& t_{i-1}t_i=01,\\
     (\r\bullet\s)^-&\textrm{if}& t_{i-1}t_i=10,\\
     \r\bullet\s&\textrm{if}& t_{i-1}t_i=11.
   \end{array}\right.
 \end{equation}

{Similarly, we can write}
\begin{equation*}
\s\bullet\t=\b_1\b_2\ldots\b_\ell,
\end{equation*}
 where each $\b_i\in\set{{\s^-},  \s, \L(\s), \L(\s)^+}$, {and it follows from} the definition of $\Phi_\s$ that
 \begin{equation} \label{eq:may8-1}
   \b_i=\left\{
   \begin{array}
     {lll}
     \L(\s)&\textrm{if}&  t_{i-1}t_i=00,\\
     \L(\s)^+&\textrm{if}& t_{i-1}t_i=01,\\
     {\s^-}&\textrm{if}& t_{i-1}t_i=10,\\
     \s&\textrm{if}& t_{i-1}t_i=11,
   \end{array}\right.
 \end{equation}
for $1<i\le \ell$. {Comparing \eqref{eq:kk-1} and \eqref{eq:may8-1} and using Lemmas \ref{lem:mon-1} and \ref{lem:lyndon-commutative} it follows that
\[\Phi_\r(\b_i)=\mathbf{B}_i\quad \forall i\ge 1.\] Moreover, the sequence $\b_1,\b_2,\dots,\b_\ell$ is connectible because $\b_1\b_2\dots\b_\ell=\Phi_\s(\t)$ arises from a walk along the directed graph in Figure \ref{fig1}. Hence, Lemma \ref{lem:connectible} yields
\[
\r\bullet(\s\bullet\t)=\Phi_\r(\b_1\b_2\ldots\b_\ell)=\Phi_\r(\b_1)\Phi_\r(\b_2)\dots\Phi_\r(\b_\ell)=\mathbf{B}_1\mathbf{B}_2\dots\mathbf{B}_\ell=(\r\bullet\s)\bullet \t,
\]
as desired.}
\end{proof}

\begin{proof}[Proof of Proposition \ref{main:substitution}]
  The proposition follows by Lemmas \ref{lem:substitution-Lyndon} and \ref{lem:lyndon-associative} and Example \ref{ex:substitution}, which shows that $\bullet$ is not commutative.
\end{proof}

\section{Critical values in a basic interval} \label{sec:basic-intervals}

In this section we will prove Theorem \ref{main:critical-basic-intervals}. 
Recall from (\ref{eq:set-substitution-Farey-words}) that $\La$ consists of all words $\cs$ of the form
\begin{equation*} \label{eq:comosition-farey}
\cs=\s_1\bullet\s_2\bullet\cdots\bullet\s_k,\quad k\in\N,
\end{equation*}
where each $\s_i\in\F$.  By Proposition \ref{main:substitution} it follows that $\Lambda\subset\LL$, and each $\cs\in\Lambda$ can be uniquely represented in the above form.
Take $\cs\in\Lambda$. As in Definition \ref{def:basic-intervals} we let $I^\cs:=[\beta_\ell^\cs, \beta_*^\cs]$ be the basic interval generated by $\cs$. Then by Lemmas \ref{lem:quasi-greedy expansion} and \ref{lem:lyndon-equivalence} it follows that
\begin{equation} \label{eq:basic-interval}
\de(\beta_\ell^\cs)=\L(\cs)^\f\quad\textrm{and}\quad \de(\beta_*^\cs)=\L(\cs)^+\cs^-\L(\cs)^\f.
\end{equation}

To prove Theorem \ref{main:critical-basic-intervals} we first prove the following proposition, which provides one of the key tools in this paper and will be used again in Section \ref{sec:cri-bifurcation-sets}.

\begin{proposition} \label{prop:countable}
For any $\cs\in\La$ the set
\begin{equation} \label{eq:8-2}
  \Ga(\cs):=\set{(x_i): \cs^\f \lle \si^n((x_i))\lle \L(\cs)^\f~\forall n\ge 0}
\end{equation}
is countable.
\end{proposition}

We point out that the specific form of $\cs$ is essential in this proposition: it is not enough to merely assume that $\cs$ is a Lyndon word. For instance, take $\cs=0010111\in \Omega_L^*$. Then $\L(\cs)=1110010$, and it is easy to see that $\Ga(\cs)\supset \{10,110\}^\N$.

For $\cs=\s\in\F$, Proposition \ref{prop:countable} follows from the following stronger result, proved in \cite[{Proposition 4.4}]{Kalle-Kong-Langeveld-Li-18}.

\begin{lemma} \label{lem:finite-symmetric}
  For any $\s=s_1\ldots s_m\in\F$, the set
  \[
  \Ga(\s):=\set{(x_i): {\s^\f}\lle \si^n((x_i))\lle \L(\s)^\f~\forall n\ge 0}
  \]
  consists of exactly $m$ different elements.
\end{lemma}

In order to reduce the technicalities in the proof of Proposition \ref{prop:countable}, we extend the definition from Lemma \ref{lem:Farey-property} (iii) and define the {\em conjugate} of any word $\cs\in\La$ by
\[
\varphi(\cs):=\overline{\L(\cs)}.
\]

\begin{lemma} \label{lem:general-conjugate}
The function $\varphi:\La\to\{0,1\}^*;\ \cs\mapsto \varphi(\cs)$ is a semigroup automorphism on $(\La,\bullet)$. That is, $\varphi$ maps $\La$ bijectively onto itself, and
\begin{equation} \label{eq:distribute-phi}
\varphi(\s_1\bullet\cdots\bullet\s_k)=\varphi(\s_1)\bullet\cdots\bullet\varphi(\s_k) \qquad \forall \s_1,\dots,\s_k\in\F.
\end{equation}
Furthermore, $\varphi$ is its own inverse:
\begin{equation} \label{eq:phi-inverse}
\varphi(\varphi(\cs))=\cs \qquad\forall \cs\in\La.
\end{equation}
\end{lemma}

\begin{proof}
We prove \eqref{eq:distribute-phi} and \eqref{eq:phi-inverse} simultaneously by induction on the degree $k$ of $\cs=\s_1\bullet\cdots\bullet\s_k$. For $k=1$, \eqref{eq:distribute-phi} is trivial and \eqref{eq:phi-inverse} follows from Lemma \ref{lem:Farey-property} (iii). Now suppose \eqref{eq:distribute-phi} and \eqref{eq:phi-inverse} both hold for any word $\cs=\s_1\bullet\cdots\bullet\s_k$ of degree $k$, and consider $\cs\bullet \r$ with $\r\in\F$. Set $\widetilde{\cs}:=\varphi(\cs)$. We claim first that, for any word $\t\in\{0,1\}^*$,
\begin{equation} \label{eq:reflection-relation}
\Phi_\cs(\overline{\t})=\overline{\Phi_{\widetilde{\cs}}(\t)}.
\end{equation}
The expression on the right is well defined since, by \eqref{eq:distribute-phi}, $\widetilde{\cs}=\varphi(\s_1)\bullet\cdots\bullet\varphi(\s_k)\in\Omega_L^*$.

Write $\ca:=\L(\cs)$ and $\widetilde{\ca}:=\L(\widetilde{\cs})$. By \eqref{eq:phi-inverse}, $\varphi(\widetilde{\cs})=\cs$, so we have
\begin{equation}\label{eq:conjugate-S}
\overline{\ca}=\widetilde{\cs}\quad \textrm{and} \quad\widetilde{\ca}=\overline{\cs}.
 \end{equation}Now note the rotational skew-symmetry in the edge labels of the directed graph in Figure \ref{fig1}: The edge path corresponding to the word $\overline{\t}$ is just the $180^\circ$ rotation about the center of the figure of the edge path corresponding to $\t$. On the other hand, replacing the vertex labels $\cs,\cs^-,\ca$ and $\ca^+$ by $\widetilde{\cs}=\overline{\ca}$, $\widetilde{\cs}^-=\overline{\ca^+}$, $\widetilde{\ca}=\overline{\cs}$ and $\widetilde{\ca}^+=\overline{\cs^-}$ respectively and rotating the whole graph by $180^\circ$, we get the original graph back except that all the vertex labels and edge labels are reflected. This implies \eqref{eq:reflection-relation}.

Now we can apply \eqref{eq:reflection-relation} to $\t=\overline{\L(\r)}$ and obtain:
\begin{equation} \label{eq:phi-distribute-Sr}
\varphi(\cs)\bullet \varphi(\r)=\widetilde{\cs}\bullet\varphi(\r)=\Phi_{\widetilde{\cs}}(\overline{\L(\r)})=\overline{\Phi_\cs(\L(\r))}
=\overline{\cs\bullet\L(\r)}=\overline{\L(\cs\bullet \r)}=\varphi(\cs\bullet \r).
\end{equation}
Since $\r\in\F$ was arbitrary, the induction hypothesis \eqref{eq:distribute-phi} and \eqref{eq:phi-distribute-Sr} give
\begin{equation} \label{eq:phi-distribute-plus-1}
\varphi(\s_1\bullet\cdots\bullet\s_k\bullet \s_{k+1})=\varphi(\s_1)\bullet\cdots\bullet\varphi(\s_k)\bullet\varphi(\s_{k+1}) \qquad \forall \s_1,\dots,\s_{k+1}\in\F.
\end{equation}
Thus, \eqref{eq:distribute-phi} holds for $k+1$ in place of $k$. Next, by Lemma \ref{lem:Farey-property} (iii), $\varphi(\s_i)\in\F$ and $\varphi(\varphi(\s_i))=\s_i$ for each $i$, so applying \eqref{eq:phi-distribute-plus-1} with $\varphi(\s_i)$ in place of $\s_i$ for each $i$, we conclude that $\varphi(\varphi(\cs'))=\cs'$ for every $\cs'\in\La$ of degree $k+1$ also.

Thus, we have proved \eqref{eq:distribute-phi} and \eqref{eq:phi-inverse} by induction. Now the remaining statements of the lemma follow immediately: by \eqref{eq:distribute-phi}, Lemma \ref{lem:Farey-property} (iii) and Proposition \ref{main:substitution} it follows that $\varphi(\cs)\in\La$ for every $\cs\in\La$, whereas \eqref{eq:phi-inverse} implies that $\varphi: \La\to \La$ is bijective. Therefore, $\varphi$ is an automorphism of $(\La,\bullet)$.
\end{proof}

Define
\[
\overline{\Ga(\cs)}:=\big\{\overline{(x_i)}: (x_i)\in\Ga(\cs)\big\}, \qquad \cs\in\La.
\]
It is clear that $\overline{\Ga(\cs)}$ has the same cardinality as $\Ga(\cs)$. Observe also by (\ref{eq:conjugate-S}) that
\begin{align}
\begin{split}
\overline{\Ga(\cs)}&=\big\{(y_i): \cs^\f\lle \sigma^n(\overline{(y_i)})\lle \L(\cs)^\f\ \forall n\geq 0\big\}\\
&=\big\{(y_i): \overline{\cs}^\f\lge \sigma^n((y_i))\lge \overline{\L(\cs)}^\f\ \forall n\geq 0\big\}\\
&=\big\{(y_i): \L(\varphi(\cs))^\f \lge \sigma^n((y_i))\lge \varphi(\cs)^\f\ \forall n\geq 0\big\}\\
&=\Ga(\varphi(\cs)).
\end{split}
\label{eq:reflection-and-conjugate}
\end{align}

\begin{proof}[Proof of Proposition \ref{prop:countable}]
For $\cs=\s\in\F$ the proposition follows from Lemma \ref{lem:finite-symmetric}. So it suffices to prove that if $\Ga(\cs)$ is countable for an $\cs\in\La$, then $\Ga(\cs\bullet\r)$ is also countable for any $\r\in\F$.

Fix $\cs\in\La$ with $\Ga(\cs)$ countable; fix $\r\in\F$, and note that $\r$ begins with $0$ and $\L(\r)$ begins with $1$. {Therefore} $\cs\bullet\r$ begins with $\cs^-$ and $\L(\cs\bullet\r)=\cs\bullet\L(\r)$ begins with $\L(\cs)^+$. So
\begin{equation} \label{eq:8-1}
(\cs\bullet\r)^\f\prec \cs^\f\quad\textrm{and}\quad \L(\cs\bullet\r)^\f\succ \L(\cs)^\f.
\end{equation}
By (\ref{eq:8-2}) and (\ref{eq:8-1}) it follows that
\[
\Ga(\cs)\subseteq\set{(x_i): (\cs\bullet\r)^\f\lle\si^n((x_i))\lle \L(\cs\bullet\r)^\f~\forall n\ge 0}=\Ga(\cs\bullet\r).
\]
Since $\Ga(\cs)$ is countable, it suffices to prove that the difference set $\Ga(\cs\bullet\r)\setminus \Ga(\cs)$ is countable.
By Proposition \ref{prop:Farey-new-form} the word $\r$ must be of one of the following four types:
\begin{itemize}
\item[{\rm(I)}]  $\r=01^p$ for some $p\in\N$;
\item[{\rm(II)}] $\r=0^p 1$ for some $p\in\N$;
\item[{\rm(III)}]
{$\r=01^p01^{p+t_1}\ldots 01^{p+t_N}01^{p+1}$ for some $p\in\N$ and  $0t_1\ldots t_N 1\in\F$}.
\item[{\rm (IV)}] {$\r=0^{p+1}10^{p+t_1}1\ldots 0^{p+t_N}10^p1$ for some $p\in\N$ and $0t_1\ldots t_N 1\in\F$};
\end{itemize}
Since the words in (II) and (IV) are the conjugates of the words in (I) and (III), respectively, it suffices by Lemma \ref{lem:general-conjugate} and the relationship \eqref{eq:reflection-and-conjugate} to consider cases (I) and (III). Let $\ca:=\L(\cs)$.

\medskip
{\bf Case I.} $\r=01^p$ for some $p\in\N$. Note that $\cs\bullet\r=\Phi_\cs(01^p)=\cs^-\ca^+\cs^{p-1}$ and $\L(\cs\bullet\r)=\cs\bullet\L(\r)=\Phi_\cs(1^p 0)=\ca^+\cs^{p-1}\cs^-$. Then $\Ga(\cs\bullet \r)$ consists of all sequences $(x_i)\in\set{0,1}^\N$ satisfying
\begin{equation}  \label{eq:31-1}
  (\cs^-\ca^+\cs^{p-1})^\f\lle \si^n((x_i))\lle (\ca^+\cs^{p-1}\cs^-)^\f\quad\forall n\ge 0.
\end{equation}
Take a sequence $(x_i)\in \Ga(\cs\bullet\r)\setminus \Ga(\cs)$. Then by (\ref{eq:8-2}) and (\ref{eq:31-1}) it follows that $x_{k+1}\ldots x_{k+m}=\cs^-$ or $\ca^+$ for some $k\ge 0$. If $x_{k+1}\ldots x_{k+m}=\cs^-$, then by taking $n=k$
 in (\ref{eq:31-1}) we obtain
 \[
 x_{k+m+1}x_{k+m+2}\ldots \lge (\ca^+\cs^{p-1}\cs^-)^\f.
 \]
 On the other hand, by taking $n=k+m$ in (\ref{eq:31-1}) we see  that the above inequality is indeed an equality. So, $x_{k+1}x_{k+2}\ldots=(\cs^-\ca^+\cs^{p-1})^\f$.

 If $x_{k+1}\ldots x_{k+m}=\ca^+$, then by taking $n=k$ in (\ref{eq:31-1}) we have
 \begin{equation}  \label{eq:31-2}
   x_{k+m+1}x_{k+m+2}\ldots\lle(\cs^{p-1}\cs^-\ca^+)^\f.
 \end{equation}
 Note by (\ref{eq:31-1}) that $x_{i+1}\ldots x_{i+m}\lge \cs^-$ for all $i\ge 0$. So  by (\ref{eq:31-2}) there must exist a $j\in\set{k+m, k+2m, \ldots, k+pm}$ such that $x_{j+1}\ldots x_{j+m}=\cs^-$. Then by the same argument as above we conclude that $x_{j+1}x_{j+2}\ldots=(\cs^-\ca^+\cs^{p-1})^\f$. So, $\Ga(\cs\bullet \r)\setminus \Ga(\cs)$ is at most countable.

\medskip
{\bf Case III.} $\r=01^p01^{p+t_1}\ldots 01^{p+t_N}01^{p+1}$, where $p\in\N$ and $\hat\r:=0 t_1\ldots t_N 1\in\F$.
Consider the substitution
\[
\eta_p:=U_1^p\circ U_0: 0\mapsto 01^p;\quad 1\mapsto 01^{p+1}.
\]
Then $\r=\eta_p(\hat\r)$ as shown in the proof of Proposition \ref{prop:Farey-new-form}.
Note by Lemma \ref{lem:Farey-property} that $ \L(\r)=1^{p+1}01^{p+t_1}01^{p+t_2}\ldots 01^{p+t_N}01^p0$ and $\L(\hat\r)=1 t_1\ldots t_N 0$. Then
\begin{equation}  \label{eq:31-11}
  \begin{split}
  \L(\r)^\f&=\si\big((01^{p+1}01^{p+t_1}\ldots 01^{p+t_N}01^p)^\f\big)
 =\si\big(\eta_p((1t_1\ldots t_N 0)^\f)\big)=\si\big(\eta_p(\L(\hat\r)^\f)\big).
  \end{split}
\end{equation}

\emph{Claim}: If $(x_i)\in \Ga(\cs\bullet \r)$ begins with $x_1\ldots x_m=\cs^-$, then there exists a unique sequence $(z_i)\in\set{0,1}^\N$ such that $(x_i)=\Phi_\cs\big(\eta_p(z_1z_2\ldots)\big)$.

Note that $\r$ begins with $01^p0$ and $\L(\r)$ begins with $1^{p+1}0$. Then $\cs\bullet\r$ begins with $\Phi_\cs(01^p0)=\cs^-\ca^+ \cs^{p-1}\cs^-$ and $\L(\cs\bullet\r)=\cs\bullet\L(\r)$ begins with $\Phi_\cs(1^{p+1}0)=\ca^+\cs^p\cs^-$. Let $(x_i)\in \Ga(\cs\bullet\r)$ with $x_1\ldots x_m=\cs^-$. Then
\begin{equation} \label{eq:31-12}
  \cs^-\ca^+\cs^{p-1}\cs^-\lle x_{n+1}\ldots x_{n+m(p+2)}\lle \ca^+\cs^p\cs^-\quad \forall n\ge 0.
\end{equation}
By taking $n=0$ in (\ref{eq:31-12}) it follows that
\begin{equation}  \label{eq:31-13}
  x_{m+1}\ldots x_{m(p+2)}\lge \ca^+\cs^{p-1}\cs^-.
\end{equation}
On the other hand, by taking $n=m$ in (\ref{eq:31-12}) we have
\begin{equation}  \label{eq:31-14}
  x_{m+1}\ldots x_{m(p+3)}\lle \ca^+\cs^p\cs^-.
\end{equation}
By (\ref{eq:31-12})--(\ref{eq:31-14}) it follows that
\[
\textrm{either}\quad x_{m+1}\ldots x_{m(p+2)}=\ca^+\cs^{p-1}\cs^-\quad\textrm{or}\quad x_{m+1}\ldots x_{m(p+3)}=\ca^+\cs^p\cs^-.
\]
Observe that in both cases we obtain a block ending with $\cs^-$. Then we can repeat the above argument indefinitely, and conclude that
\begin{equation}  \label{eq:31-15}
  (x_i)\in\set{\cs^-\ca^+\cs^{p-1}, \cs^-\ca^+\cs^p}^\N=\set{\Phi_\cs(\eta_p(0)), \Phi_\cs(\eta_p(1))}^\N.
\end{equation}
Since $\eta_p(0)=01^p$ and $\eta_p(1)=01^{p+1}$,  it follows from (\ref{eq:31-15}) and Lemma \ref{lem:connectible} that
\begin{equation}\label{eq:31-16}
(x_i)=\Phi_\cs(\eta_p(z_1))\Phi_\cs(\eta_p(z_2))\ldots =\Phi_\cs\big(\eta_p(z_1)\eta_p(z_2)\ldots\big)=\Phi_\cs\big(\eta_p(z_1z_2\ldots)\big)
\end{equation}
for some sequence $(z_i)\in\set{0,1}^\N$. The uniqueness of $(z_i)$ follows by the definition of the substitutions $\eta_p$ and $\Phi_\cs$. This proves the claim.

Now take a sequence $(x_i)\in \Ga(\cs\bullet \r)\setminus \Ga(\cs)$. Then by (\ref{eq:8-2}) and (\ref{eq:31-12}) we can find an $n_0\ge 0$ such that $x_{n_0+1}\ldots x_{n_0+m}=\cs^-$ or $\ca^+$. If  $x_{n_0+1}\ldots x_{n_0+m}=\ca^+$, then by (\ref{eq:31-12}) there must exist $n_1>n_0$  such that $x_{n_1+1}\ldots x_{n_1+m}=\cs^-$. So, without loss of generality, we may assume  $x_{n_0+1}\ldots x_{n_0+m}=\cs^-$. Then by the claim there is a unique sequence $(z_i)\in\set{0,1}^\N$ such that $x_{n_0+1}x_{n_0+2}\ldots=\Phi_\cs(\eta_p(z_1z_2\ldots))\in \Ga(\cs\bullet \r)$. By the definition of $\Ga(\cs\bullet \r)$ it follows that
\begin{equation} \label{eq:kk-3}
(\cs\bullet\r)^\f \lle\si^n\big(\Phi_{\cs}(\eta_p(z_1z_2\ldots))\big)\lle \L(\cs\bullet\r)^\f\quad\forall n\ge 0.
\end{equation}
Note by $\r=\eta_p(\hat{\r})$ and Lemma {\ref{lem:connectible} (ii)} that $(\cs\bullet\r)^\f=\Phi_\cs(\r^\f)=\Phi_\cs(\eta_p(\hat\r^\f))$. {Similarly}, by Lemma {\ref{lem:connectible} (ii), Lemma} \ref{lem:lyndon-commutative} and (\ref{eq:31-11}) it follows that $\L(\cs\bullet\r)^\f=\Phi_\cs(\L(\r)^\f)=\Phi_\cs\big(\si(\eta_p(\L(\hat\r)^\f))\big)$. {Thus}, (\ref{eq:kk-3}) can be rewritten as
\begin{equation*} \label{eq:28-5}
 \Phi_\cs(\eta_p(\hat\r^\f)) \lle\si^n\big(\Phi_\cs(\eta_p(z_1z_2\ldots))\big)\lle\Phi_\cs\big(\si(\eta_p(\L(\hat\r)^\f))\big)\quad\forall n\ge 0.
\end{equation*}
By Lemma \ref{lem:monotonicity-Phi} this implies that
\begin{equation} \label{eq:28-61}
\eta_p(\hat\r^\f)\lle \sigma^n\big(\eta_p(z_1z_2\ldots)\big)\lle \si\big(\eta_p(\L(\hat\r)^\f)\big)\quad\forall n\ge 0.
\end{equation}
Note  $(c_i)\prec(d_i)$ is equivalent to $\eta_p((c_i))\prec \eta_p((d_i))$. So, by (\ref{eq:28-61}) and the definition of $\eta_p$ it follows that
\begin{equation*} \label{eq:31-17}
\hat\r^\f\lle \si^n(z_1z_2\ldots)\lle \L(\hat\r)^\f\quad\forall n\ge 0,
\end{equation*}
{and hence} $(z_i)\in \Ga(\hat \r)$. Since $\hat\r\in\F$, we know that $\Ga(\hat \r)$ is finite by Lemma \ref{lem:finite-symmetric}. Hence there are only countably many choices for the sequence $(z_i)$, and thus by the claim there are only countably many choices for the tail sequence of $(x_i)$. Therefore, $\Ga(\cs\bullet \r)\setminus\Ga(\cs)$ is at most countable.
\end{proof}

Recall from Lemma \ref{lem:dim-survivorset} the symbolic survivor set $\K_\beta(t)$. To prove Theorem \ref{main:critical-basic-intervals} we also recall the following result from \cite[Lemma 3.7]{Kalle-Kong-Langeveld-Li-18}.

\begin{lemma} \label{lem:simplification-K-beta}
Let $\beta\in(1,2]$ and $t\in[0,1)$. If $\si^m(\de(\beta))\lle b(t,\beta)$, then
\[
\K_\beta(t)=\set{(d_i): b(t,\beta)\lle\si^n((d_i))\lle (\de_1(\beta)\ldots \de_m(\beta)^-)^\f~\forall n\ge 0}.
\]
\end{lemma}

\begin{proof}[Proof of Theorem \ref{main:critical-basic-intervals}]
{That the basic intervals $I^\cs, \cs\in\La$ are pairwise disjoint will be shown in Proposition \ref{prop:basic-interval-disjoint} below. In what follows, we fix a basic interval $I^\cs=[\beta_\ell^\cs, \beta_*^\cs]$.
Take $\beta\in I^\cs$}, and let $t^*=(\Phi_\cs(0^\f))_\beta= (\cs^-\ca^\f)_\beta$, where $\ca=\L(\cs)=A_1\ldots A_m$.  Then by (\ref{eq:basic-interval}) and Lemma \ref{lem:quasi-greedy expansion} it follows that
\begin{equation}  \label{eq:a2-1}
  \ca^\f=\de(\beta_\ell^\cs)\lle \de(\beta)\lle \de(\beta_*^\cs)=\ca^+\cs^-\ca^\f.
\end{equation}
Since $\cs=\s_1\bullet\cdots\bullet\s_k$ with each $\s_i\in\F$, by Proposition \ref{main:substitution} we have $\cs\in\LL$. If $\beta\in(\beta_\ell^\cs, \beta_*^\cs]$, then by Lemma \ref{lem:lyndon-equivalence} we have
\[
\si^n(\cs^-\ca^\f)\lle \ca^\f\prec\de(\beta)\quad\forall n\ge 0;
\]
and by Lemma \ref{lem:greedy-expansion} it follows that $\cs^-\ca^\f$ is the greedy $\beta$-expansion of $t^*$, i.e., $b(t^*, \beta)=\cs^-\ca^\f$. If $\beta=\beta_\ell^\cs$, then $\si^n(\cs^-\ca^\f)\lle \ca^\f=\de(\beta)$ for all $n\ge 0$; and in this case one can verify that the greedy $\beta$-expansion of $t^*$ is given by $b(t^*, \beta)=\cs0^\f$.

First we prove $\tau(\beta)\ge t^*$. Note that $\ca=A_1\ldots A_m$. Let
$t_N:= ((\cs^-\ca^N A_1\ldots A_j)^\f)_\beta,$ where the index $j\in\set{1,\ldots, m}$ satisfies
$\cs=\S(\ca)=A_{j+1}\ldots A_{m} A_1\ldots A_j.$
Note {by Lemma \ref{lem:lyndon-equivalence}} that $A_{i+1}\ldots A_{m}\prec A_1\ldots A_{m-i}$ for any $0<i<m$. Then by (\ref{eq:a2-1}) one can verify that
\begin{align*} \label{eq:9-1}
\si^n\big((\cs^-\ca^N A_1\ldots A_j)^\f\big)&=\si^n\left(\cs^-\ca^{N+1}\big(A_1\ldots A_j^-\ca^{N+1}\big)^\f\right)\\
&\prec \ca^\f\lle\de(\beta)\qquad \forall n\ge 0.
\end{align*}
So, $b(t_N, \beta)=(\cs^-\ca^N A_1\ldots A_j)^\f$. This implies that any sequence $(x_i)$ constructed by arbitrarily concatenating blocks of the form
\[
\cs^-\ca^k A_1\ldots A_j,\quad  k>N
\]
satisfies $(\cs^-\ca^N A_1\ldots A_j)^\f\lle\si^n((x_i))\prec\de(\beta)$ for all $n\ge 0$. So,
\[
\set{\cs^-\ca^{N+1}A_1\ldots A_j, \cs^- \ca^{N+2}A_1\ldots A_j}^\N\subset\K_\beta(t_N).
\]
By Lemma \ref{lem:dim-survivorset} this implies that $\dim_H K_\beta(t_N)>0$ for all $N\ge 1$. Thus $\tau(\beta)\ge t_N$ for all $N\ge 1$. Note that $t_N\nearrow t^*$ as $N\to\f$. We then conclude that $\tau(\beta)\ge t^*$.

Next we prove $\tau(\beta)\le t^*$. By (\ref{eq:a2-1}) and Lemma \ref{lem:simplification-K-beta} it follows  that
\begin{equation} \label{eq:a2-2}
\begin{split}
\K_\beta(t^*)&\subset\set{(x_i): \cs^-\ca^\f\lle \si^n((x_i))\prec \ca^+\cs^-\ca^\f~\forall n\ge 0}\\
&=\set{(x_i): \cs^-\ca^\f\lle \si^n((x_i))\lle \ca^\f~\forall n\ge 0}=:\Ga.
\end{split}
\end{equation}
Note by Proposition \ref{prop:countable} that $\Ga(\cs)=\set{(x_i): \cs^\f\lle\si^n((x_i))\lle \ca^\f~\forall n\ge 0}$ is a countable subset of $\Ga$. Furthermore, any sequence in the difference set $\Ga\setminus\Ga(\cs)$ must end with $\cs^-\ca^\f$. {As a result}, $\Ga$ is also countable. By (\ref{eq:a2-2}) this implies that $\dim_H K_\beta(t^*)=0$, and thus $\tau(\beta)\le t^*$. This completes the proof.
\end{proof}

\section{Geometrical structure of the basic intervals and {exceptional} sets} \label{sec:geometric-structure}

In this section we will prove Theorem \ref{main:geometrical-structure-basic-intervals}.
The proof will be split into two subsections. In Subsection \ref{subsec:bifurcation-sets} we {demonstrate} the tree structure of the Lyndon intervals $J^\cs, \cs\in\La$ and the relative {exceptional} sets $E^\cs, \cs\in\La$, {from which it follows that the basic intervals $I^\cs=[\beta_\ell^\cs, \beta_*^\cs]$, $\cs\in\La$ are pairwise disjoint.} We show that each relative {exceptional} set $E^\cs$ has zero box-counting dimension, and the {exceptional} set $E$ has zero packing dimension. In Subsection \ref{sec:exceptional-set} we prove that the {infinitely Farey} set $E_\f$ has zero Hausdorff dimension.

\subsection{Tree structure of the Lyndon intervals and relative exceptional sets} \label{subsec:bifurcation-sets}

 Given $\cs\in\La$, recall from Definitions \ref{def:basic-intervals} and \ref{def:lyndon-intervals} the basic interval $I^\cs=[\beta_\ell^\cs,\beta_*^\cs]$ and the Lyndon interval $J^\cs=[\beta_\ell^\cs, \beta_r^\cs]$  generated by $\cs$, respectively. Then by Lemmas \ref{lem:quasi-greedy expansion} and \ref{lem:lyndon-equivalence} it follows that
\begin{equation} \label{eq:10-1}
\de(\beta_\ell^\cs)=\L(\cs)^\f,\quad \de(\beta_*^\cs)=\L(\cs)^+\cs^-\L(\cs)^\f\quad\textrm{and}\quad \de(\beta_r^\cs)=\L(\cs)^+ \cs^\f.
\end{equation}
First we show that the Lyndon intervals $J^\cs, \cs\in\La$ have a tree structure.

\begin{proposition}  \label{prop:basic-interval-disjoint}
  Let $\cs\in{\La}$. Then $I^\cs\subset J^{\cs}$. Furthermore,
  \begin{itemize}
  \item [{\rm(i)}] for any $\r\in\F$ we have $J^{\cs\bullet\r}\subset J^\cs\setminus I^\cs$;
  \item[{\rm(ii)}] for any two different words $\r,\r'\in\F$ we have $J^{\cs\bullet\r}\cap J^{\cs\bullet\r'}=\emptyset$.
  \end{itemize}
\end{proposition}

\begin{proof}
Let $I^\cs=[\beta_\ell^\cs, \beta_*^\cs]$ and $J^{\cs}=[\beta_\ell^\cs, \beta_r^\cs]$.
Then by (\ref{eq:10-1}) it follows that
\[
\de(\beta_*^\cs)=\L(\cs)^+\cs^-\L(\cs)^\f\prec \L(\cs)^+\cs^\f=\de(\beta_r^\cs),
\]
which implies $\beta_*^\cs<\beta_r^\cs$ by Lemma \ref{lem:quasi-greedy expansion}. So $I^\cs \subset J^\cs$.

For {(i),} let $\r\in\F$. Then $\r$ begins with digit $0$ and ends with digit $1$. By {Lemma \ref{lem:connectible} (ii) and Lemma} \ref{lem:lyndon-commutative} this implies that
 \begin{align*}
  \de(\beta_\ell^{\cs\bullet\r})&= \L(\cs\bullet\r)^\f=(\cs\bullet\L(\r))^\f=\Phi_\cs(\L(\r)^\f)\\
  &\succ \Phi_\cs(10^\f)=\L(\cs)^+ \cs^-\L(\cs)^\f=\de(\beta_*^\cs).
 \end{align*}
So, $\beta_\ell^{\cs\bullet\r}>{\beta_*^\cs}$. Furthermore, by Lemmas \ref{lem:mon-1}, {\ref{lem:connectible} (ii)} and \ref{lem:lyndon-commutative} it follows that
  \begin{align*}
  \de(\beta_r^{\cs\bullet\r})&=\L(\cs\bullet\r)^+(\cs\bullet\r)^\f\\
  &=\Phi_\cs(\L(\r)^+)\Phi_\cs(\r^\f)=\Phi_\cs(\L(\r)^+\r^\f)\\
  &\prec \Phi_\cs(1^\f)=\L(\cs)^+ \cs^\f=\de(\beta_r^\cs).
  \end{align*}
  This proves $\beta_r^{\cs\bullet\r}<\beta_r^\cs$. Hence, $J^{\cs\bullet\r}=[\beta_\ell^{\cs\bullet\r}, \beta_r^{\cs\bullet\r}]\subset(\beta_*^\cs, \beta_r^\cs]=J^\cs\setminus I^\cs$.

  Next we prove {(ii)}. Let $\r, \r'$ be two different Farey words in $\F$. By Lemma \ref{lem:farey-interval} (i) it follows that $J^{\r}\cap J^{\r'}=\emptyset$. Write $J^{\r}=[\beta_\ell^\r, \beta_r^\r]$ and $J^{\r'}=[\beta_\ell^{\r'}, \beta_r^{\r'}]$. Since $J^\r$ and $J^{\r'}$ are disjoint, we may assume $\beta_r^\r<\beta_\ell^{\r'}$. By (\ref{eq:10-1})  and Lemma \ref{lem:quasi-greedy expansion} it follows that
  \begin{equation} \label{eq:basic-interval-1}
  \L(\r)^+\r^\f=\de(\beta_r^\r)\prec \de(\beta_\ell^{\r'})=\L(\r')^\f.
  \end{equation}
  Then by (\ref{eq:10-1}), (\ref{eq:basic-interval-1}) and Lemma \ref{lem:monotonicity-Phi} we obtain that
  \begin{align*}
  \de(\beta_r^{\cs\bullet\r})&=\L(S\bullet\r)^+(S\bullet\r)^\f=\Phi_\cs(\L(\r)^+\r^\f)\\
	&\prec\Phi_\cs(\L(\r')^\f)=\L(\cs\bullet\r')^\f=\de(\beta_\ell^{\cs\bullet\r'}).
  \end{align*}
  It follows that $\beta_r^{\cs\bullet\r}<\beta_\ell^{\cs\bullet\r'}$, and hence $J^{\cs\bullet\r}\cap J^{\cs\bullet\r'}=\emptyset$.
\end{proof}

\begin{remark} \label{rem:lyndon-interval}
 Proposition \ref{prop:basic-interval-disjoint} implies that the Lyndon intervals $J^\cs, \cs\in\La$ have a tree structure. More precisely, we say $J^{\mathbf R}$ is an \emph{offspring} of $J^\cs$ if there exists a word $\mathbf T\in\La$ such that $\mathbf R=\cs\bullet\mathbf T$. Then any offspring of $J^\cs$ is a subset of $J^\cs$. Furthermore, if $J^{\cs'}$ is not an offspring of $J^\cs$ {and $J^\cs$ is not an offspring of $J^{\cs'}$}, then Proposition \ref{prop:basic-interval-disjoint} implies that $J^{\cs'}\cap J^\cs=\emptyset$.
 Consequently, the basic intervals $I^\cs, \cs\in\La$ are pairwise disjoint. 
\end{remark}%

Recall from Section \ref{sec: Introduction} the {exceptional} set $E=(1,2]\setminus\bigcup_{\r\in\F}J^\r$ and the relative {exceptional} sets $E^\cs=(J^\cs\setminus I^\cs)\setminus\bigcup_{\r\in\F}J^{\cs\bullet\r}$ with $\cs\in\La$. Next we will show that $E$ is bijectively mapped to $E^\cs$ via the map
\begin{equation} \label{eq:a5-1}
\Psi_\cs: (1,2]\to  J^\cs\setminus I^\cs=(\beta_*^\cs, \beta_r^\cs];\quad \beta\mapsto \de^{-1}\circ\Phi_\cs\circ\de(\beta),
\end{equation}
where $\de(\beta)$ is the quasi-greedy $\beta$-expansion of $1$.

We mention that $\Psi_\cs$ is not surjective, which somewhat complicates the proof of Proposition \ref{lem:geometry-basic-interval-4} below. For example, let $\cs=011$. Then $\de(\beta_*^\cs)=111010(110)^\f$ and $\de(\beta_r^\cs)=111(011)^\f$. Take $\beta\in(\beta_*^\cs, \beta_r^\cs]$ such that $\de(\beta)=1110110^\f$. One can verify that $\beta\notin\Psi_\cs((1,2])$, since $\de(\beta)$ cannot be written as a concatenation of words from $\{\cs,\cs^-,\L(\cs),\L(\cs)^+\}$.

\begin{lemma}  \label{lem:geometry-basic-interval-2}
For any $\cs\in\La$ the map $\Psi_\cs$ is well-defined and strictly increasing.
\end{lemma}

\begin{proof}
Let $\cs\in\La$ with $\ca=\L(\cs)$.  First we show that the map $\Psi_\cs: (1,2]\to J^\cs\setminus I^\cs=(\beta_*^\cs, \beta_r^\cs]$  is well-defined. Note that
\begin{equation} \label{eq:a5-2}
\de(\beta_*^\cs)=\ca^+ \cs^-\ca^\f=\Phi_\cs(10^\f)\quad\textrm{and}\quad \de(\beta_r^\cs)=\ca^+\cs^\f=\Phi_\cs(1^\f).
\end{equation}
Take $\beta\in(1,2]$. Then $10^\f\prec \de(\beta)\lle 1^\f$. By Lemma \ref{lem:monotonicity-Phi} and (\ref{eq:a5-2}) it follows that
\[
\de(\beta_*^\cs)\prec\Phi_\cs(\de(\beta))\lle \de(\beta_r^\cs).
\]
Thus, by Lemma \ref{lem:quasi-greedy expansion} it suffices to prove that
\begin{equation} \label{eq:a5-3}
  \si^n(\Phi_\cs(\de(\beta)))\lle\Phi_\cs(\de(\beta))\quad\forall n\ge 0.
\end{equation}
Note by Lemma \ref{lem:quasi-greedy expansion} that $\si^n(\de(\beta))\lle \de(\beta)$ for all $n\ge 0$, and $\de(\beta)$ begins with digit $1$. Thus, (\ref{eq:a5-3}) follows by Lemma \ref{lem:cri-exception-1} (i), and we conclude that the map $\Psi_\cs$ is well-defined.

The monotonicity of $\Psi_\cs =\de^{-1}\circ\Phi_\cs\circ\de$ follows since both maps $\de$ and $\Phi_\cs$ are strictly increasing by Lemmas \ref{lem:quasi-greedy expansion} and \ref{lem:monotonicity-Phi}, respectively. This completes the proof.
\end{proof}

 \begin{proposition}  \label{lem:geometry-basic-interval-4}
   For any $\cs\in\La$ we have $\Psi_\cs(E)=E^\cs$.
 \end{proposition}

\begin{proof}
We {first} prove that
\begin{equation} \label{eq:3-1}
  \Psi_\cs(\beta_\ell^\r)=\beta_\ell^{\cs\bullet\r}\quad\textrm{and}\quad \Psi_\cs(\beta_r^\r)=\beta_r^{\cs\bullet\r}, \qquad \forall \r\in{\LL}.
  \end{equation}
 Observe that $\de(\beta_\ell^\r)=\L(\r)^\f$. Then by {Lemma \ref{lem:connectible} (ii) and Lemma} \ref{lem:lyndon-commutative} it follows that
  \begin{align*}
  \Phi_\cs(\de(\beta_\ell^\r))=\Phi_\cs(\L(\r)^\f)=(\cs\bullet\L(\r))^\f=\L(\cs\bullet \r)^\f=\de(\beta_\ell^{\cs\bullet \r}),
  \end{align*}
  {so} $\Psi_\cs(\beta_\ell^\r)=\beta_\ell^{\cs\bullet \r}$. Similarly, {since $\de(\beta_r^\r)=\L(\r)^+\r^\f$, Lemmas \ref{lem:mon-1}, {\ref{lem:connectible} (ii)} and \ref{lem:lyndon-commutative} imply} that
  \begin{align*}
  \Phi_\cs(\de(\beta_r^\r))&=\Phi_\cs(\L(\r)^+\r^\f)=\Phi_\cs(\L(\r)^+)\Phi_\cs(\r^\f)\\
  &=(\cs\bullet\L(\r))^+\Phi_\cs(\r)^\f=\L(\cs\bullet\r)^+(\cs\bullet \r)^\f=\de(\beta_r^{\cs\bullet\r}).
  \end{align*}
  We conclude that $\Psi_\cs(\beta_r^{\r})=\beta_r^{\cs\bullet\r}$. This proves (\ref{eq:3-1}).
	
Note by Lemma \ref{lem:quasi-greedy expansion} (ii) that the map $\beta\mapsto \de(\beta)$ is left continuous in $(1,2]$, and is right continuous at a point $\beta_0$ if and only if $\de(\beta_0)$ is {\em not} periodic. Hence by Lemma \ref{lem:periodic-qg-expansion} it follows that the map $\beta\mapsto \de(\beta)$ is continuous at $\beta_0$ if and only if $\de(\beta_0)$ is not of the form $\L(\r)^\f$ for a Lyndon word $\r$. Since $\Phi_\cs$ is clearly continuous with respect to the order topology and the map $\delta^{-1}$ is continuous, it follows that $\Psi_\cs$ is continuous at $\beta_0$ if and only if $\de(\beta_0)$ is not of the form $\L(\r)^\f$ for a Lyndon word $\r$. Moreover, $\Psi_\cs$ is left continuous everywhere.

Note that if $\de(\beta_0)=\L(\r)^\f$, then by Lemma \ref{lem:quasi-greedy expansion} (ii) it follows that as $\beta$ decreases to $\beta_0$ the sequence $\de(\beta)$ converges to $\L(\r)^+0^\f$ with respect to the order topology, i.e., $\lim_{\beta\searrow \beta_0}\de(\beta)=\L(\r)^+0^\f$.
Since $\Psi_\cs$ is increasing, it follows that
\[
\Psi_\cs((1,2])=(J^\cs\backslash I^\cs)\backslash \bigcup_{\r\in\Omega_L^*} (p^{\r},q^{\r}],
\]
where $\de(p^\r)=\Phi_\cs(\L(\r)^\f)$ and $\de(q^\r)=\Phi_\cs(\L(\r)^+0^\f)$. Note that $(p^\r,q^\r]\subset J^{\cs\bullet \r}$.
By Lemma \ref{lem:farey-interval} (iii) there is a (unique) Farey word $\hat{\r}$ such that {$J^\r\subset J^{\hat r}$. Applying \eqref{eq:3-1} to both $\r$ and $\hat{\r}$ and using Lemma \ref{lem:geometry-basic-interval-2} we conclude that $J^{\cs\bullet\r}\subset J^{\cs\bullet\hat{\r}}$. Hence, $(p^\r,q^\r]\subset J^{\cs\bullet \hat{\r}}$.}
Therefore, if $\beta\in E^\cs=(J^\cs\backslash I^\cs)\backslash \bigcup_{\r\in\Omega_F^*} J^{\cs\bullet \r}$, then $\beta$ lies in the range of $\Psi_\cs$. This implies that
\begin{equation} \label{eq:inclusion-1}
E^\cs\subset\Psi_\cs((1,2]).
\end{equation}

	Now assume first that $\beta\in\Psi_\cs(E)$. Then $\beta=\Psi_\cs({\hat\beta})$, where ${\hat\beta}\not\in J^\r$ for any $\r\in\F$. Hence $\beta\not\in J^{\cs\bullet \r}$ for any $\r\in\F$ by \eqref{eq:3-1} and since $\Psi_\cs$ is increasing. {Therefore}, $\beta\in E^\cs$.
	
	Conversely, suppose $\beta\in E^\cs$. By (\ref{eq:inclusion-1}), $\beta=\Psi_\cs({\hat\beta})$ for some ${\hat\beta}\in(1,2]$. If ${\hat\beta}\in J^\r$ for some $\r\in\F$, then $\beta\in\Psi_\cs(J^\r)\subset J^{\cs\bullet\r}$ by \eqref{eq:3-1}, contradicting that $\beta\in E^\cs$. Hence, ${\hat\beta}\in E$ and then $\beta\in\Psi_\cs(E)$. This completes the proof.
	\end{proof}

Kalle et al.~proved in \cite[Theorem C]{Kalle-Kong-Langeveld-Li-18} that the Farey intervals $J^\r, \r\in\F$ cover the whole interval $(1,2]$ up to a set of zero Hausdorff dimension. Here we strengthen this result and show that the {exceptional} set $E$ is uncountable and has zero packing dimension. Furthermore, we show that each relative {exceptional} set $E^\cs$ is uncountable and has zero box-counting dimension. The proof uses the following simple lemma.

\begin{lemma} \label{lem:Lyndon-interval-length}
Let $J^\cs=[\beta_\ell^\cs,\beta_r^\cs]=:[p,q]$ be any Lyndon interval. Then {the length of $J^\cs$ satisfies}
\[
|J^\cs|\leq \frac{q}{q-1}q^{-|\cs|}.
\]
\end{lemma}

\begin{proof}
Since $\de(p)=\L(\cs)^\f$ and $\de(q)=\L(\cs)^+\cs^\f$, we have
\[
\big(\L(\cs)^+0^\f\big)_p=1=\big(\L(\cs)^+\cs^\f\big)_q=:((c_i))_q.
\]
It follows that
\[
|J^\cs|=q-p=\sum_{i=1}^{|\cs|}\frac{c_i}{q^{i-1}}+\sum_{i=|\cs|+1}^\f\frac{c_i}{q^{i-1}}-\sum_{i=1}^{|\cs|}\frac{c_i}{p^{i-1}} \le\sum_{i=|\cs|+1}^\f\frac{1}{q^{i-1}}=\frac{q}{q-1}q^{-|\cs|},
\]
as required.
\end{proof}

\begin{proposition} \label{prop:geometry-basic-interval-2}
\mbox{}
\begin{enumerate}[{\rm(i)}]
\item The {exceptional} set
  \[
  E=(1,2]\setminus\bigcup_{\r\in\F}J^{\r}
  \]
 is uncountable and has zero packing dimension. 
\item [{\rm(ii)}] For any $\cs\in\La$ the relative {exceptional} set
  \[
  E^\cs=(J^\cs\setminus I^\cs)\setminus \bigcup_{\r\in\F}J^{\cs\bullet\r}
  \]
 is uncountable and has zero box-counting dimension.  
\end{enumerate}
\end{proposition}

\begin{proof}
(i) First we prove $\dim_P E=0$.
Let $\rho_N\in(1,2]$ such that  $\de(\rho_N)=(10^{N-1})^\f$. Then by  Lemma \ref{lem:quasi-greedy expansion}  it follows that $\rho_N\searrow 1$ as $N\to\f$. Thus
 $E=\bigcup_{N=1}^\f(E\cap[\rho_N,2])$. By {the countable stability of packing dimension} {(cf.~\cite{Falconer_1990})} it suffices to prove that 
   \begin{equation} \label{eq:20-1}
   \dim_B(E\cap[\rho_N,2])=0\quad\textrm{for all}\quad N\in\N.
   \end{equation}

  Let $N\in\N$. Take a Farey interval $J^\s:=[p, q]\subset[\rho_N, 2]$ with $\s=s_1\ldots s_m\in\F$ such that
	\begin{equation} \label{eq:mlb}
  m> N+2-3\log_2(\rho_N-1).
  \end{equation}
   Write $\L(\s)=a_1\ldots a_m$. Then  $\de(p)=(a_1\ldots a_m)^\f$. Since $p\ge \rho_N$, by Lemma \ref{lem:quasi-greedy expansion}  we have  $(a_1\ldots a_m)^\f=\de(p)\lge\de(\rho_N)=(10^{N-1})^\f$, which implies that $a_1\ldots a_{N+1}\lge 10^{N-1}1$. Then by Proposition \ref{prop:Farey-new-form} we conclude that
  \begin{equation} \label{eq:slb}
  s_1\ldots s_{N+1}\lge 0^{N}1.
  \end{equation}
Note that
\[
(\L(\s)^+ 0^\f)_p=1=(\L(\s)^+{\s^\f})_q=:((c_i))_q.
\]
So, by (\ref{eq:slb}) it follows that
  \begin{align*}
   \sum_{i=1}^m\frac{c_i}{p^i}=1=\sum_{i=1}^m\frac{c_i}{q^i}+\sum_{i=m+1}^\f\frac{c_i}{q^i}>\sum_{i=1}^m \frac{c_i}{q^i}+\frac{1}{q^{m+N+1}},
     \end{align*}
  which implies
  \[
  \frac{1}{q^{m+N+1}}<\sum_{i=1}^\f\left(\frac{1}{p^i}-\frac{1}{q^i}\right)=\frac{q-p}{(p-1)(q-1)}.
  \]
 Whence,
  \begin{equation}  \label{eq:20-2}
    |J^\s|=q-p>\frac{(p-1)(q-1)}{q^{N+1}}q^{-m}\ge \frac{(\rho_N-1)^2}{2^{N+1}}q^{-m}.
  \end{equation}
On the other hand, by Lemma \ref{lem:Lyndon-interval-length} it follows that
\begin{equation} \label{eq:20-3}
|J^\s|\le \frac{q}{q-1}q^{-m}\le\frac{2}{\rho_N-1}q^{-m} {\leq \frac{2}{\rho_N-1}\rho_N^{-m}}.
\end{equation}

Now we list all of the Farey intervals in $[\rho_N,2]$ in a decreasing order according to their length, say $J^{\s_1}, J^{\s_2}, \ldots$. In other words, $|J^{\s_i}|\ge |J^{\s_j}|$ for any $i<j$. For a Farey interval $J^\s$, if $J^\s=J^{\s_k}$ we then define its \emph{order index} as $o(J^\s)=k$.

Set $C_N:=2\log 2/\log \rho_N$. Let $J^{\s'}$ be a Farey interval with $|\s'|>C_N m$. Then by (\ref{eq:mlb}), (\ref{eq:20-2}) and (\ref{eq:20-3}) it follows  that
 \begin{equation*} 
   |J^{\s'}|\le\frac{2}{\rho_N-1}\rho_N^{-C_N m}=\frac{2}{\rho_N-1}2^{-2m}<\frac{(\rho_N-1)^2}{2^{N+1}}2^{-m}\le |J^\s|.
 \end{equation*}
 This implies that
 \begin{equation}\label{eq:20-5}
o(J^\s)\le\sum_{k=1}^{\lfloor C_N m\rfloor}\#\set{\s'\in\F: |\s'|=k}\le\sum_{k=1}^{\lfloor C_N m\rfloor}(k-1)< C_N^2 m^2,
\end{equation}
where the second inequality follows by (\ref{eq:kk-6}) since the number of non-degenerate Farey words of length $k$ is at most $k-1$ {(see \cite[Proposition 2.3]{Carminati-Isola-Tiozzo-2018})}. Together with (\ref{eq:20-3}), (\ref{eq:20-5}) implies that
 \begin{align*}
   \liminf_{i\to\f}\frac{-\log|J^{\s_i}|}{o(J^{\s_i})}&=+\f.
  \end{align*}
 Note that $[\rho_N,2]\cap E=[\rho_N,2]\setminus\bigcup_{\s\in\F}J^\s$. So, by \cite[Proposition 3.6]{Falconer_1997} we conclude {\eqref{eq:20-1}}. This proves $\dim_P E=0$.

 Next we prove that $E$ is uncountable. For $\s\in\F$ let $\hat J^\s=(\beta_\ell^\s, \beta_r^\s)$ be the interior of the Farey interval $J^\s=[\beta_\ell^\s,\beta_r^\s]$. By Lemma \ref{lem:farey-interval} (i) it follows that {the compact set}
 \[
 \hat E:=[1,2]\setminus\bigcup_{\s\in\F}\hat J^\s
 \]
{is nonempty and has no isolated points. Hence, $\hat E$ is a perfect set and is therefore uncountable. Since $\hat{E}\backslash E$ is countable, it follows that $E$ is uncountable as well.}

 (ii) In a similar way we prove $\dim_B E^\cs=0$. Note that $E^\cs=(\beta_*^\cs, \beta_r^\cs]\setminus\bigcup_{\r\in\F}J^{S\bullet\r}$.
Fix a Farey word $\r=r_1\ldots r_m$.
   Then the Lyndon interval $J^{\cs\bullet\r}=[\beta_\ell^{\cs\bullet\r}, \beta_r^{\cs\bullet\r}]=:[p_\r, q_\r]$ satisfies
  \begin{equation*} \label{eq:kd-2}
  {(\L(\cs\bullet\r)^+ 0^\f)_{p_\r}=1=(\L(\cs\bullet\r)^+ \; (\cs\bullet\r)^\f)_{q_\r}=:((d_i))_{q_\r}.}
  \end{equation*}
  So,
  \[
  \sum_{i=1}^{m|\cs|}\frac{d_i}{p_\r^i}=1=\sum_{i=1}^{m|\cs|}\frac{d_i}{q_\r^i}+\sum_{i=m|\cs|+1}^\f\frac{d_i}{q_\r^i}>\sum_{i=1}^{m|\cs|}\frac{d_i}{q_\r^i}+\frac{1}{q_\r^{(m+1)|\cs|+1}},
  \]
  where the inequality follows by observing that $\cs\in\LL$ and thus $\cs\bullet\r^\f\lge 0^{|\cs|}10^\f$. Therefore,
  \[
  \frac{1}{q_\r^{(m+1)|\cs|+1}}\le \sum_{i=1}^\f\left(\frac{1}{p_\r^i}-\frac{1}{q_\r^i}\right)=\frac{q_\r-p_\r}{(p_\r-1)(q_\r-1)},
  \]
which implies
\begin{equation*} \label{eq:kd-3}
|J^{\cs\bullet\r}|=q_\r-p_\r\ge \frac{(p_\r-1)(q_\r-1)}{q_\r^{|\cs|+1}}q_\r^{-|\cs\bullet\r|}>\frac{(\beta_*^\cs-1)^2}{2^{|\cs|+1}}q_\r^{-|\cs\bullet\r|}.
\end{equation*}
On the other hand, by {Lemma \ref{lem:Lyndon-interval-length}} it follows that
\begin{equation*} \label{eq:kd-4}
|J^{\cs\bullet\r}|\leq \frac{q_\r}{q_\r-1}q_\r^{-m|\cs|}\le\frac{2}{\beta_*^\cs-1}q_\r^{-|\cs\bullet\r|}.
\end{equation*}
Now let $(J^{\cs\bullet \r_i})$ be an enumeration of the intervals $J^{\cs\bullet \r}$, $\r\in\F$, arranged in order by decreasing length. Then by a similar argument as in (i) above, we obtain
\[
 \liminf_{i\to\f}\frac{-\log|J^{\cs\bullet\r_i}|}{\log o(J^{\cs\bullet\r_i})}{=+\f}.
\]
Thus, $\dim_B E^\cs=0$.

Finally, since {we showed in (i) that $E$ is uncountable, we conclude by Lemma \ref{lem:geometry-basic-interval-2} and Proposition \ref{lem:geometry-basic-interval-4}} that $E^\cs=\Psi_\cs(E)$ is also uncountable. This completes the proof.
\end{proof}

 \subsection{{The infinitely Farey} set} \label{sec:exceptional-set}
{Recall from (\ref{eq:E_s}) that
 \begin{equation*}
     E_\f=\bigcap_{n=1}^\f\bigcup_{\cs\in\La(n)}J^\cs,
 \end{equation*}
where
 \begin{equation*} \label{eq:kk-5}
 \La(n)=\set{\s_1\bullet\s_2\bullet\cdots\bullet\s_n: ~~ \s_i\in\F~\textrm{for all }1\le i\le n}.
 \end{equation*}
In particular, $\La(1)=\F$ and $\La=\bigcup_{n=1}^\f\La(n)$.
 Note that $(1,2]=E\cup\bigcup_{\s\in\F}J^\s$. Furthermore, {for each word $\cs\in \La$ we have}
 \[
 J^\cs\setminus I^\cs=E^\cs\cup\bigcup_{\r\in\F}J^{\cs\bullet\r}.
 \]
 By iteration of the above equation we {obtain} the following partition of the interval $(1,2]$.

\begin{lemma} \label{lem:presentation-J-s}
 The interval $(1,2]$ can be partitioned as
 \[
 (1,2]=E\cup E_\f\cup\bigcup_{\cs\in\La}E^\cs\cup\bigcup_{\cs\in\La}I^\cs.
 \]
 \end{lemma}

 To complete the proof of Theorem \ref{main:geometrical-structure-basic-intervals} we still need the following dimension result for $E_\f$.}

 \begin{proposition}  \label{prop:dim-exceptional-set}
    We have $\dim_H E_\f=0$.
 \end{proposition}

 \begin{proof}
Note by \eqref{eq:E_s} that
\begin{equation} \label{eq:E-infinity-cover}
E_\f=\bigcap_{n=1}^\f\bigcup_{\cs\in\La(n)}J^\cs \subset \bigcap_{n=1}^\f\bigcup_{\cs\in\La:|\cs|\geq n}J^\cs.
\end{equation}
  This suggests covering $E_\f$ by the intervals $J^\cs$ for $\cs\in\La$ with $|\cs|\geq n$ for a sufficiently large $n$. To this end, we first estimate the diameter of $J^\cs$.
Take $\cs\in \La$ with $|\cs|=m$, and write $J^\cs=[p,q]$. Then
$\de(p)=\L(\cs)^\f$ and $\de(q)=\L(\cs)^+ \cs^\f$,
and it follows {from Lemma \ref{lem:Lyndon-interval-length}} that
\begin{equation} \label{eq:H-diameter}
{|J^\cs|\leq \frac{q}{q-1}\cdot q^{-m}.}
\end{equation}
{Let $(\beta_n)$ be an arbitrary sequence in $(1,2)$ decreasing to $1$.
We will use \eqref{eq:H-diameter} to show that $\dim_H (E_\f\cap(\beta_n,2])=0$ for all $n\in\N$, so the result will follow from the countable stability of Hausdorff dimension {(cf.~\cite{Falconer_1990})}. Fix $n\in\N$. Observe that, if $J^\cs=[p,q]$ intersects $(\beta_n,2]$, then $q>\beta_n$ and so by \eqref{eq:H-diameter},
\begin{equation} \label{eq:J-diameter-bound}
|J^\cs|\leq \frac{2}{\beta_n-1}\beta_n^{-m}=:C_n\beta_n^{-m}.
\end{equation}
}

Next, we count how many words $\cs\in\La$ there are with $|\cs|=m$. Call this number $N_m$. Observe that if $\cs=\s_1\bullet\s_2\bullet\cdots\bullet \s_k$ and $|\s_i|=l_i$ for $i=1,\dots,k$, then $|\cs|=l_1 l_2\cdots l_k$.
Note by \cite[Proposition 2.3]{Carminati-Isola-Tiozzo-2018} that $\#\set{\r\in\F: |\r|=l}\le l-1$ for any $l\ge 2$. Thus, for any given tuple $(l_1,\dots,l_k)$, the number of possible choices for the words $\s_1,\dots,\s_k$ is at most $l_1 l_2\cdots l_k=|\cs|=m$. It remains to estimate how many {\em ordered factorizations} of $m$ there are, that is to estimate the number
\[
f_m:=\#\{(l_1,\dots,l_k): k\in\N, \quad l_i\in\N_{\geq 2}\ \forall i\quad \mbox{and}\quad l_1l_2\cdots l_k=m\}.
\]
By considering the possible values of $l_1$, it is easy to see that $f_m$ satisfies the recursion 
\[
f_m=\sum_{d|m, d>1}f_{m/d},
\]
where we set $f_1:=1$. (See \cite{Hille_1936}.) We claim that $f_m\leq m^2$. This is trivial for $m=1$, so let $m\geq 2$ and assume $f_n\leq n^2$ for all $n<m$; then
\begin{equation*}
f_m=\sum_{d|m, d>1}f_{m/d}\leq \sum_{d|m, d>1}\left(\frac{m}{d}\right)^2
\leq m^2\sum_{d=2}^\f \frac{1}{d^2}=m^2\left(\frac{\pi^2}{6}-1\right)<m^2.
\end{equation*}
This proves the claim, and we thus conclude that $N_m\leq m^3$. 
{Now, given $\ep>0$ and $\de>0$, choose $N$ large enough so that $C_n (\beta_n)^{-N}<\de$. Using \eqref{eq:E-infinity-cover} and \eqref{eq:J-diameter-bound}, we obtain}
\[
\mathcal H_\de^\ep(E_\f{\cap(\beta_n,2]})\leq \sum_{\cs\in\La: |\cs|\geq N,\ J^\cs\cap(\beta_n,2]\neq\emptyset} |J^\cs|^\ep\leq \sum_{m=N}^\f m^3 {C_n^\ep\beta_n^{-m\ep}}\to 0
\]
as $N\to\infty$. This shows that $\dim_H {(E_\f\cap(\beta_n,2])}=0$, {as desired}.
\end{proof}

\begin{proof}[Proof of Theorem \ref{main:geometrical-structure-basic-intervals}]
  The theorem follows by  Proposition \ref{prop:geometry-basic-interval-2}, Lemma \ref{lem:presentation-J-s} and Proposition \ref{prop:dim-exceptional-set}.
\end{proof}

 \section{Critical values in the exceptional sets} \label{sec:cri-bifurcation-sets}

 By Proposition \ref{lem:critical-non-Fareyinterval}, Theorem \ref{main:critical-basic-intervals} and Theorem \ref{main:geometrical-structure-basic-intervals} it suffices to determine the critical value $\tau(\beta)$ for
 \[
 \beta\in \bigcup_{\cs\in\La}E^\cs\cup E_\f.
 \]

 First we {compute} $\tau(\beta)$ for $\beta\in\bigcup_{\cs\in \La}E^\cs$. Recall from Lemma  \ref{lem:geometry-basic-interval-2} and  Proposition \ref{lem:geometry-basic-interval-4} that for each $\cs\in\La$ the map $\Psi_\cs$ bijectively maps the {exceptional} set $E=(1,2]\setminus\bigcup_{\s\in\F}J^\s$ to the relative {exceptional} set $E^\cs= (J^\cs\setminus I^\cs)\setminus\bigcup_{\r\in\F}J^{\cs\bullet\r}$.

\begin{lemma} \label{lem:beta-in-E}
Let $\hat{\beta}\in E\backslash\{2\}$ with $\de(\hat{\beta})=\de_1\de_2\ldots$. Also let $\cs\in\La$, and set $\beta:=\Psi_\cs(\hat{\beta})$. Then
\begin{enumerate}[{\rm(i)}]
\item  $b(\tau(\hat{\beta}),\hat{\beta})=0\de_2\de_3\ldots$; and
\item the map $\hat{t}\mapsto \big(\Phi_\cs(b(\hat{t},\hat{\beta}))\big)_{\beta}$ is continuous at $\tau(\hat{\beta})$.
\end{enumerate}
\end{lemma}

\begin{proof}
First we prove (i). Note by Proposition \ref{lem:critical-non-Fareyinterval} that $\tau(\hat\beta)=1-{1/\hat\beta}=(0\de_2\de_3\ldots)_{\hat\beta}$. {So} by Lemma \ref{lem:greedy-expansion} it suffices to {verify} that
\begin{equation}
  \label{eq:greedy-expansion}
  \si^n(0\de_2\de_3\ldots)\prec \de_1\de_2\ldots\quad\forall ~n\ge 0.
\end{equation}
By Lemma \ref{lem:quasi-greedy expansion}, it is immediate that $\si^n(0\de_2\de_3\ldots)\lle {(\de_i)}$ for all $n\geq 0$. If equality holds for some $n$, then $\de(\hat{\beta}){=(\de_i)}$ is periodic with period $m\geq 2$ (since $\hat{\beta}\neq 2$), so by Lemma \ref{lem:periodic-qg-expansion}, $\de(\hat{\beta})=\L(\r)^\f$ for some Lyndon word $\r$. This implies $\hat{\beta}={\beta_\ell^{\r}}\in J^{\r}$. But then by Lemma \ref{lem:farey-interval}, $\hat{\beta}\in J^\s$ for some Farey word $\s$, and so $\hat{\beta}\not\in E$, a contradiction. This proves (\ref{eq:greedy-expansion}), and then yields statement (i).

For (ii), note by Lemma \ref{lem:greedy-expansion} (ii) that the map $\hat{t}\mapsto b(\hat{t},\hat{\beta})$ is continuous at all points $\hat{t}$ {for which} $b(\hat t,\hat \beta)$ {does not end} with $0^\f$. Furthermore, the map $\Phi_\cs$ is continuous with respect to the order topology. So, by statement (i) it follows that the map $\hat{t}\mapsto \big(\Phi_\cs(b(\hat{t},\hat{\beta}))\big)_{\beta}$ is continuous at $\tau(\hat{\beta})$, completing the proof.
\end{proof}

\begin{proposition}  \label{prop:cric-exception-1}
   Let  $\cs\in\La$. Then for any $\beta\in E^\cs$ we have
   \[
   \tau(\beta)= \Big(\Phi_\cs{\big(0\de_2\de_3\ldots\big)}\Big)_\beta,
   \]
   where $1\de_2\de_3\ldots$ is the quasi-greedy expansion of $1$ in base $\hat\beta:=\Psi_\cs^{-1}(\beta)$.
	 \end{proposition}

\begin{proof}
Let $\beta\in E^\cs$. Then by Lemma  \ref{lem:geometry-basic-interval-2} and Proposition \ref{lem:geometry-basic-interval-4} there exists a unique $\hat\beta\in E$ such that ${\hat\beta}=\Psi_\cs^{-1}(\beta)\in E$, in other words, $\de(\beta)=\Phi_\cs(\de({\hat\beta}))$. Write $\de({\hat\beta})=\de_1\de_2\ldots$, and set $t^*:=(\Phi_\cs(0\de_2\de_3\ldots))_\beta$. We will show that $\tau(\beta)=t^*$, by proving that $h_{top}(\K_\beta(t))>0$ for $t<t^*$, and $\K_\beta(t)$ is countable for $t>t^*$. We consider separately the two cases: (i) $\hat{\beta}<2$ and (ii) $\hat{\beta}=2$.

\medskip

{\bf Case I.} $\hat{\beta}<2$. First, for notational convenience, we define the map
\[
\Theta_{\cs,\hat{\beta}}: (0,1)\to (0,1); \quad \hat{t}\mapsto \big(\Phi_\cs(b(\hat{t},\hat{\beta}))\big)_{\beta}.
\]
Since by Lemma \ref{lem:beta-in-E} the map $\Theta_{\cs,\hat{\beta}}$ is continuous at $\tau(\hat{\beta})$ and $t^*=\Theta_{\cs,\hat{\beta}}(\tau(\hat\beta))=(\Phi_\cs(0\de_2\de_3\ldots))_\beta$, it is, by the monotonicity of the set-valued map $t\mapsto \K_\beta(t)$, sufficient to prove the following two things:
\begin{equation} \label{eq:first-thing}
\hat{t}<\tau(\hat{\beta}) \qquad \Longrightarrow \qquad h_{top}\big(\K_\beta(\Theta_{\cs,\hat{\beta}}(\hat{t}))\big)>0,
\end{equation}
and
\begin{equation} \label{eq:second-thing}
\hat{t}>\tau(\hat{\beta}) \qquad \Longrightarrow \qquad \K_\beta(\Theta_{\cs,\hat{\beta}}(\hat{t}))\ \mbox{is countable}.
\end{equation}

First, take $\hat{t}<\tau(\hat{\beta})$ and set $t:=\Theta_{\cs,\hat{\beta}}(\hat{t})=\big(\Phi_\cs(b(\hat t, \hat \beta))\big)_\beta$. Since $\si^n(b(\hat{t},\hat{\beta}))\prec\de(\hat{\beta})$ for all $n\geq 0$, Lemma \ref{lem:cri-exception-1} implies that
\[
\si^n\big(\Phi_\cs(b(\hat{t},\hat{\beta}))\big)\prec \Phi_\cs(\de(\hat{\beta}))=\de(\beta) \quad\forall n\geq 0.
\]
Hence, $b(t,\beta)=\Phi_\cs(b(\hat{t},\hat{\beta}))$. Now
\begin{align*}
\Phi_\cs\big(\K_{\hat\beta}(\hat{t})\big)&=\{\Phi_\cs((x_i)): b(\hat{t},\hat{\beta})\lle \si^n((x_i))\prec\de(\hat{\beta})\quad\forall n\geq 0\}\\
&\subset \{\Phi_\cs((x_i)): \Phi_\cs\big(b(\hat{t},\hat{\beta})\big)\lle \si^n\big(\Phi_\cs((x_i))\big)\prec\Phi_\cs(\de(\hat{\beta})) \quad\forall n\geq 0\}\\
&=\{\Phi_\cs((x_i)): b(t,\beta)\lle \si^n\big(\Phi_\cs((x_i))\big)\prec \de(\beta) \quad\forall n\geq 0\}\\
&\subset \{(y_i): b(t,\beta)\lle \si^n((y_i)) \prec \de(\beta) \quad\forall n\geq 0\}=\K_{\beta}(t),
\end{align*}
where the first inclusion again follows by Lemma \ref{lem:cri-exception-1}. We deduce that
\[
h_{top}\big(\K_{\beta}(t)\big)\geq h_{top}\big(\Phi_\cs(\K_{\hat\beta}(\hat{t}))\big)=|\cs|^{-1}h_{top}\big(\K_{\hat\beta}(\hat{t}))>0,
\]
where the last inequality follows since $\hat{t}<\tau(\hat{\beta})$. This gives \eqref{eq:first-thing}.

Next, let $\hat{t}>\tau(\hat{\beta})$ and set $t:=\Theta_{\cs,\hat{\beta}}(\hat{t})$. Then by the same argument as above we have $b(t,\beta)=\Phi_\cs(b(\hat t,\hat \beta))$. Since $\hat{\beta}\in E$, there exists a sequence of Farey intervals $J^{\r_k}=[\beta_\ell^{\r_k}, \beta_r^{\r_k}]$ with $\r_k\in\F$ such that $\hat q_k:=\beta_\ell^{\r_k}\searrow \hat\beta$ as $k\to\f$.

We claim that $b(\hat{t},\hat{\beta})\succ (\r_k)^\f$ for all sufficiently large $k$. This can be seen as follows. As explained in the proof of Lemma \ref{lem:beta-in-E}, $\de(\hat{\beta})$ is not periodic, and therefore by Lemma \ref{lem:quasi-greedy expansion} (ii) the map $\beta'\mapsto \de(\beta')$ is continuous at $\hat{\beta}$ (where we use $\beta'$ to denote a generic base). This implies $\de(\hat{q}_k)\searrow \de(\hat{\beta})=1\de_2\de_3\ldots$ as $k\to\f$. But $\de(\hat{q}_k)=\de(\beta_\ell^{\r_k})=\L(\r_k)^\f$, and by Lemma \ref{lem:Farey-property} $\L(\r_k)$ is the word obtained from $\r_k$ by flipping the first and last digits. Thus, $(\r_k)^\f$ converges to $0\de_2\de_3\ldots$ in the order topology. Since $\hat{t}>\tau(\hat{\beta})$ implies $b(\hat{t},\hat{\beta})\succ b(\tau(\hat{\beta}),\hat{\beta})=0\de_2\de_3\ldots$, the claim follows.

We can now deduce that, for all sufficiently large $k$,
\begin{align*}
\K_\beta(t)&=\set{(y_i):b(t,\beta)\lle\si^n((y_i))\prec \de(\beta)~\forall n\ge 0}\\
  &=\set{(y_i):\Phi_\cs(b(\hat t, \hat\beta))\lle\si^n((y_i))\prec \Phi_\cs(\de(\hat\beta))~\forall n\ge 0}\\
  &\subset \set{(y_i): \Phi_\cs((\r_k)^\f)\lle \si^n((y_i))\prec \Phi_\cs(\L(\r_k)^\f)~\forall n\ge 0}\\
	&=\set{(y_i): (\cs\bullet\r_k)^\f\lle \si^n((y_i))\prec \L(\cs\bullet\r_k)^\f~\forall n\ge 0},
\end{align*}
where the inclusion follows using the claim and $\de(\hat{\beta})\prec\de(\hat{q}_k)=\L(\r_k)^\f$.
Hence, $\K_\beta(t)$ is countable by Proposition \ref{prop:countable}. This establishes \eqref{eq:second-thing}.

\medskip
{\bf Case II.} $\hat{\beta}=2$. In this case, $\de(\hat{\beta})=1^\f$, so $\de(\beta)=\Phi_\cs(\de(\hat{\beta}))=\L(\cs)^+\cs^\f$ and $t^*=\big(\Phi_\cs(01^\f)\big)_\beta=(\cs^-\L(\cs)^+\cs^\f)_\beta=(\cs 0^\f)_\beta$.
{Recall} that $\beta=\beta_r^\cs$ is the right endpoint of the Lyndon interval $J^\cs$.

If $t<t^*$, then $b(t,\beta)\prec b(t^*,\beta)=\cs 0^\f$, so by Lemma \ref{lem:greedy-expansion} {(iii)} there exists $k\in \N$ such that $b(t,\beta)\lle \cs^-\L(\cs)^+\cs^k 0^\f$. It follows that
\begin{align*}
\K_\beta(t)&=\{(x_i): b(t,\beta)\lle \si^n((x_i))\prec\de(\beta)\ \forall n\geq 0\}\\
&\supset \big\{\L(\cs)^+\cs^k\cs^-,\L(\cs)^+\cs^{k+1}\cs^-\big\}^\N,
\end{align*}
and hence $h_{top}(\K_\beta(t))>0$.

Now suppose $t>t^*$. Then {$b(t,\beta)\succ b(t^*,\beta)=\cs 0^\f$, so by Lemma \ref{lem:simplification-K-beta}},
\begin{align*}
\K_\beta(t)&=\{(x_i): b(t,\beta)\lle \si^n((x_i))\prec\L(\cs)^+\cs^\f\ \forall n\geq 0\}\\
&\subset \{(x_i): \cs 0^\f\lle  \si^n((x_i)){\prec} \L(\cs)^+\cs^\f\ \forall n\geq 0\}\\
&=\{(x_i): \cs^\f\lle \si^n((x_i))\lle \L(\cs)^+\cs^\f\ \forall n\geq 0\}\\
&= \{(x_i): \cs^\f\lle \si^n((x_i))\lle \L(\cs)^\f\ \forall n\geq 0\},
\end{align*}
where the second equality follows by using $\cs\in\LL$, {so $\si^n((x_i))\lge \cs 0^\f$ for all $n\ge 0$ if and only if} $\si^n((x_i))\lge \cs^\f$ for all $n\ge 0$.
Therefore, $\K_\beta(t)$ is countable by Proposition \ref{prop:countable}. This completes the proof.
\end{proof}

Next we will determine the critical value $\tau(\beta)$ for $\beta\in {E_\f}$. Recall from (\ref{eq:E_s}) that
 \[
 {E_\f=\bigcap_{n=1}^\f\bigcup_{\cs\in\La(n)}J^\cs,}
 \]
 where for each $n\in\N$ the Lyndon intervals $J^\cs, \cs\in{\La}(n)$ are pairwise disjoint. {Thus,} for any $\beta\in {E_\f}$ there exists a unique sequence of words $(\s_k)$ {with each $\s_k\in\F$} such that
\[
\set{\beta}=\bigcap_{n=1}^\f J^{\s_1\bullet\cdots\bullet\s_n}.
\]
We call $(\s_k)$ the \emph{coding} of $\beta$.

\begin{proposition}  \label{prop:cri-exception-2}
For any $\beta\in {E_\f}$ with its coding $(\s_k)$ we have
\[
\tau(\beta) =\lim_{n\to\f} (\s_1\bullet\cdots\bullet\s_n 0^\f)_\beta.
\]
\end{proposition}

\begin{proof}
Take $\beta\in E_\f$.  For $k\ge 1$ let $\cs_k:=\s_1\bullet\cdots\bullet\s_k$, and write $t_k:=(\cs_k 0^\f)_\beta$. Note that $\beta\in J^{\cs_k}=[\beta_\ell^{\cs_k}, \beta_r^{\cs_k}]$ for all $k\ge 1$. Hence
\begin{equation} \label{eq:6-2}
\de(\beta)\succ\de(\beta_\ell^{\cs_k})=\L(\cs_k)^\f,
\end{equation}
which implies that $b(t_k,\beta)=\cs_k 0^\f$ for all $k\ge 1$.  Observe that $\cs_{k+1}=\cs_k\bullet\s_{k+1}$ begins with $\cs_k^-$. Therefore,
  \[
  t_{k+1}= (\cs_{k+1}0^\f)_\beta< (\cs_k 0^\f)_\beta=t_k,
  \]
  so the sequence $(t_k)$ is decreasing. Since $t_k\ge 0$ for all $k\ge 1$, the limit
 $t^*:=\lim_{k\to\f}t_k$
 exists.  We will now show that $\tau(\beta)=t^*$.

  First we prove $\tau(\beta)\le t^*$. Since $t_k$ decreases to $t^*$ as $k\to\f$, it suffices to prove that $\tau(\beta)\le t_k$ for all $k\ge 1$.
 Let $q_k:=\beta_r^{\cs_k}$ for all $k\geq 1$. Then $q_k>\beta$ since $\beta\in J^{\cs_k}$, and
$q_k\searrow \beta$ as $k\to\f$. Set $t_k':=(\cs_k 0^\f)_{q_k}$. Since $q_k>\beta$, one can verify that $b(t_k', q_k)=\cs_k0^\f=b(t_k,\beta)$.  So,
  \begin{equation} \label{eq:kdr-2}
  \begin{split}
    \K_\beta(t_k)&=\set{(x_i):b(t_k,\beta)\lle\si^n((x_i))\prec \de(\beta)~\forall n\ge 0}\\
    &\subset\set{(x_i): \cs_k 0^\f \lle\si^n((x_i))\prec \de(q_k)~\forall n\ge 0}=\K_{q_k}(t_k').
  \end{split}
  \end{equation}
Note by Proposition \ref{prop:cric-exception-1} and Case II of its proof that $\tau(q_k)=(\cs_k 0^\f)_{q_k}=t_k'$.
  This implies that $\dim_H K_{q_k}(t_k')=0$, and thus by (\ref{eq:kdr-2}) we have $\dim_H K_\beta(t_k)=0$. {Hence} $\tau(\beta)\le t_k$ for any $k\ge 1$. Letting $k\to\f$ we obtain that $\tau(\beta)\le t^*$.

Next we prove $\tau(\beta)\ge t^*$. Note that $\beta=\Psi_{\cs_k}(\beta_k)$, where $\beta_k\in E_\f$ has coding $(\s_{k+1},\s_{k+2},\dots)$. {Let $a(\hat t,\beta_k)$ denote the quasi-greedy expansion of $\hat t$ in base $\beta_k$ (cf.~\cite[Lemma 2.3]{DeVries-Komornik-2011}). Observe that the map $\hat t\mapsto a(\hat t, \beta_k)$ is {strictly increasing and} left continuous everywhere in $(0,1)$, and {thus}} the map $\hat t\mapsto (\Phi_{\cs_k}(a(\hat t,\beta_k)))_\beta$ is also left continuous in $(0,1)$. So, by the same argument as in the proof of (\ref{eq:first-thing}) it follows that
\[
\tau(\beta)\geq \big(\Phi_{\cs_k}\big(a(\tau(\beta_k),\beta_k)\big)\big)_\beta\geq \big(\Phi_{\cs_k}(0^\f)\big)_\beta>(\cs_k^-0^\f)_\beta
\]
for every $k\in\N$, and letting $k\to\f$ gives $\tau(\beta)\geq t^*$.
\end{proof}

To illustrate Proposition \ref{prop:cri-exception-2} we construct in each Farey interval $J^\s$ a transcendental base $\beta\in E_\f$ and give {an} explicit formula for the critical value $\tau(\beta)$. Recall from \cite{Allouche_Shallit_1999} that the classical \emph{Thue-Morse sequence} $(\theta_i)_{i=0}^\f=01101001\ldots$ is defined recursively as follows. Let $\theta_0=0$; and if $\theta_0\ldots \theta_{2^n-1}$ is defined for some $n\ge 0$, then
\begin{equation} \label{eq:Thue-Morse-definition}
\theta_{2^n}\ldots \theta_{2^{n+1}-1}=\overline{\theta_0\ldots \theta_{2^n-1}}.
\end{equation}
 By the definition of $(\theta_i)$ it follows that
\begin{equation} \label{eq:propoerty-Thue-Morse-sequence}
  \theta_{2k+1}=1-\theta_k,\quad  \theta_{2k}=\theta_k\quad\textrm{for any }k\ge 0.
\end{equation}
{Komornik and Loreti \cite{KomornikLoreti1998} showed that}
\begin{equation} \label{eq:thue-morse-inequality}
  \theta_{i+1}\theta_{i+2}\ldots \prec \theta_1\theta_2\ldots\quad\forall i\ge 1.
\end{equation}

\begin{proposition} \label{prop:cri-exception-3}
  Given $\s=s_1\ldots s_m=0\c 1\in\F$, let {$\beta:=\beta_\f^\s\in (1,2]$} such that
  \[
  (\theta_1 \c\theta_2\,\theta_3 \c\theta_4\,\ldots\,\theta_{2k+1}\c\theta_{2k+2}\ldots)_{\beta}=1.
  \]
  Then $\beta\in E_\f\cap J^\s$ is transcendental, and
  \[
  \tau(\beta)=\frac{2\sum_{j=2}^m s_j\beta^{m-j}+\beta^{m-1}-\beta^m}{\beta^m-1}.
  \]
\end{proposition}

We point out that in the above proposition, $\c$ may be the empty word.
To prove the transcendence of $\beta$ we recall the following result {due to} Mahler \cite{Mahler_1976}.

\begin{lemma}[Mahler, 1976] \label{lem:mahler}
If $z$ is an algebraic number in the open unit disc, then the number
$
Z:=\sum_{i=1}^\f\theta_i z^i
$
is transcendental.
\end{lemma}

\begin{proof}[Proof of Proposition \ref{prop:cri-exception-3}]
Let $\s=s_1\ldots s_m=0\c 1\in\F$. First we prove that
\begin{equation}  \label{eq:exam-1}
  \de(\beta)=\theta_1 \c\theta_2\,\theta_3 \c\theta_4\,\ldots\,\theta_{2k+1}\c\theta_{2k+2}\ldots=:(\de_i).
\end{equation}
By Lemma \ref{lem:quasi-greedy expansion} it suffices to prove that
$\si^n((\de_i))\lle (\de_i)$ for all $n\ge 1$. Note by Lemma \ref{lem:Farey-property} that $\de_1\ldots \de_m=1\c 1=\L(\s)^+=:{a_1\ldots a_m^+}$. Take $n\in\N$, and write $n=mk+j$ with $k\in\N\cup\{0\}$ and {$j\in\{1,2,\dots,m\}$}. We will prove $\si^n((\de_i))\prec (\de_i)$ in the following three cases.

Case I. {$j\in\set{1,2,\ldots, m-2}$}. Note by (\ref{eq:propoerty-Thue-Morse-sequence}) that $\theta_{2k}=1-\theta_{2k+1}$. {This implies that} $\si^n((\de_i))$ begins with either $a_{j+1}\ldots a_m$ or $a_{j+1}\ldots a_{m}^+s_1\ldots s_{m-1}$. By Lemma \ref{lem:lyndon-equivalence} it follows that $\si^n((\de_i))\prec (\de_i)$.

Case II. {$j=m-1$}. Then $\si^n((\de_i))$ begins with $\theta_{2k}\theta_{2k+1}\c$ for some $k\in\N$. If $\theta_{2k}=0$, then it is clear that $\si^n((\de_i))\prec (\de_i)$ since $\de_1=1$. Otherwise, (\ref{eq:propoerty-Thue-Morse-sequence})  implies that $\theta_{2k}\theta_{2k+1}\c=10\c=1 s_1\ldots s_{m-1}$. Hence,
\[
\si^n((\de_i))=1 s_1\ldots s_{m-1}\,\de_{n+m+1}\de_{n+m+2}\ldots\prec 1 a_2\ldots a_m^+ \de_{m+1}\de_{m+2}=(\de_i),
\]
where the strict inequality follows since $s_1\ldots s_{m-1}\lle a_2\ldots a_m$.

Case III. {$j=m$}. Then
\begin{align*}
\si^n((\de_i))&={\theta_{2k+1}\c\theta_{2k+2}\,\theta_{2k+3}\c\theta_{2k+4}\ldots} \prec \theta_1\c\theta_2\,\theta_3\c\theta_4\ldots=(\de_i),
\end{align*}
where the strict inequality is a consequence of (\ref{eq:thue-morse-inequality}).

Therefore, by Cases I--III we establish (\ref{eq:exam-1}). Next we show that $\beta\in E_\f$. For $k\in\N$ let $\cs_k:=\s_1\bullet\s_2\bullet\cdots\bullet\s_k$ with $\s_1=\s$ and $\s_i=01$ for all $2\le i\le k$. Then $\cs_k\in\La$ for all $k\in\N$. So it suffices to show that $\beta\in J^{\cs_k}$ for all $k\ge 1$. First we claim that
\begin{equation} \label{eq:exam-2}
  \cs_k=\overline{\theta_1}\c\overline{\theta_2}\,\overline{\theta_3}\c\overline{\theta_4}\ldots \overline{\theta_{2^k-1}}\c\overline{\theta_{2^k}}^+,\qquad \L({\cs_k})=\theta_1\c\theta_2\,\theta_3\c\theta_4\,\ldots\theta_{2^k-1}\c\theta_{2^k}^-
\end{equation}
for all $k\ge 1$.

Since $\cs_1=\s=0\c 1=\overline{\theta_1}\c\overline{\theta_2}^+$ and $\L(\cs_1)=1\c 0=\theta_1\c\theta_2^-$, (\ref{eq:exam-2}) holds for $k=1$. Now suppose (\ref{eq:exam-2}) holds for a given $k\in\N$. Then
\begin{align*}
  \cs_{k+1}=\cs_k\bullet(01) =\cs_k^-\L(\cs_k)^+&=\overline{\theta_1}\c\overline{\theta_2}\ldots \overline{\theta_{2^k-1}}\c\overline{\theta_{2^k}}\;\theta_1\c\theta_2\,\ldots\theta_{2^k-1}\c\theta_{2^k}\\
  &=\overline{\theta_1}\c\overline{\theta_2}\ldots\overline{\theta_{2^{k+1}-1}}\c\overline{\theta_{2^{k+1}}}^+,
\end{align*}
where the last equality follows since, by the definition of $(\theta_i)$ in \eqref{eq:Thue-Morse-definition}, $\theta_{2^k+1}\ldots \theta_{2^{k+1}}=\overline{\theta_1\ldots \theta_{2^k}}^+$. Similarly, by the induction hypothesis and Lemma \ref{lem:lyndon-commutative} we obtain
\begin{align*}
  \L(\cs_{k+1})=\L(\cs_k\bullet (01))=\cs_k\bullet\L(01)&=\cs_k\bullet(10)=\L(\cs_k)^+\cs_k^-\\
  &=\theta_1\c\theta_2\,\ldots\theta_{2^k-1}\c\theta_{2^k}\,\overline{\theta_1}\c\overline{\theta_2}\ldots \overline{\theta_{2^k-1}}\c\overline{\theta_{2^k}}\\
  &=\theta_1\c\theta_2\,\ldots\theta_{2^{k+1}-1}\c\theta_{2^{k+1}}^-.
\end{align*}
Hence, by induction (\ref{eq:exam-2}) holds for all $k\ge 1$.

Next, recall that the Lyndon interval $J^{\cs_k}=[\beta_\ell^{\cs_k}, \beta_r^{\cs_k}]$ satisfies
\begin{equation}  \label{eq:exam-3}
  \de(\beta_\ell^{\cs_k})=\L(\cs_k)^\f\quad\textrm{and}\quad \de(\beta_r^{\cs_k})=\L(\cs_k)^+\cs_k^\f.
\end{equation}
By (\ref{eq:exam-1}) and (\ref{eq:exam-2}) it follows that $\de(\beta)$ begins with $\L(\cs_k)^+$, so by \eqref{eq:exam-3}, $\de(\beta)\succ\L(\cs_k)^\f=\de(\beta_\ell^{\cs_k})$. Thus, $\beta>\beta_\ell^{\cs_k}$ for all $k\ge 1$. On the other hand, by (\ref{eq:exam-1}), (\ref{eq:exam-2}), and Lemmas \ref{lem:lyndon-commutative} and \ref{lem:lyndon-associative} we see that $\de(\beta)$ also begins with
\begin{align*}
\L(\cs_{k+2})^+&=\L\big(\cs_k\bullet(01\bullet 01)\big)^+=\L\big(\cs_k\bullet(0011)\big)^+\\
&=\big(\cs_k\bullet\L(0011)\big)^+=\big(\cs_k\bullet (1100)\big)^+ =\L(\cs_k)^+\cs_k\cs_k^-\L(\cs_k)^+,
\end{align*}
which is strictly smaller than a prefix of $\de(\beta_r^{\cs_k})^+={\L(\cs_k)^+}\cs_k^\f$. This implies that $\beta<\beta_r^{\cs_k}$ for all $k\ge 1$. Hence, $\beta\in J^{\cs_k}$ for all $k\ge 1$, and thus $\beta\in E_\f\cap J^\s$.
Furthermore, by (\ref{eq:exam-1}), (\ref{eq:exam-2}) and Proposition \ref{prop:cri-exception-2} it follows that
\begin{align*}
  \tau({\beta})&=\lim_{k\to\f}(\cs_k 0^\f)_{\beta}=\big(\overline{\theta_1}\c\overline{\theta_2}\,\overline{\theta_3}\c\overline{\theta_4}\ldots\big)_{\beta}\\
  &=\sum_{k=0}^\f\left(\frac{1}{\beta^{m k+1}}+2\sum_{j=2}^{m-1}\frac{s_j}{\beta^{m k+j}}+\frac{1}{\beta^{m k+m}}\right)-({\theta_1}\c{\theta_2}{\theta_3}\c\theta_4\ldots)_{\beta}\\
  &=\frac{\beta^{m-1}+2\sum_{j=2}^{m-1}s_j\beta^{m-j}+1}{\beta^m-1}-1\\
  &=\frac{2\sum_{j=2}^m s_j\beta^{m-j}+\beta^{m-1}-\beta^m}{\beta^m-1},
\end{align*}
where we recall that $\c=s_2\ldots s_{m-1}$, and the last equality uses that $s_m=1$.

Finally,  the transcendence of $\beta$ follows by  using (\ref{eq:exam-1}), Lemma \ref{lem:mahler} and a similar  argument as in the proof of \cite[Proposition 5.2]{Kong-Li-2020}.
\end{proof}

\begin{remark}\label{rem:komornik-Loreti-transcendence}\mbox{}

\begin{itemize}
 \item When $\s=01$, the base $\beta_\f^{01}\approx 1.78723$ given in Proposition \ref{prop:cri-exception-3} is the \emph{Komornik-Loreti constant} (cf.~\cite{KomornikLoreti1998}), whose transcendence was first proved by Allouche and Cosnard \cite{Allouche_Cosnard_2000}. {In this case we obtain $\tau(\beta_\f^{01})=(2-\beta_\f^{01})/(\beta_\f^{01}-1)\approx .270274$.}
 \item When $\s=001$, the base $\beta_\f^{001}\approx 1.55356$ is a critical value for the fat Sierpinski gaskets studied by Li and the second author in \cite{Kong-Li-2020}. In this case we have $\tau(\beta_\f^{001})\approx .241471$.
\end{itemize}
\end{remark}

By Proposition \ref{prop:cri-exception-3} and numerical calculation we give in Table \ref{table-1} the triples $(\s, \beta_\f^\s, \tau(\beta_\f^\s))$ for all $\s\in F_3^*\subset\F$.
Based on Proposition \ref{prop:cri-exception-3} we conjecture that each base $\beta\in E_\f$ is transcendental.

\begin{center}
\begin{table}[h!]
\begin{tabular}{|c|c|c|c |c|c|c|c|}
  \hline
  $\s=$&  0001 & 001 & 00101   &01& 01011&011& 0111  \\\hline
 $\beta_\f^\s\approx$& 1.43577 & 1.55356 &1.59998& 1.78723 & 1.83502 &1.91988 & 1.96452 \\\hline
 $\tau(\beta_\f^\s)\approx$& 0.218562& 0.241471 &0.336114 & 0.270274 & 0.432175 & 0.40305 & 0.455933 \\
  \hline
\end{tabular}

\bigskip
\caption{The triples $(\s, \beta_\f^\s, \tau(\beta_\f^\s))$ with $\s\in F_3^*\subset\F$.}
\label{table-1}
\end{table}
\end{center}

\section{C\`adl\`ag {property} of the critical value function} \label{sec:continuity}

In this section we prove Proposition \ref{prop:discontinuities} and Theorem \ref{main:critical-devils-staircase}. Recall by Lemma \ref{lem:presentation-J-s} that the interval $(1,2]$ can be partitioned as
\begin{equation} \label{eq:partition-2}
(1,2]=E\cup\bigcup_{\s\in\F}J^\s=E \cup E_\f\cup\bigcup_{\cs\in\La}E^\cs\cup\bigcup_{\cs\in\La}I^\cs.
\end{equation}
Here we emphasize that the {exceptional} set $E$, the relative {exceptional} sets $E^\cs$ and the {infinitely Farey set} $E_\f$ featuring in (\ref{eq:partition-2}) all have Lebesgue measure zero in view of Propositions \ref{prop:geometry-basic-interval-2} and \ref{prop:dim-exceptional-set}. Hence the basic intervals $I^\cs, \cs\in\La$, and then certainly the Lyndon intervals $J^\cs, \cs\in\La$, are dense in $(1,2]$. This allows for approximation of points in $E, E^\cs$ and $E_\f$ by left and/or right endpoints of such Lyndon intervals.
We also recall from Theorem \ref{main:critical-basic-intervals} that for any basic interval $I^\cs=[\beta_\ell^\cs, \beta_*^\cs]$ with $\de(\beta_\ell^\cs)=\L(\cs)^\f$ and $\de(\beta_*^\cs)=\L(\cs)^+ \cs^-\L(\cs)^\f$, the critical value is given by
\begin{equation} \label{eq:cri-basic-interval}
\tau(\beta)=(\Phi_\cs(0^\f))_\beta=(\cs^-\L(\cs)^\f)_\beta\quad\textrm{for any}\quad \beta\in I^\cs.
\end{equation}
{Moreover}, by Proposition \ref{prop:cric-exception-1} it follows that for each $\beta\in E^\cs$ we have
\begin{equation} \label{eq:cri-relative-bifur}
\tau(\beta)={(\Phi_\cs(0\de_2\de_3\ldots))_\beta},
\end{equation}
where $1\de_2\de_3\ldots=\de({\hat\beta})$ with ${\hat\beta}={\Psi_\cs}^{-1}(\beta)\in E$. In particular, when $\beta\in E$ we have $\tau(\beta)=1-\frac{1}{\beta}$ (see Proposition \ref{lem:critical-non-Fareyinterval}). When $\beta\in E_\f$, {it follows by Proposition \ref{prop:cri-exception-2} that}
\begin{equation} \label{eq:cri-infinite-bifur}
\tau(\beta)=\lim_{n\to\f}(\s_1\bullet\s_2\bullet\cdots\bullet\s_n 0^\f)_\beta,
\end{equation}
where $(\s_k)$ is the unique coding of $\beta$ (that is $\beta\in J^{\s_1\bullet\cdots\bullet\s_k}$ for all $k\in\N$).

From (\ref{eq:cri-basic-interval}) it is clear that the critical value function $\tau$ is continuous inside each basic interval $I^\cs=[\beta_\ell^\cs, \beta_*^\cs]$. So, in view of (\ref{eq:partition-2}), we still need to consider the continuity of $\tau$ for $\beta\in {E\cup E_\f\cup\bigcup_{\cs\in \La}(E^\cs\setminus\set{\beta_r^\cs})}$, the left continuity of $\tau$ at $\beta=\beta_\ell^\cs$ and $\beta=\beta_r^\cs$, and the right continuity of $\tau$ at $\beta=\beta_*^\cs$.
We need the following lemma:

\begin{lemma} \label{lem:no-periodic-deltas}
If $\beta\in(1,2)$ and $\de(\beta)$ is periodic, then $\beta\in\bigcup_{\cs\in\La} I^\cs$.
\end{lemma}

\begin{proof}
Assume $\de(\beta)$ is periodic. In view of (\ref{eq:partition-2}) it suffices to prove
\[
\beta\notin E\cup E_\f\cup\bigcup_{\cs\in\La}E^\cs.
\]
 First we prove $\beta\notin E$. By Lemma \ref{lem:periodic-qg-expansion}, $\de(\beta)=\L(\cs')^\f$ for some Lyndon word $\cs'$. This means $\beta$ is the left endpoint of a Lyndon interval, so by Lemma \ref{lem:farey-interval}, $\beta\in J^\s$ for some Farey word $\s$. Hence $\beta\not\in E$.

Next, suppose $\beta\in E^\cs$ for some $\cs\in\La$. {Clearly, $\beta\ne \beta_r^\cs$ since $\de(\beta_r^\cs)=\L(\cs)^+\cs^\f$ is not periodic. Thus $\beta\in (\beta_*^\cs, \beta_r^\cs)\setminus\bigcup_{r\in\F}J^{\cs\bullet\r}$. So, by Proposition \ref{prop:geometry-basic-interval-2}} there is a sequence $(\r_k)$ of Farey words such that $\beta_\ell^{\cs\bullet\r_k}\searrow \beta$. {Write $\de(\beta)=(a_1\ldots a_n)^\f$ with minimal period $n$. Then} by Lemma \ref{lem:quasi-greedy expansion} {it follows that} $\de(\beta_\ell^{\cs\bullet\r_k})\searrow a_1\ldots a_n^+ 0^\f$, so for all sufficiently large $k$, $\de(\beta_\ell^{\cs\bullet\r_k})$ contains a block of more than $2m$ zeros, where $m:=|\cs|$. But this is impossible, since $\de(\beta_\ell^{\cs\bullet\r_k})=\L(\cs\bullet \r_k)^\f$ is a concatenation of blocks from $\cs, \cs^-, \L(\cs)$ and $\L(\cs)^+$. These blocks all have length $m$, and only $\cs^-$ could possibly consist of all zeros, while $\cs^-$ can only be followed by $\L(\cs)$ or $\L(\cs)^+$. Thus, $\de(\beta_\ell^{\cs\bullet\r_k})$ cannot contain a block of $2m$ zeros. This contradiction shows that $\beta\not\in E^\cs$.

Finally, suppose $\beta\in E_\f$. Then there is a sequence $(\cs_k)$ of words in $\La$ such that $\beta$ lies in the interior of $J^{\cs_k}$ for each $k$. Note that {$\beta=\beta_\ell^{\cs'}$} is the left endpoint of $J^{\cs'}$. Thus, $J^{\cs'}\cap J^{\cs_k}\neq\emptyset$, and therefore by Lemma \ref{lem:farey-interval}, it must be the case that $J^{\cs'}\subset J^{\cs_k}$ for all $k$. But this is impossible, since $|J^{\cs_k}|\to 0$. Hence, $\beta\not\in E_\f$.
\end{proof}

\begin{proof}[Proof of Proposition \ref{prop:discontinuities}]
First
fix $\beta_0\in(1,2]\setminus\set{\beta_r^\cs: \cs\in\La}$. It is sufficient to prove that

\bigskip
$(*)$ for each $N\in\N$ there exists $r>0$ such that if $\beta\in(1,2]$ satisfies $|\beta-\beta_0|<r$, then

\ \ \ \ \ there is a word
$s_1\ldots s_N$ such that $\tau(\beta)$ has a $\beta$-expansion beginning with $s_1\ldots s_N$,

\ \ \ \ \ and $\tau(\beta_0)$
has a $\beta_0$-expansion beginning with $s_1\ldots s_N$.

\bigskip
For, if $\tau(\beta)=(s_1\ldots s_N c_1 c_2\ldots)_\beta$ and $\tau(\beta_0)=(s_1\ldots s_N d_1 d_2\ldots)_{\beta_0}$, then
\begin{align*}
|\tau(\beta)-\tau(\beta_0)|&\leq \big|(s_1\ldots s_N c_1 c_2\ldots)_\beta-(s_1\ldots s_N c_1 c_2\ldots)_{\beta_0}\big|\\
&\qquad + \big|(s_1\ldots s_N c_1 c_2\ldots)_{\beta_0}-(s_1\ldots s_N d_1 d_2\ldots)_{\beta_0}\big|\\
&\leq \sum_{i=1}^\f \left|\frac{1}{\beta^i}-\frac{1}{\beta_0^i}\right|+\sum_{i=N+1}^\f \frac{1}{\beta_0^i}\\
&=\frac{|\beta-\beta_0|}{(\beta-1)(\beta_0-1)}+\frac{1}{(\beta_0-1)\beta_0^N}\\
&<\frac{r}{(\beta-1)(\beta_0-1)}+\frac{1}{(\beta_0-1)\beta_0^N},
\end{align*}
and this can be made as small as desired by choosing $N$ sufficiently large and $r$ sufficiently small. In view of (\ref{eq:partition-2}) we prove $(*)$ by considering several cases.

\medskip
{\bf Case I.} $\beta_0\in(\beta_\ell^\cs,\beta_*^\cs)$ for some basic interval $I^\cs=[\beta_\ell^\cs,\beta_*^\cs]$ with $\cs\in\La$. It is clear from Theorem \ref{main:critical-basic-intervals} that $(*)$ holds in this case.

\medskip
{\bf Case II.} $\beta_0\in E$. Then by Proposition \ref{prop:geometry-basic-interval-2} there exists a sequence of Farey intervals $J^{\s_k}=[\beta_\ell^{\s_k}, \beta_r^{\s_k}]$ such that $\beta_\ell^{\s_k}\to\beta_0$ as $k\to\f$. Furthermore, $|J^{\s_k}|\to 0$ as $k\to\f$. This implies that the length $|\s_k|$ of the Farey word $\s_k$ goes to infinity as $k\to\f$. Let $N\in\N$ be given. We can choose $r>0$ small enough so that if a Farey interval $J^\s$ intersects $(\beta_0-r,\beta_0+r)$, then $|\s|>N$ and
\begin{equation}  \label{eq:14-2}
    \de_1(\beta)\ldots \de_N(\beta)=\de_1(\beta_0)\ldots \de_N(\beta_0).
  \end{equation}
We can guaranteee \eqref{eq:14-2} because for $\beta_0\in E\setminus\set{2}$ the expansion $\de(\beta_0)$ is not periodic by Lemma \ref{lem:no-periodic-deltas}, so by Lemma \ref{lem:quasi-greedy expansion} the map $\beta\mapsto \de(\beta)$ is continuous at $\beta_0$. Furthermore, for $\beta_0=2$ the map $\beta\mapsto\de(\beta)$ is left continuous at $\beta_0$.

Let $\beta\in(\beta_0-r, \beta_0+r)$. By {(\ref{eq:partition-2})} we have either $\beta\in E$ or $\beta\in J^\s$ for some $\s\in\F$.
 If $\beta\in E$, then by Proposition \ref{lem:critical-non-Fareyinterval} it follows that $\tau(\beta)=1-\frac{1}{\beta}=(0\de_2(\beta)\de_3(\beta)\ldots)_\beta$ and
 \begin{equation}\label{eq:kkk-1}
 \tau(\beta_0)=1-\frac{1}{\beta_0}=\big(0\de_2(\beta_0)\de_3(\beta_0)\ldots\big)_{\beta_0},
 \end{equation}
so $(*)$ holds by \eqref{eq:14-2}.

 Next we assume $\beta\in J^\s$ with $\s=s_1\ldots s_m\in\F$. By our choice of $r$ it follows that $m=|\s|>N$, and (\ref{eq:14-2}) holds. Since $\beta\in J^\s=[\beta_\ell^\s, \beta_r^\s]$, we have
 $\L(\s)^\f\lle \de(\beta)\lle \L(\s)^+{\s^\f}$. Write $\L(\s)=a_1\dots a_m$; then by (\ref{eq:14-2}) it follows that
 \begin{equation}  \label{eq:14-4}
   \de_1(\beta_0)\ldots\de_N(\beta_0)=\de_1(\beta)\ldots \de_N(\beta)=a_1\ldots a_N.
 \end{equation}
 Observe by (\ref{eq:partition-2})--(\ref{eq:cri-infinite-bifur}) and Lemma \ref{lem:Farey-property} that $\tau(\beta)$ has a $\beta$ expansion beginning with ${\s^-}=0a_2\ldots a_m$. Hence, by (\ref{eq:14-4}) $\tau(\beta)$ has a $\beta$-expansion with prefix
 \begin{equation} \label{eq:14-5}
 s_1\ldots s_N=0a_2\ldots a_N=0\de_2(\beta_0)\ldots \de_N(\beta_0).
 \end{equation}
  This, together with (\ref{eq:kkk-1}), gives $(*)$.

\medskip
{\bf Case III.} $\beta_0\in E^\cs\backslash \{\beta_r^\cs\}$ for some $\cs\in\La$. The proof is similar to that of Case II, but there are some extra details involving the substitution operator.
By Proposition \ref{prop:geometry-basic-interval-2} it follows that
\begin{equation}  \label{eq:cont-0}
  {\beta_0}\in(\beta_*^\cs, \beta_r^\cs)\setminus\bigcup_{\r\in\F}J^{\cs\bullet\r},
  \end{equation}
  and there exists a sequence $(\r_k)$ of Farey words such that $\beta_{\ell}^{\cs\bullet\r_k}\to{\beta_0}$ as $k\to\f$. This implies that $|\r_k|\to \f$ as $k\to\f$.

  Let $N\in\N$ be given; without loss of generality we may assume that $N=N'|\cs|$ for some integer $N'$. By Lemma \ref{lem:no-periodic-deltas} we can choose $r>0$ sufficiently small so that if a Lyndon interval $J^{\cs\bullet\r}$ intersects $(\beta_0-r,\beta_0+r)$, then $|\r|>N'$ and
  \begin{equation}  \label{eq:cont-2}
    \de_1(\beta)\ldots \de_{N'|\cs|}(\beta)=\de_1(\beta_0)\ldots\de_{N'|\cs|}(\beta_0)=\cs\bullet(\de_1({\hat\beta_0})\ldots\de_{N'}({\hat\beta_0}))
  \end{equation}
  for any $\beta\in(\beta_0-r,\beta_0+r)$, where ${\hat\beta_0}=\Psi_\cs^{-1}(\beta_0)\in E$.

  Now take $\beta\in(\beta_0-r, \beta_0+r)$. By (\ref{eq:cont-0}) it follows that either $\beta\in E^\cs$ or $\beta\in J^{\cs\bullet\r}$ for some $\r\in\F$. If $\beta\in E^\cs$, we let ${\hat\beta}=\Psi_\cs^{-1}(\beta)$. Then (\ref{eq:cont-2}) yields
  \[
  \cs\bullet(\de_1(\hat\beta)\ldots\de_{N'}(\hat\beta))=\de_1(\beta)\ldots\de_{N'|\cs|}(\beta) =\cs\bullet(\de_1(\hat{\beta}_0)\ldots\de_{N'}(\hat{\beta}_0)),
  \]
  which implies $\de_1(\hat\beta)\ldots \de_{N'}(\hat\beta)=\de_1(\hat{\beta}_0)\ldots\de_{N'}(\hat{\beta}_0)$. So, by (\ref{eq:cri-relative-bifur}) it follows that $\tau(\beta)$ has a $\beta$-expansion with a prefix
  \[
  \cs\bullet(0\de_2(\hat\beta)\ldots \de_{N'}((\hat\beta))=\cs\bullet(0\de_2(\hat{\beta}_0)\ldots \de_{N'}(\hat{\beta}_0)),
  \]
  which coincides with a prefix of a $\beta_0$-expansion of $\tau(\beta_0)$. This gives $(*)$.
	
  If $\beta\in J^{\cs\bullet\r}$ for some $\r=r_1\ldots r_m\in\F$, then by our assumption on $r$ we have $m>N'$. Furthermore,
  \[
  (\cs\bullet\L(\r))^\f=\L(\cs\bullet\r)^\f\lle \de(\beta)\lle \L(\cs\bullet\r)^+(\cs\bullet\r)^\f=(\cs\bullet\L(\r))^+(\cs\bullet\r)^\f.
  \]
  Therefore, writing $\L(\r)=b_1\ldots b_m$, it follows from (\ref{eq:cont-2}) that
  \begin{equation}  \label{eq:cont-3}
    \cs\bullet(b_1\ldots b_{N'})=\de_1(\beta)\ldots\de_{N'|\cs|}(\beta)=\cs\bullet(\de_1(\hat{\beta}_0)\ldots \de_{N'}(\hat{\beta}_0)).
  \end{equation}
This shows that $b_1\ldots b_{N'}=\de_1(\hat{\beta}_0)\ldots\de_{N'}(\hat{\beta}_0)$. Now observe by (\ref{eq:partition-2})--(\ref{eq:cri-infinite-bifur}) and Lemma \ref{lem:Farey-property} that $\tau(\beta)$ has a $\beta$-expansion beginning with $(\cs\bullet\r)^-=\cs\bullet(0 b_2\ldots b_m)$. So, by (\ref{eq:cri-relative-bifur}) and (\ref{eq:cont-3}) it follows that $\tau(\beta)$ has a $\beta$-expansion with a prefix
  \[
  \cs\bullet(0b_2\ldots b_{N'})=\cs\bullet(0\de_2(\hat{\beta}_0)\ldots \de_{N'}(\hat{\beta}_0))
  \]
  of length $N$, which is also a prefix of a $\beta_0$-expansion of $\tau(\beta_0)$. This again gives $(*)$.

%
%
%
	
\medskip
{\bf Case IV.} $\beta_0=\beta_\ell^\cs$ for some $\cs\in\La$. Here the right continuity of $\tau$ at $\beta_0$ follows from Theorem \ref{main:critical-basic-intervals}. The left-continuity {can be seen as follows.
If $\cs\in\F$, then $\tau(\beta_0)$ has a $\beta_0$-expansion $0\de_2(\beta_\ell^\cs)\de_3(\beta_\ell^\cs)\ldots$, and the left continuity follows by the argument in Case II, using the left-continuity of the map $\beta\mapsto \de(\beta)$ at $\beta_0$. Otherwise, $\cs=\cs'\bullet \r$ for some $\cs'\in\La$ and $\r\in\F$, and $\tau(\beta_0)$ has a $\beta_0$-expansion $\Phi_{\cs'}(0\de_2(\beta_\ell^\r)\de_3(\beta_\ell^\r)\ldots)$. In this case the left continuity follows from the argument in Case III.}


\medskip
{\bf Case V.} $\beta_0=\beta_*^\cs$ for some $\cs\in\La$. Here the left continuity at $\beta_0$ follows from Theorem \ref{main:critical-basic-intervals}. The right-continuity can be seen as follows. First note that 
\[
\de(\beta_*^\cs)=\L(\cs)^+\cs^-\L(\cs)^\f, \qquad\mbox{and}\qquad \tau(\beta_*^\cs)=\big(\cs^-\L(\cs)^\f\big)_{\beta_*^\cs}. 
\]
Observe also that $\beta_*^\cs=\lim_{\hat{\beta}\searrow 1}\Psi_\cs(\hat{\beta})$, so by Lemma \ref{lem:quasi-greedy expansion}, $\Phi_\cs(\de(\hat{\beta}))\searrow \de(\beta_*^\cs)$ as $\hat{\beta}\searrow 1$.

Now let $N\in\N$ be given. As in Case III, we may assume that $N=N'|\cs|$ for some integer $N'$. We choose $r>0$ small enough so that, if $\r\in\F$ and $J^{\cs\bullet\r}$ intersects $(\beta_0,\beta_0+r)$, then $|\r|>N'$. Now take $\beta\in (\beta_0,\beta_0+r)$. If $\beta\in E^\cs$, then $\beta=\Psi_\cs(\hat{\beta})$ for some $\hat{\beta}\in E$, and $\de(\beta)=\Phi_\cs(\de(\hat{\beta}))$. Note that, if $\beta\searrow \beta_0$, then $\hat{\beta}\searrow 1$ and so $\de(\hat{\beta})\searrow 10^\f$. Thus, we may also assume $r$ is small enough so that $\de(\hat{\beta})$ begins with $10^{N'-1}$. Then $\tau(\beta)$ has a $\beta$-expansion beginning with $\Phi_\cs(0^{N'})=\cs^-\L(\cs)^{N'-1}$, which is also a prefix of length $N$ of a $\beta_0$-expansion of $\tau(\beta_0)$.

Similarly, if $\beta\in J^{\cs\bullet \r}$ then we may assume $r$ is small enough so that $\r$ begins with $0^{N'-1}$. Then $\tau(\beta)$ has a $\beta$-expansion beginning with $(\cs\bullet \r)^-$, and therefore beginning with $\cs^-\L(\cs)^{N'-1}$, and we conclude as above.

%

\medskip
{\bf Case VI.} $\beta_0\in E_\f$. Then there exists a sequence $(\s_k)$ of Farey words such that
 \begin{equation*} \label{eq:cont-4}
 \set{\beta} =\bigcap_{k=1}^\f J^{\cs_k},
 \end{equation*}
 where $\cs_k:=\s_1\bullet\s_2\bullet\cdots\bullet\s_k$.
 Note that $J^{\cs_k}\supset J^{\cs_{k+1}}$ for all $k\ge 1$, and $|\cs_k|\to\f$ as $k\to\infty$.

Let $N\in\N$ be given, and choose $k$ so large that $|\cs_k|>N$.
Choose $r>0$ sufficiently small so that $(\beta_0-r, \beta_0+r)\subset J^{\cs_k}$. Take $\beta\in(\beta_0-r, \beta_0+r)$. Then $\beta\in J^{\cs_k}$, so by (\ref{eq:partition-2})--(\ref{eq:cri-infinite-bifur}) it follows that $\tau(\beta)$ has
 a $\beta$-expansion beginning with $\cs_k^-$, which is also a prefix (of length at least $N$) of a $\beta_0$-expansion of $\tau(\beta_0)$. Hence, we obtain $(*)$.

 \medskip
Finally, we consider $\beta_0=\beta_r^\cs$ for $\cs\in\La$. The left continuity at $\beta_0$ (i.e. the analog of $(*)$ for $\beta\in(\beta_0-r,\beta_0)$) follows just as in Case III. The jump at $\beta_0$ (i.e. \eqref{eq:upward-jump}) can be seen as follows. {Since $\tau(\beta_0)=(\cs 0^\f)_{\beta_0}<(\cs^\f)_{\beta_0}$ by Proposition \ref{prop:cric-exception-1} (or rather, Case II of its proof),} it suffices to {show that}
\begin{equation} \label{eq:right-discontinuity}
\lim_{\beta\searrow\beta_0}\tau(\beta)=(\cs^\f)_{\beta_0}.
\end{equation}

First assume $\cs\in\F$. Then $\de(\beta_0)=\L(\cs)^+\cs^\f$, and by Lemma \ref{lem:Farey-property} it follows that $\cs^\f=0\de_2(\beta_0)\de_3(\beta_0)\ldots$. So, by the same argument as in Case II we obtain (\ref{eq:right-discontinuity}).

Next suppose $\cs=\cs'\bullet\r$ for some $\cs'\in\La$ and $\r\in\F$. Then
\[
\de(\beta_0)=\L(\cs)^+\cs^\f=(\cs'\bullet\L(\r))^+(\cs'\bullet\r)^\f=\Phi_{\cs'}(\L(\r)^+{\r^\f})=\Phi_{\cs'}(\de(\hat \beta_0)),
\]
where $\hat\beta_0=\Psi_{\cs'}^{-1}(\beta_0)\in E$. This implies that
\begin{equation*} 
\Phi_{\cs'}(0\de_2(\hat\beta_0)\de_3(\hat\beta_0)\ldots)=\Phi_{\cs'}({\r^\f})=(\cs'\bullet\r)^\f=\cs^\f.
\end{equation*}
{The same argument as in Case III then yields} (\ref{eq:right-discontinuity}).
\end{proof}

\begin{proof}[Proof of Theorem \ref{main:critical-devils-staircase}]
The theorem follows by Proposition \ref{prop:discontinuities}, Theorem \ref{main:critical-basic-intervals} and Theorem \ref{main:geometrical-structure-basic-intervals}.
\end{proof}

\section*{Acknowledgements}
The first author  is partially supported by Simons Foundation grant \# 709869. The second author thanks Charlene Kalle for some useful discussions when she visited Chongqing University in November {2019}. He is supported by NSFC No.~11971079 and the Fundamental
and Frontier Research Project of Chongqing No.~cstc2019jcyj-msxmX0338 and No.~cx2019067.


\end{document}